\documentclass{amsart}
\usepackage{amssymb}

\usepackage{graphicx}
\usepackage{amscd}

\newtheorem{theorem}{Theorem}
\theoremstyle{plain}

\newtheorem{corollary}{Corollary}

\newtheorem{definition}{Definition}

\newtheorem{lemma}{Lemma}

\newtheorem{proposition}{Proposition}

\numberwithin{equation}{section}

\input{tcilatex}
\input{diagrams.tex}
\newarrow{Onto}{-}{-}{-}{-}{->>}
\newarrow{Into}{C}{-}{-}{-}{>}
\newarrow{Corr}{<}{-}{-}{-}{>}
\newarrow{Ev}{|}{-}{-}{-}{>}
\begin{document}
\title{\emph{Twisting} of Quantum Spaces and Twisted coHom Objects}
\author{Sergio D. Grillo}
\address{\textit{Centro At\'{o}mico Bariloche and Instituto Balseiro}\\
\textit{\ 8400-S. C. de Bariloche}\\
\textit{\ Argentina}}
\email{sergio@cabtep2.cnea.gov.ar}
\date{May 2002}

\begin{abstract}
Twisting process for quantum linear spaces is defined. It consists in a
particular kind of globally defined deformations on finitely generated
algebras. Given a quantum space $\mathcal{A}=\left( \mathbf{A}_{1},\mathbf{A}%
\right) $, a multiplicative cosimplicial quasicomplex $\mathsf{C}^{\bullet }%
\left[ \mathbf{A}_{1}\right] $ in the category \textrm{Grp} is associated to 
$\mathbf{A}_{1}$, in such a way that for every $n\in \mathbb{N}$ a subclass
of linear automorphisms of $\mathbf{A}^{\otimes n}$ is obtained from the
groups $\mathsf{C}^{n}\left[ \mathbf{A}_{1}\right] $. Among the elements of
this subclass, the counital 2-cocycles are those which define the twist
transformations. In these terms, the twisted internal coHom objects $%
\underline{hom}^{\Upsilon }\left[ \mathcal{B},\mathcal{A}\right] $,
constructed in \cite{gm}, can be described as twisting of the proper coHom
objects $\underline{hom}\left[ \mathcal{B},\mathcal{A}\right] $. Moreover,
twisted tensor products $\mathcal{A}\circ _{\tau }\mathcal{B}$, in terms of
which above objects were built up, can be seen as particular 2-cocycle
twisting of $\mathcal{A}\circ \mathcal{B}$, enabling us to generalize the
mentioned construction. The quasicomplexes $\mathsf{C}^{\bullet }\left[ 
\mathbf{V}\right] $ are further studied, showing for instance that, when $%
\mathbf{V}$ is a coalgebra the quasicomplexes related to Drinfeld twisting,
corresponding to bialgebras generated by $\mathbf{V}$, are subobjects of $%
\mathsf{C}^{\bullet }\left[ \mathbf{V}\right] $.
\end{abstract}

\maketitle

\section*{Introduction}

In a previous paper \cite{gm}, combining the ideas of twisted tensor products%
\textbf{\ (}TTP) \cite{cap} and internal coHom objects, we have built up on
the monoidal category $\left( \mathrm{CA},\circ ,\mathcal{K}\right) $ of 
\emph{conic algebras} or \emph{conic }quantum spaces, i.e. finitely
generated graded algebras (which constitutes a full subcategory of the
-general- quantum linear spaces \cite{man0}\cite{man1}), functions $\left( 
\mathcal{B},\mathcal{A}\right) \mapsto \underline{hom}^{\Upsilon }\left[ 
\mathcal{B},\mathcal{A}\right] \in \mathrm{CA}$, one for each collection of
automorphisms $\Upsilon =\left\{ \sigma _{\mathcal{A}}:\mathcal{A}\backsimeq 
\mathcal{A}\right\} _{\mathcal{A}\in \mathrm{CA}}$, defining $\mathrm{CA}%
^{op}$-based categories distinct from the one related to its proper internal
coHom objects $\underline{hom}\left[ \mathcal{B},\mathcal{A}\right] =%
\mathcal{B}\triangleright \mathcal{A}$ (in the sense that, in general, $%
\underline{hom}\left[ \mathcal{B},\mathcal{A}\right] \ncong \underline{hom}%
^{\Upsilon }\left[ \mathcal{B},\mathcal{A}\right] $). The opposite objects 
\begin{equation*}
\underline{hom}^{\Upsilon }\left[ \mathcal{B},\mathcal{A}\right] ^{op}\doteq 
\underline{Hom}^{\Upsilon }\left[ \mathcal{B}^{op},\mathcal{A}^{op}\right]
\in \mathrm{CA}^{op}
\end{equation*}
represent `spaces of morphisms' whose coordinate rings, given precisely by
the algebras $\underline{hom}^{\Upsilon }\left[ \mathcal{B},\mathcal{A}%
\right] $, do not commute with the ones of their respective domain $\mathcal{%
B}^{op}$, in the sense that they give rise to \emph{twisted coevaluation}
arrows $\mathcal{A}\rightarrow \underline{hom}^{\Upsilon }\left[ \mathcal{B},%
\mathcal{A}\right] \circ _{\tau }\mathcal{B}$, where $\circ _{\tau }$ is a%
\textbf{\ }TTP constructed in terms of $\sigma _{\mathcal{A}},\sigma _{%
\mathcal{B}}\in \Upsilon $. This is why we have called them \emph{twisted}
coHom objects.

We have seen in that paper the objects $\underline{end}^{\Upsilon }\left[ 
\mathcal{A}\right] =\underline{hom}^{\Upsilon }\left[ \mathcal{A},\mathcal{A}%
\right] $, which are endowed with a bialgebra structure, are counital
2-cocycle twisting \cite{drin} of the coEnd objects $\underline{end}\left[ 
\mathcal{A}\right] =\mathcal{A}\triangleright \mathcal{A}$. In other words,
those objects are related by a (non formal) bialgebra \emph{deformation}.
However, we have not been able to make an analogous claim relating the
objects $\underline{hom}^{\Upsilon }\left[ \mathcal{B},\mathcal{A}\right] $
and $\underline{hom}\left[ \mathcal{B},\mathcal{A}\right] $ for $\mathcal{B}%
\neq \mathcal{A}$ yet. To this end, we will develope in this paper a
suitable kind of globally defined deformations of quantum spaces, which we
shall call \emph{twisting of quantum spaces}, in such a way that a
multiplicative cosimplicial quasicomplex structure $\mathsf{C}^{\bullet }$
can be related to them, being 2-cocycles the \emph{well-behaved}
deformations. We shall see the objects $\underline{hom}^{\Upsilon }\left[ 
\mathcal{B},\mathcal{A}\right] $ are twisting by 2-cocycles of $\underline{%
hom}\left[ \mathcal{B},\mathcal{A}\right] $. In particular, for the coEnd
objects, those transformations define the previously cited 2-cocycle
bialgebra twisting.

The main aim of our work is to construct and analyze the mentioned
cosimplicial quasicomplex structure and related algebra deformation process.
We invoke the concept of twisted coHom objects just because they were our
main motivation to this paper.

\bigskip

Although our twist transformation will be defined on the hole class of
quantum spaces, i.e. the category $\mathrm{FGA}$ of finitely generated
algebras, they are mainly designed to be applied on the conic ones. This is
why, among other things, the first chapter is dedicated to them. Thus, in \S 
\textbf{1} we make a brief review about conic quantum spaces, recalling
their monoidal structures and related internal (co)Hom objects. For a more
extended treatment, reference \cite{gm} can be consulted.

In order to introduce the twisting process in a (quasi)cohomological
framework, in \S \textbf{2} we endow\textbf{\ }every linear space $\mathbf{V}
$ (or every tensor algebra $\mathbf{V}^{\otimes }$) with a multiplicative
cosimplicial quasicomplex structure $\left( \mathsf{C}^{\bullet }\left[ 
\mathbf{V}\right] ,\partial \right) $. Recall that formal deformations \cite
{ger} of an associative algebra $\mathbf{A}$ are controlled by the
Hochschild complex \cite{ho} of $\mathbf{A}$. Given a pair $\mathcal{A}%
=\left( \mathbf{A}_{1},\mathbf{A}\right) $ and a (counital) 2-cocycle $\psi $
inside certain subclass of $\mathsf{C}^{2}\left[ \mathbf{A}_{1}\right] $,
the \emph{admissible} 2-cochains, we define in \S \textbf{3} a deformation
of $\mathcal{A}$ as a new quantum space $\mathcal{A}_{\psi }$, the \emph{%
twisting} of $\mathcal{A}$ by $\psi $. The mentioned subclass defines a
subgroup of linear automorphisms of $\mathbf{A}\otimes \mathbf{A}$. Deformed
quantum spaces $\mathcal{A}_{\psi }$ and $\mathcal{A}_{\varphi }$ are
isomorphic \emph{iff }the cocycles $\psi $ and $\varphi $ are cohomologous
through some admissible 1-cochain. We show that twist transformation can be
composed and have inverse. Moreover, an equivalence relation between quantum
spaces can be defined, namely: $\mathcal{A}$ and $\mathcal{B}$ are \emph{%
twist or gauge related} when there exists $\psi $ such that $\mathcal{A}%
_{\psi }\backsimeq \mathcal{B}$ (as in the case of bialgebra twisting). We
also show the twisting of bialgebras in the category $\mathrm{FGA}$ are
particular twisting of quantum spaces.

Finally, we investigate in \S \textbf{4} the behavior of twist
transformations in relation to certain functorial structures in $\mathrm{FGA}
$. We see, for example, that products $\mathcal{A}_{\psi }\circ \mathcal{B}%
_{\varphi }$ are twisting of products $\mathcal{A}\circ \mathcal{B}$.
Analogous claim is also true, when restricted to particular kind of twist
transformations, for functors $!$, $\bullet $, and products $\odot $, $%
\triangleleft $, $\triangleright $ and $\diamond $ defined in \cite{gm}.
This study lead us to a better understanding of the construction of twisted
coHom objects and, in turn, a way to generalize it.

\bigskip

We often adopt definitions and notation extracted form Mac Lane's book \cite
{mac}. $\Bbbk $ indicates some of the numeric fields, $\mathbb{R}$ or $%
\mathbb{C}$. The usual tensor product on $\Bbbk \mathrm{-Alg}\equiv \mathrm{%
Alg}$\textrm{\ }and $\mathrm{Vct}_{\Bbbk }\equiv \mathrm{Vct}$\textrm{\ }%
(the categories of unital associative $\Bbbk $-algebras and of $\Bbbk $%
-vector spaces, respectively) is denoted by $\otimes $. \textrm{Grp} denotes
the category of groups and their homomorphisms.

\section{Conic quantum spaces}

In this chapter we recall definition, functorial structures and internal
coHom objects (standard and twisted ones) of a subclass of quantum spaces
that we have called \emph{conic quantum spaces} \cite{gm}, introducing the
needed notation for the rest of the paper. We suppose the reader is familiar
with the concepts of quadratic and general quantum spaces \cite{man0}\cite
{man1}.

Conic quantum spaces define a category $\mathrm{CA}$ such that the full
inclusions $\mathrm{QA}\subset \mathrm{CA}\subset \mathrm{FGA}$ holds, being 
$\mathrm{FGA}$ and $\mathrm{QA}$\textrm{\ }the categories of finitely
generated algebras and of quadratic algebras, respectively. Remember that 
\textrm{FGA }is formed out by pairs $\left( \mathbf{A}_{1},\mathbf{A}\right) 
$ where $\mathbf{A}$ is a unital algebra generated by a finite linear
subspace $\mathbf{A}_{1}\subset \mathbf{A}$, and its morphisms are algebra
morphisms $\mathbf{A}\rightarrow \mathbf{B}$ such that restricted to $%
\mathbf{A}_{1}\subset \mathbf{A}$ give linear maps $\mathbf{A}%
_{1}\rightarrow \mathbf{B}_{1}$ ($\subset \mathbf{B}$). The objects of $%
\mathrm{CA}$ are pairs $\mathcal{A}\doteq \left( \mathbf{A}_{1},\mathbf{A}%
\right) \in \mathrm{FGA}$ such that $\mathbf{A}$ is a graded algebra with
grading given by 
\begin{equation*}
\mathbf{A}=\bigoplus\nolimits_{n\in \mathbb{N}_{0}}\,\mathbf{A}_{n};\;\;\;%
\mathbf{A}_{n}=\Pi \left( \mathbf{A}_{1}^{\otimes n}\right) ;\;\;\;\mathbb{N}%
_{0}=\mathbb{N}\cup \left\{ 0\right\} ,
\end{equation*}
where $\Pi =\Pi _{\mathbf{A}}:\mathbf{A}_{1}^{\otimes }\twoheadrightarrow 
\mathbf{A}$ is the canonical epimorphism of algebras such that restricted to 
$\mathbf{A}_{1}$ gives the inclusion $\mathbf{A}_{1}\hookrightarrow \mathbf{A%
}$. We are denoting by $\mathbf{A}_{1}^{\otimes }$ the tensor algebra of $%
\mathbf{A}_{1}$, i.e. $\mathbf{A}_{1}^{\otimes }=\bigoplus_{n\in \mathbb{N}%
_{0}}\mathbf{A}_{1}^{\otimes n}$. In the general case the related algebras
are just filtered, more precisely 
\begin{equation}
\mathbf{A}=\bigcup\nolimits_{n\in \mathbb{N}_{0}}\mathbf{F}_{n};\;\mathbf{F}%
_{n}=\Pi \left( \bigoplus\nolimits_{i=0}^{n}\mathbf{A}_{1}^{\otimes
i}\right) .  \label{filt}
\end{equation}
Alternatively, the objects of $\mathrm{CA}$ can be described as those pairs $%
\mathcal{A}$ in $\mathrm{FGA}$ such that $\ker \Pi $ is a graded bilateral
ideal $\ker \Pi _{\mathbf{A}}=\bigoplus\nolimits_{n\in \mathbb{N}_{0}}%
\mathbf{I}_{n}$ with $\mathbf{I}_{n}\subset \mathbf{A}_{1}^{\otimes n}$. In
these terms, given another pair $\left( \mathbf{B}_{1},\mathbf{B}\right) $
with $\ker \Pi _{\mathbf{B}}=\bigoplus_{n\in \mathbb{N}_{0}}\mathbf{J}_{n}$,
the arrows $\left( \mathbf{A}_{1},\mathbf{A}\right) \rightarrow \left( 
\mathbf{B}_{1},\mathbf{B}\right) $ in $\mathrm{CA}$ are characterized by
linear maps $\alpha _{1}:\mathbf{A}_{1}\rightarrow \mathbf{B}_{1}$ such that 
$\alpha _{1}^{\otimes n}\left( \mathbf{I}_{n}\right) \subset \mathbf{J}_{n}$%
, being 
\begin{equation*}
\alpha _{1}^{\otimes }=\bigoplus\nolimits_{n\in \mathbb{N}_{0}}\alpha
_{1}^{\otimes n}:\mathbf{A}_{1}^{\otimes }\rightarrow \mathbf{B}%
_{1}^{\otimes }
\end{equation*}
the unique extension of $\alpha _{1}$ to $\mathbf{A}_{1}^{\otimes }$ as a
morphism of algebras. The algebra map $\alpha :\mathbf{A}\rightarrow \mathbf{%
B}$ that defines any morphism in \textrm{FGA }is the unique algebra
homomorphism such that $\alpha \,\Pi _{\mathbf{A}}=\Pi _{\mathbf{B}}\,\alpha
_{1}^{\otimes }$. In particular, $\alpha _{1}=\left. \alpha \right| _{%
\mathbf{A}_{1}}$. Returning to the ideal of a conic quantum space, since the
restriction of $\Pi _{\mathbf{A}}$ to $\mathbf{A}_{1}$ is the identity, we
always have $\mathbf{I}_{1}=\left\{ 0\right\} $. In addition, $\mathbf{I}%
_{0}\neq \left\{ 0\right\} $ \emph{iff }$\mathbf{A}=\left\{ 0\right\} $.
Thus, generically $\ker \Pi _{\mathbf{A}}=\bigoplus\nolimits_{n\geq 2}%
\mathbf{I}_{n}$, unless $\mathbf{A}=\left\{ 0\right\} $.

\bigskip

Examples of conic quantum spaces are, beside quadratics, those with $\ker
\Pi =I\left[ \mathbf{Y}\right] $, where $\mathbf{Y}\subset \mathbf{A}%
_{1}^{\otimes m}$ for some $m>2$, indicating by $I\left[ \mathbf{X}\right] $
the bilateral ideal generated algebraically by $\mathbf{X}\subset \mathbf{A}%
_{1}^{\otimes }$; i.e. $I\left[ \mathbf{X}\right] =\mathbf{A}_{1}^{\otimes
}\otimes \mathbf{X}\otimes \mathbf{A}_{1}^{\otimes }$. We have called them $%
m $\textbf{-}\emph{th quantum spaces}, and denoted $\mathrm{CA}^{m}$ the
full subcategory of $\mathrm{FGA}$ which has these pairs as objects. Thus $%
\mathrm{QA}=\mathrm{CA}^{2}$.

\bigskip

We recall that from the geometric point of view, the quantum spaces actually
are the opposite objects $\mathcal{A}^{op}$ to the pairs $\mathcal{A}=\left( 
\mathbf{A}_{1},\mathbf{A}\right) \in \mathrm{FGA}$, being the latter their
(generically) \emph{noncommutative coordinate rings}. That is to say $%
\mathrm{FGA}^{op}$, which we indicate $\mathsf{QLS}$, is the category of
quantum spaces, and $\mathrm{FGA}$\textrm{\ }the category of their
coordinate rings. Nevertheless, since the duality between $\mathsf{QLS}$%
\textrm{\ }and $\mathrm{FGA}$, for quantum spaces\emph{\ }we understand the
objects of each one of these categories. Categories $\mathrm{CA}$ and $%
\mathrm{CA}^{m}$ give rise to full subcategories of $\mathsf{QLS}$.

\subsection{Functorial structures and internal coHom objects}

The bifunctor $\circ $ defines an (strict and symmetric) monoidal structure
with unit $\mathcal{K}=\left( \Bbbk ,\Bbbk ^{\otimes }\right) $ in $\mathrm{%
CA}$. Recall that $\mathcal{A}\circ \mathcal{B}=\left( \mathbf{A}_{1}\otimes 
\mathbf{B}_{1},\mathbf{A}\circ \mathbf{B}\right) $, being $\mathbf{A}\circ 
\mathbf{B}$ the subalgebra of $\mathbf{A}\otimes \mathbf{B}$ generated by $%
\mathbf{A}_{1}\otimes \mathbf{B}_{1}$, and that the unit for $\circ $ is $%
\mathcal{I}\doteq \left( \Bbbk ,\Bbbk \right) $ in $\mathrm{FGA}$. If $%
\mathcal{A},\mathcal{B}\in \mathrm{CA}$, the algebra $\mathbf{A}\circ 
\mathbf{B}$ effectively defines an object of $\mathrm{CA}$ due to the kernel
of $\Pi _{\circ }:\left[ \mathbf{A}_{1}\otimes \mathbf{B}_{1}\right]
^{\otimes }\twoheadrightarrow \mathbf{A}\circ \mathbf{B}$, using the
(canonical) homogeneous isomorphism of graded algebras 
\begin{equation}
\left[ \mathbf{A}_{1}\otimes \mathbf{B}_{1}\right] ^{\otimes }\backsimeq %
\left[ \mathbf{A}_{1}\otimes \mathbf{B}_{1}\right] ^{\widehat{\otimes }%
}\doteq \,\bigoplus\nolimits_{n\in \mathbb{N}_{0}}\left( \mathbf{A}%
_{1}^{\otimes n}\otimes \mathbf{B}_{1}^{\otimes n}\right) ,  \label{canis}
\end{equation}
being the latter the subalgebra of $\mathbf{A}_{1}^{\otimes }\otimes \mathbf{%
B}_{1}^{\otimes }$ generated by $\mathbf{A}_{1}\otimes \mathbf{B}_{1}$, is a
graded bilateral ideal isomorphic to 
\begin{equation*}
\bigoplus\nolimits_{n\in \mathbb{N}_{0}}\left( \mathbf{I}_{n}\otimes \mathbf{%
B}_{1}^{\otimes n}+\mathbf{A}_{1}^{\otimes n}\otimes \mathbf{J}_{n}\right)
\;\left( \subset \left[ \mathbf{A}_{1}\otimes \mathbf{B}_{1}\right] ^{%
\widehat{\otimes }}\right) .
\end{equation*}
Note that $\left[ \mathbf{A}_{1}\otimes \mathbf{B}_{1}\right] ^{\widehat{%
\otimes }}=\mathbf{A}_{1}^{\otimes }\circ \mathbf{B}_{1}^{\otimes }$. We
frequently identify the latter algebra and $\left[ \mathbf{A}_{1}\otimes 
\mathbf{B}_{1}\right] ^{\otimes }$. In particular, $\Pi _{\circ }$ will also
be understood as a map $\mathbf{A}_{1}^{\otimes }\circ \mathbf{B}%
_{1}^{\otimes }\twoheadrightarrow \mathbf{A}\circ \mathbf{B}$.

\bigskip

For quadratic algebras there is another monoidal structure $\bullet $ and a
functor $!$ such that 
\begin{equation}
\mathcal{A}^{!!}\backsimeq \mathcal{A},\;\;\left( \mathcal{A}\circ \mathcal{B%
}\right) ^{!}\backsimeq \mathcal{A}^{!}\bullet \mathcal{B}^{!},\;\;\left( 
\mathcal{A}\bullet \mathcal{B}\right) ^{!}\backsimeq \mathcal{A}^{!}\circ 
\mathcal{B}^{!},\;\;\mathcal{K}^{!}\backsimeq \mathcal{U},  \label{p}
\end{equation}
being $\mathcal{U}=\left( \Bbbk ,\mathbf{U}\right) $, $\mathbf{U}=\left.
\Bbbk \left[ e\right] \right/ I\left[ e^{2}\right] $, a unit object for $%
\bullet $. Analogous functors to $\bullet $ and $!$ can be given in $\mathrm{%
CA}$, namely 
\begin{equation*}
\mathbf{A}\odot \mathbf{B}\doteq \left. \mathbf{A}_{1}^{\otimes }\circ 
\mathbf{B}_{1}^{\otimes }\right/ \bigoplus\nolimits_{n\in \mathbb{N}_{0}}%
\mathbf{I}_{n}\otimes \mathbf{J}_{n}
\end{equation*}
and 
\begin{equation}
\mathbf{A}^{!}\doteq \left. \mathbf{A}_{1}^{\ast \otimes }\right/ \mathbf{I}%
^{\dagger };\;\;\;\mathbf{I}^{\dagger }\doteq \bigoplus\nolimits_{n\in 
\mathbb{N}_{0}}\mathbf{I}_{n}^{\dagger }\doteq I\left[ \bigoplus\nolimits_{n%
\geq 2}\mathbf{I}_{n}^{\perp }\right] ,  \label{idag}
\end{equation}
with $\mathbf{I}_{n}^{\perp }=\left\{ x\in \mathbf{A}_{1}^{\ast \otimes
n}:\left\langle x,y\right\rangle =0,\;\forall y\in \mathbf{I}_{n}\right\} $,
for\ $n\geq 2$. We shall take $\mathbf{I}_{0,1}^{\perp }=\left\{ 0\right\} $%
. A unit element for $\left( \mathrm{CA},\odot \right) \ $is $\mathcal{U}$
as for $\mathrm{QA}$, but the restriction of $\odot $ to $\mathrm{QA}$ does
not coincide with the functor $\bullet $. On the other hand, the functor $!$
does coincide with the corresponding to the quadratic cases, although $%
!^{2}\ncong id_{\mathrm{CA}}$. The preserved properties are, with respect to 
$\left( \ref{p}\right) $, 
\begin{equation*}
\left( \mathcal{A}\odot \mathcal{B}\right) ^{!}\backsimeq \mathcal{A}%
^{!}\circ \mathcal{B}^{!},\qquad \mathcal{K}^{!}\backsimeq \mathcal{U}%
,\qquad \;\mathcal{K}\backsimeq \mathcal{U}^{!},
\end{equation*}
while $\left( \mathcal{A}\circ \mathcal{B}\right) ^{!}\ncong \mathcal{A}%
^{!}\odot \mathcal{B}^{!}$ and $\mathcal{A}^{!!}\ncong \mathcal{A}$. Some
kind of covariant mixing of $\odot $ and $!$ can be defined, namely 
\begin{equation*}
\begin{array}{c}
\triangleright :\mathrm{CA}^{op}\times \mathrm{CA}\rightarrow \mathrm{CA}%
,\;\;\triangleleft :\mathrm{CA}\times \mathrm{CA}^{op}\rightarrow \mathrm{CA}%
, \\ 
\\ 
\diamond :\mathrm{CA}^{op}\times \mathrm{CA}^{op}\rightarrow \mathrm{CA},
\end{array}
\end{equation*}
satisfying 
\begin{equation}
\triangleright \,=\bullet \,\left( !\times id\right) ,\;\triangleleft
\,=\bullet \,\left( id\times !\right) ,\;\diamond \,=\bullet \,\left(
!\times !\right)  \label{compa}
\end{equation}
when restricted to $\mathrm{QA}$. They are given on objects by (identifying
each object with its opposite) 
\begin{eqnarray*}
\mathbf{A}\triangleright \mathbf{B} &\doteq &\left. \mathbf{A}_{1}^{\ast
\otimes }\circ \mathbf{B}_{1}^{\otimes }\right/ I\left[ \bigoplus%
\nolimits_{n\in \mathbb{N}_{0}}\mathbf{I}_{n}^{\perp }\otimes \mathbf{J}_{n}%
\right] , \\
&& \\
\mathbf{A}\triangleleft \mathbf{B} &\doteq &\left. \mathbf{A}_{1}^{\otimes
}\circ \mathbf{B}_{1}^{\ast \otimes }\right/ I\left[ \bigoplus\nolimits_{n%
\in \mathbb{N}_{0}}\mathbf{I}_{n}\otimes \mathbf{J}_{n}^{\perp }\right] , \\
&& \\
\mathbf{A}\diamond \mathbf{B} &\doteq &\left. \mathbf{A}_{1}^{\ast \otimes
}\circ \mathbf{B}_{1}^{\ast \otimes }\right/ I\left[ \bigoplus\nolimits_{n%
\in \mathbb{N}_{0}}\mathbf{I}_{n}^{\perp }\otimes \mathbf{J}_{n}^{\perp }%
\right] .
\end{eqnarray*}
The functors $\triangleright $ and $\triangleleft $ have $\mathcal{K}$ as
left and as a right unit, respectively, in the sense that $\mathcal{K}%
\triangleright \mathcal{A}\backsimeq \mathcal{A}$ and $\mathcal{A}%
\triangleleft \mathcal{K}\backsimeq \mathcal{A}$ for any $\mathcal{A}$ in $%
\mathrm{CA}$ (on the other hand, $\mathcal{A}\triangleright \mathcal{U}%
\backsimeq \mathcal{A}^{!}$, $\mathcal{U}\triangleleft \mathcal{A}\backsimeq 
\mathcal{A}^{!}$ and $\mathcal{K}\diamond \mathcal{A}\backsimeq \mathcal{A}%
\diamond \mathcal{K}\backsimeq \mathcal{A}^{!}$).

While the internal coHom objects of $\mathrm{QA}$ are given by $\underline{%
hom}\left[ \mathcal{B},\mathcal{A}\right] =\mathcal{B}^{!}\bullet \mathcal{A}
$, for $\mathrm{CA}$, as it might be expected from Eq. $\left( \ref{compa}%
\right) $, we have the following result.

\begin{theorem}
The category $\left( \mathrm{CA},\circ ,\mathcal{K}\right) $ has \textbf{%
internal coHom objects} given by $\underline{hom}\left[ \mathcal{B},\mathcal{%
A}\right] =\mathcal{B}\triangleright \mathcal{A}$ with coevaluation arrow $%
\mathcal{A}\rightarrow \left( \mathcal{B}\triangleright \mathcal{A}\right)
\circ \mathcal{B}$ defined by the map $a_{i}\mapsto \left( b^{j}\otimes
a_{i}\right) \otimes b_{j}$ \emph{(}sum over repeated indices is understood%
\emph{)}, where $\left\{ a_{i}\right\} $, $\left\{ b_{j}\right\} $ and $%
\left\{ b^{j}\right\} $ are basis of $\mathbf{A}_{1}$, $\mathbf{B}_{1}$ and $%
\mathbf{B}_{1}^{\ast }$, respectively.\ \ \ $\blacksquare $
\end{theorem}

The proof is given in \cite{gm}.

\subsection{The twisted internal coHom objects}

An internal coHom object in $\mathrm{FGA}$ is an initial object of the comma
category $\left( \mathcal{A}\downarrow \mathrm{FGA}\circ \mathcal{B}\right) $%
. The objects of each $\left( \mathcal{A}\downarrow \mathrm{FGA}\circ 
\mathcal{B}\right) $, \emph{diagrams }in the Manin terminology, are pairs $%
\left\langle \varphi ,\mathcal{H}\right\rangle _{_{\mathcal{A}\mathbf{,}%
\mathcal{B}}}=\left\langle \varphi ,\mathcal{H}\right\rangle $, with $%
\mathcal{H}\in \mathrm{FGA}$ and $\varphi $ a morphism $\mathcal{A}%
\rightarrow $ $\mathcal{H}\circ \mathcal{B}$; and its arrows $\left\langle
\varphi ,\mathcal{H}\right\rangle \rightarrow \left\langle \varphi ^{\prime
},\mathcal{H}^{\prime }\right\rangle $ are given by morphisms $\alpha :%
\mathcal{H}\rightarrow \mathcal{H}^{\prime }$ satisfying $\varphi ^{\prime
}=\left( \alpha \circ I_{B}\right) \,\varphi $. For every such category $%
\left\langle \varphi ,\mathcal{H}\right\rangle \mapsto \mathcal{H}$ defines
an embedding $\mathsf{U}:\left( \mathcal{A}\downarrow \mathrm{FGA}\circ 
\mathcal{B}\right) \hookrightarrow \mathrm{FGA}$. The disjoint union of the
family $\left\{ \left( \mathcal{A}\downarrow \mathrm{FGA}\circ \mathcal{B}%
\right) \right\} _{\mathcal{A},\mathcal{B}\in \mathrm{FGA}}$, namely $%
\mathrm{FGA}^{\circ }$, has a semigroupoid structure given by the partial
product functor $\widehat{\circ }$. $\mathsf{U}$ extends to an obvious
embedding $\mathrm{FGA}^{\circ }\hookrightarrow \mathrm{FGA}$, which we
shall also call $\mathsf{U}$, and that satisfies 
\begin{equation}
\mathsf{U}\,\widehat{\circ }=\circ \,\left( \mathsf{U}\times \mathsf{U}%
\right) \;\;\;and\;\;\;\mathsf{U}\left\langle \ell _{\mathcal{A}},\mathcal{I}%
\right\rangle =\mathcal{I}.  \label{morr}
\end{equation}
So $\mathsf{U}:\mathrm{FGA}^{\circ }\hookrightarrow \mathrm{FGA}$ is an
embedding of categories with unital associative partial products. The same
is true for $\mathrm{CA}$ and every $\mathrm{CA}^{m}$, but changing $%
\mathcal{I}$ by $\mathcal{K}$ in Eq. $\left( \ref{morr}\right) $. We have
seen in \cite{gm} that it is possible to arrive at the notions of
cocomposition and coidentity from the semigroupoid structure in $\mathrm{FGA}%
^{\circ }$, the initiality of each $\underline{hom}\left[ \mathcal{B},%
\mathcal{A}\right] \in \left( \mathcal{A}\downarrow \mathrm{FGA}\circ 
\mathcal{B}\right) $, and the existence of a functor $\mathsf{U}$ satisfying
Eq. $\left( \ref{morr}\right) $. Then, the function $\left( \mathcal{B},%
\mathcal{A}\right) \mapsto \underline{hom}\left[ \mathcal{B},\mathcal{A}%
\right] $, together with the corresponding cocomposition and coidentity
arrows, defines an $\mathrm{FGA}$-cobased (or dually $\mathrm{FGA}^{op}$%
-based) category with objects also in $\mathrm{FGA}$. Such a cobased
category has in addition the usual notion of coevaluation $\mathcal{A}%
\rightarrow \underline{hom}\left[ \mathcal{B},\mathcal{A}\right] \circ 
\mathcal{B}$.

Based on these ideas, we construct in that paper a family of categories $%
\Upsilon ^{\mathcal{A},\mathcal{B}}$ with initial objects, denoted $%
\underline{hom}^{\Upsilon }\left[ \mathcal{B},\mathcal{A}\right] $, and with
a family of embedding $\Upsilon ^{\mathcal{A},\mathcal{B}}\hookrightarrow 
\mathrm{FGA}$. The objects of each $\Upsilon ^{\mathcal{A},\mathcal{B}}$ are
essentially arrows $\mathcal{A}\rightarrow \mathcal{H}\circ _{\tau }\mathcal{%
B}$, where $\circ _{\tau }$ denotes a twisted tensor product between $%
\mathcal{H}$ and $\mathcal{B}$. Each category is defined by a linear
isomorphism $\widehat{\tau }_{\mathcal{A},\mathcal{B}}:\mathbf{B}_{1}\otimes 
\mathbf{B}_{1}^{\ast }\otimes \mathbf{A}_{1}\backsimeq \mathbf{B}_{1}\otimes 
\mathbf{B}_{1}^{\ast }\otimes \mathbf{A}_{1}$ from which the twisted
products $\circ _{\tau }$ are built up. Restricting $\mathcal{A},\mathcal{B}$
to $\mathrm{CA}$ and considering $\widehat{\tau }_{\mathcal{A},\mathcal{B}%
}=id\otimes \sigma _{\mathcal{B}}^{\ast -1}\otimes \sigma _{\mathcal{A}}$,
the disjoint union $\Upsilon ^{\cdot }$ of these categories has a
semigroupoid structure such that an equation like $\left( \ref{morr}\right) $
holds for the extended functor $\Upsilon ^{\cdot }\hookrightarrow \mathrm{CA}
$. Thus, $\left( \mathcal{B},\mathcal{A}\right) \mapsto \underline{hom}%
^{\Upsilon }\left[ \mathcal{B},\mathcal{A}\right] $ gives rise to an $%
\mathrm{CA}$-cobased category with the additional notion of coevaluation $%
\mathcal{A}\rightarrow \underline{hom}^{\Upsilon }\left[ \mathcal{B},%
\mathcal{A}\right] \circ _{\tau }\mathcal{B}$. We will not give the details
of this construction, referring the interested reader to \cite{gm}.
Nevertheless, in the rest of this section we review the main properties of
the objects $\underline{hom}^{\Upsilon }\left[ \mathcal{B},\mathcal{A}\right]
$ and its relationship to the proper coHom objects $\mathcal{B}%
\triangleright \mathcal{A}$ of $\mathrm{CA}$.

\bigskip

Let $\left\{ \sigma _{\mathcal{A}}:\mathbf{A}_{1}\backsimeq \mathbf{A}%
_{1}\right\} _{\mathcal{A}\in \mathrm{CA}}$ be a collection defining $%
\Upsilon ^{\cdot }$, such that each map $\sigma _{\mathcal{A}}$ can be
extended to an automorphism $\mathcal{A}\backsimeq \mathcal{A}$. Consider a
couple $\mathcal{A}$ and $\mathcal{B}$ of quantum spaces given by algebras $%
\mathbf{A}\backsimeq \left. \mathbf{A}_{1}^{\otimes }\right/ \mathbf{I}$ and 
$\mathbf{B}\backsimeq \left. \mathbf{B}_{1}^{\otimes }\right/ \mathbf{J}$,
being $\mathbf{I}$ and $\mathbf{J}$ the ideals linearly generated by 
\begin{equation}
R_{\lambda _{n}}^{k_{1}...k_{n}}\;a_{k_{1}}...a_{k_{n}}\in \mathbf{I}%
_{n},\;\;\;\;S_{\mu _{n}}^{k_{1}...k_{n}}\;b_{k_{1}}...b_{k_{n}}\in \mathbf{J%
}_{n},  \label{I}
\end{equation}
with $\lambda _{n}$ and $\mu _{n}$ in index sets $\Lambda _{n}$ and $\Phi
_{n}$, respectively. Writing $\phi \doteq \sigma _{\mathcal{A}}$ and $\rho
\doteq \sigma _{\mathcal{B}}^{\ast -1}$, the initial objects of each
category $\Upsilon ^{\mathcal{A},\mathcal{B}}$ can be defined by a conic
quantum space 
\begin{equation*}
\underline{hom}^{\Upsilon }\left[ \mathcal{B},\mathcal{A}\right] =\mathcal{B}%
^{\Upsilon }\triangleright \mathcal{A}^{\Upsilon },\;\;\mathcal{A}^{\Upsilon
}\doteq \left( \mathbf{A}_{1},\left. \mathbf{A}_{1}^{\otimes }\right/ 
\mathbf{I}_{\sigma }\right) ,\;\;\mathcal{B}^{\Upsilon }\doteq \left( 
\mathbf{B}_{1},\left. \mathbf{B}_{1}^{\otimes }\right/ \mathbf{J}_{\sigma
}\right) ,
\end{equation*}
with $\mathbf{I}_{\sigma }$ and $\mathbf{J}_{\sigma }$ linearly generated by 
\begin{equation*}
\left\{ \left\{ _{_{\sigma }}R_{\lambda
_{n}}^{k_{1}...k_{n}}\;a_{k_{1}}...a_{k_{n}}\right\} _{\lambda _{n}\in
\Lambda _{n}}\right\} _{n\in \mathbb{N}_{0}},\;\;\left\{ \left\{ _{_{\sigma
}}S_{\mu _{n}}^{k_{1}...k_{n}}\;b_{k_{1}}...b_{k_{n}}\right\} _{\mu _{n}\in
\Phi _{n}}\right\} _{n\in \mathbb{N}_{0}},
\end{equation*}
where 
\begin{eqnarray}
_{_{\sigma }}R_{\lambda _{n}}^{k_{1}...k_{n}} &\doteq &R_{\lambda
_{n}}^{k_{1}j_{2}...j_{n}}\;\phi _{j_{2}}^{k_{2}}\left( \phi ^{2}\right)
_{j_{3}}^{k_{3}}...\left( \phi ^{n-1}\right) _{j_{n}}^{k_{n}},  \label{II} \\
&&  \notag \\
\;\;_{_{\sigma }}S_{\mu _{n}}^{k_{1}...k_{n}} &\doteq &S_{\mu
_{n}}^{k_{1}j_{2}...j_{n}}\;\left( \rho ^{-1}\right) _{j_{2}}^{k_{2}}\left(
\rho ^{-2}\right) _{j_{3}}^{k_{3}}...\left( \rho ^{1-n}\right)
_{j_{n}}^{k_{n}}.  \label{III}
\end{eqnarray}
On the other hand, if each $\mathbf{J}_{n}^{\perp }\subset \mathbf{B}%
_{1}^{\ast \otimes n}$ is the span of 
\begin{equation}
\left\{ b^{k_{1}}...b^{k_{n}}\;\left( S^{\bot }\right)
_{k_{1}...k_{n}}^{\omega _{n}}\right\} _{\omega _{n}\in \Omega _{n}},
\label{IV}
\end{equation}
$\left( \mathbf{J}_{\sigma }\right) _{n}^{\perp }$ is spanned by the set 
\begin{equation}
\left\{ b^{k_{1}}b^{k_{2}}...b^{k_{n}}\;\left( _{_{\sigma }}S^{\bot }\right)
_{k_{1}...k_{n}}^{\omega _{n}}\right\} _{\omega _{n}\in \Omega _{n}},
\label{V}
\end{equation}
with 
\begin{equation}
\left( _{_{\sigma }}S^{\bot }\right) _{k_{1}...k_{n}}^{\omega _{n}}\doteq
\rho _{k_{2}}^{j_{2}}\left( \rho ^{2}\right) _{k_{3}}^{j_{3}}...\left( \rho
^{n-1}\right) _{k_{n}}^{j_{n}}\;\left( S^{\bot }\right)
_{k_{1}j_{2}...j_{n}}^{\omega _{n}}.  \label{VI}
\end{equation}
Therefore, $\underline{hom}^{\Upsilon }\left[ \mathcal{B},\mathcal{A}\right] 
$ will be the algebra generated by $z_{i}^{j}=b^{j}\otimes a_{i}$ and
quotient by the ideal algebraically generated by 
\begin{equation*}
\left\{ \left\{ _{\sigma }R_{\lambda
_{n}}^{k_{1}...k_{n}}\;z_{k_{1}}^{j_{1}}\cdot
z_{k_{2}}^{j_{2}}\;...\;z_{k_{n}}^{j_{n}}\;\left( _{\sigma }S^{\bot }\right)
_{j_{1}...j_{n}}^{\omega _{n}}\right\} _{\omega _{n}\in \Omega
_{n}}^{\lambda _{n}\in \Lambda _{n}}\right\} _{n\in \mathbb{N}_{0}}.
\end{equation*}

From the semigroupoid structure of $\Upsilon ^{\cdot }$ and the embedding $%
\Upsilon ^{\mathcal{A},\mathcal{B}}\hookrightarrow \mathrm{CA}$, there exist
arrows 
\begin{equation*}
\underline{hom}^{\Upsilon }\left[ \mathcal{B},\mathcal{A}\right] \rightarrow 
\underline{hom}^{\Upsilon }\left[ \mathcal{C},\mathcal{A}\right] \circ 
\underline{hom}^{\Upsilon }\left[ \mathcal{B},\mathcal{C}\right] ,\;\;\;%
\underline{end}^{\Upsilon }\left[ \mathcal{A}\right] \rightarrow \mathcal{K},
\end{equation*}
giving us the notions of cocomposition and coidentity we have just
mentioned. Of course, these arrows define a bialgebra structure on $%
\underline{end}^{\Upsilon }\left[ \mathcal{A}\right] $. Moreover, there
exists a counital 2-cocycle $\chi _{\mathcal{A}}=\chi $, such that $%
\underline{end}^{\Upsilon }\left[ \mathcal{A}\right] \backsimeq \left( 
\mathcal{A}\triangleright \mathcal{A}\right) _{\chi }$. That means $%
\underline{end}^{\Upsilon }\left[ \mathcal{A}\right] $ is a twisting by $%
\chi $ of $\underline{end}\left[ \mathcal{A}\right] $ \cite{drin}. Moreover,
from the actions $\mathcal{A}\rightarrow \underline{end}\left[ \mathcal{A}%
\right] \circ \mathcal{A}$ (the standard coevaluation maps for $\mathrm{CA}$%
), the above twisting can be translated to each $\mathcal{A}$ defining a new
quantum space $\mathcal{A}_{\chi }$ isomorphic to $\mathcal{A}^{\Upsilon }$
(c.f. \cite{maj}, page 54). Finally, 
\begin{equation*}
\underline{end}^{\Upsilon }\left[ \mathcal{A}\right] \backsimeq \mathcal{A}%
_{\chi }\triangleright \mathcal{A}_{\chi }=\underline{end}\left[ \mathcal{A}%
_{\chi }\right] \backsimeq \underline{end}\left[ \mathcal{A}\right] _{\chi }.
\end{equation*}

This paper was mainly motived to make valid the last equation for every
coHom object. Note that for $\underline{hom}^{\Upsilon }\left[ \mathcal{B},%
\mathcal{A}\right] $, with $\mathcal{B}\neq \mathcal{A}$, only a part of
this equation is valid, i.e. 
\begin{equation*}
\underline{hom}^{\Upsilon }\left[ \mathcal{B},\mathcal{A}\right] \backsimeq 
\mathcal{B}_{\chi }\triangleright \mathcal{A}_{\chi }=\underline{hom}\left[ 
\mathcal{B}_{\chi },\mathcal{A}_{\chi }\right] ,
\end{equation*}
since the last part has non sense for a quantum space which has not a
bialgebra structure. To do that we shall define along the next chapters a
twisting process for all quantum spaces.

\section{Cosimplicial quasicomplexes for tensor algebras}

Let us consider a (general) quantum space $\left( \mathbf{A}_{1},\mathbf{A}%
\right) $ with related algebra structure $\left( m,\eta \right) $. Any
linear endomorphism $\Xi \in End_{\mathrm{Vct}}\left[ \mathbf{A}\otimes 
\mathbf{A}\right] $ defines a new product\footnote{%
Since $m$ is unital, and therefore a surjective linear map, any product $%
m^{\prime }$ on $\mathbf{A}$ can be obtained in that way. Indeed, $m^{\prime
}=m_{\Xi }$ with $\Xi =\omega \,m^{\prime }$, being $\omega $ some right
inverse of $m$.} $m_{\Xi }\doteq m\,\Xi $ on $\mathbf{A}$, which can be
called \emph{twisting or deformation of }$m$ \emph{by }$\Xi $. Let us
suppose $\Xi $ is such that $m_{\Xi }$ is associative. If in addition $\Xi $
satisfies 
\begin{equation}
\Xi \,\left( \eta \otimes I\right) =\eta \otimes I\;\;\;and\;\;\;\Xi
\,\left( I\otimes \eta \right) =I\otimes \eta ,  \label{couny}
\end{equation}
then $\eta $ is also a unit for $m_{\Xi }$ and we have a new algebra
structure $\left( m_{\Xi },\eta \right) $ on $\mathbf{A}$, namely the \emph{%
twisting or deformation }$\mathbf{A}_{\Xi }$ \emph{of} $\mathbf{A}$ \emph{by 
}$\Xi $. In general, the pair $\left( \mathbf{A}_{1},\mathbf{A}_{\Xi
}\right) $ is not a quantum space, because the vector space generated
algebraically by $\mathbf{A}_{1}$ through the product $m_{\Xi }$ is not all
of $\mathbf{A}=\mathbf{A}_{\Xi }$ (the equality is for the underlying vector
spaces). Just take $\Xi \equiv 0$. Things change when we restrict ourself to
the set $Aut_{\mathrm{Vct}}\left[ \mathbf{A}\otimes \mathbf{A}\right] $, as
we shall see in \S \textbf{3.1}. On the other hand, if $\left( \mathbf{A}%
_{1},\mathbf{A}\right) $ is a conic quantum space with $\mathbf{A}%
=\bigoplus_{n\in \mathbb{N}_{0}}\mathbf{A}_{n}$, we can ask $\Xi $ to be a 
\emph{bihomogeneous} of degree cero, i.e. $\Xi \left( \mathbf{A}_{r}\otimes 
\mathbf{A}_{s}\right) \subset \mathbf{A}_{r}\otimes \mathbf{A}_{s}$, $%
\forall r,s\in \mathbb{N}_{0}$. In such a case $m_{\Xi }$ gives rise clearly
to the same gradation as $m$. Thus, our objects of interest are homogeneous
linear automorphisms satisfying $\left( \ref{couny}\right) $.

But, what kind of condition(s) must be imposed on $\Xi $ to insure $m_{\Xi }$
is associative (when $m$ is). We shall show immediately that there is a
cosimplicial quasicomplex structure related to every tensor algebra of a
given vector space. Since every linear map in $\mathbf{A}$ can be defined by
one in $\mathbf{A}_{1}^{\otimes }$ through the epimorphism $\mathbf{A}%
_{1}^{\otimes }\twoheadrightarrow \mathbf{A}$, we will study the maps on $%
\mathbf{A}\otimes \mathbf{A}$ in terms of those on $\mathbf{A}_{1}^{\otimes
}\otimes \mathbf{A}_{1}^{\otimes }$ and the mentioned quasicomplex
structure. As one might expect, 2-cocycles will define maps whose related
twisted products are associative.

Our construction of the cosimplicial object and the associated cochain
quasicomplex $\mathsf{C}^{\bullet }$ was inspired by the one developed in 
\cite{maj} and related to twisting of bialgebras,\footnote{%
A similar construction can be also found in the work of Davidov \cite{dav},
on twisting of monoidal structures.} which we shall call $\mathsf{G}%
^{\bullet }$. At the end of this chapter the relationship between these
quasicomplexes will be analyzed.

\subsection{The cosimplicial object}

Consider $\mathbf{V}\in \mathrm{Vct}$ and its associated tensor algebra $%
\mathbf{V}^{\otimes }=\bigoplus_{n\in \mathbb{N}_{0}}\mathbf{V}^{\otimes n}$%
, $\mathbf{V}^{\otimes 0}=\Bbbk $. We shall call 
\begin{equation*}
\mathsf{C}^{n}\left[ \mathbf{V}\right] \subset Aut_{\mathrm{Vct}}\left[
\left( \mathbf{V}^{\otimes }\right) ^{\otimes n}\right] ,\;\;\;n\in \mathbb{N%
}_{0},
\end{equation*}
the subgroup of \emph{n-homogeneous} of degree cero, or simply homogeneous,
linear automorphisms of $\left( \mathbf{V}^{\otimes }\right) ^{\otimes n}$,
i.e. those linear maps $\mathbb{\psi }:\left( \mathbf{V}^{\otimes }\right)
^{\otimes n}\backsimeq \left( \mathbf{V}^{\otimes }\right) ^{\otimes n}$
such that the equalities of sets 
\begin{equation}
\psi \left( \mathbf{V}^{\otimes r_{1}}\otimes ...\otimes \mathbf{V}^{\otimes
r_{n}}\right) =\mathbf{V}^{\otimes r_{1}}\otimes ...\otimes \mathbf{V}%
^{\otimes r_{n}}  \label{incc}
\end{equation}
hold for every $r_{k}\in \mathbb{N}_{0}$ and $k=1...n$. In the finite
dimensional case, such conditions are equivalent to the inclusions 
\begin{equation*}
\psi \left( \mathbf{V}^{\otimes r_{1}}\otimes ...\otimes \mathbf{V}^{\otimes
r_{n}}\right) \subset \mathbf{V}^{\otimes r_{1}}\otimes ...\otimes \mathbf{V}%
^{\otimes r_{n}}.
\end{equation*}
The unit element of $\mathsf{C}^{n}\left[ \mathbf{V}\right] $ is $\mathbb{I}%
^{\otimes n}$, being $\mathbb{I}$ the identity map of $Aut_{\mathrm{Vct}}%
\left[ \mathbf{V}^{\otimes }\right] $. In terms of $n$-fold multi-index $%
R=\left( r_{1},...,r_{n}\right) \in \mathbb{N}_{0}^{\times n}$, we can write 
\begin{equation}
\left( \mathbf{V}^{\otimes }\right) ^{\otimes n}=\bigoplus\nolimits_{R\in 
\mathbb{N}_{0}^{\times n}}\mathbf{V}^{\otimes R};\;\;\;\mathbf{V}^{\otimes
R}\doteq \mathbf{V}^{\otimes r_{1}}\otimes ...\otimes \mathbf{V}^{\otimes
r_{n}},  \label{tee}
\end{equation}
and 
\begin{equation*}
\psi =\bigoplus\nolimits_{R\in \mathbb{N}_{0}^{\times n}}\psi
_{R};\;\;\;\psi _{R}:\mathbf{V}^{\otimes R}\backsimeq \mathbf{V}^{\otimes R}.
\end{equation*}
By $\mathbb{N}_{0}^{\times 0}$ we understand the \emph{singleton} set. For
the units we write $\mathbb{I}^{\otimes n}=\bigoplus_{R\in \mathbb{N}%
_{0}^{\times n}}\mathbb{I}_{R}$ and $\mathbb{I}_{1}=\mathbb{I}$. The
composition of two maps, namely $\psi $ and $\varphi $, takes the form 
\begin{equation*}
\psi \,\varphi =\bigoplus\nolimits_{R\in \mathbb{N}_{0}^{\times n}}\left[
\psi \,\varphi \right] _{R}=\bigoplus\nolimits_{R\in \mathbb{N}_{0}^{\times
n}}\psi _{R}\,\varphi _{R}.
\end{equation*}
This notation lead us to the identification 
\begin{equation}
\mathsf{C}^{n}\left[ \mathbf{V}\right] =\times _{R\in \mathbb{N}_{0}^{\times
n}}Aut_{\mathrm{Vct}}\left[ \mathbf{V}^{\otimes R}\right] .  \label{igd}
\end{equation}
In other terms, each element of $\mathsf{C}^{n}\left[ \mathbf{V}\right] $
can be seen as a section of a fiber bundle with base $\mathbb{N}_{0}^{\times
n}$ and fibers $Aut_{\mathrm{Vct}}\left[ \mathbf{V}^{\otimes R}\right] $.

From now on, whenever we are considering a fixed generic vector space $%
\mathbf{V}$, we shall write $\mathsf{C}^{n}\left[ \mathbf{V}\right] =\mathsf{%
C}^{n}$, just for brevity. Based on groups $\mathsf{C}^{n}$, we now define a
cosimplicial object in the category $\mathrm{Grp}$ as follows.

\bigskip

The multiplication $m_{\otimes }$ of $\mathbf{V}^{\otimes }$ restricted to
each subspace $\mathbf{V}^{\otimes r}\otimes \mathbf{V}^{\otimes s}\subset 
\mathbf{V}^{\otimes }\otimes \mathbf{V}^{\otimes }$ is an injective map
which has as image the vector subspace $\mathbf{V}^{\otimes r+s}\subset 
\mathbf{V}^{\otimes }$. Then, $m_{\otimes }$ defines canonical isomorphisms $%
\mathbf{V}^{\otimes r}\otimes \mathbf{V}^{\otimes s}\backsimeq \mathbf{V}%
^{\otimes r+s}$. On the other hand, consider the functions 
\begin{equation*}
D_{i}^{n}=I\times ...\times D\times ...\times I:\mathbb{N}_{0}^{\times
n+1}\rightarrow \mathbb{N}_{0}^{\times n};\;\;i\in \left\{ 1,...,n\right\} ,
\end{equation*}
where the sum on integers $D:\left( n,m\right) \mapsto n+m$ acts on the $i$%
-th and $i+1$-th factors of $\mathbb{N}_{0}^{\times n+1}$. In terms of these
functions, given $R=\left( r_{1},...,r_{n+1}\right) $ with $n>0$, the
restrictions to $\mathbf{V}^{\otimes R}$ of the maps 
\begin{equation}
m_{i}^{n}=\underset{i\;factors}{\underbrace{I\otimes ...\otimes m_{\otimes }}%
}\otimes \underset{n-i\;factors}{\underbrace{I\otimes ...\otimes I}}:\left( 
\mathbf{V}^{\otimes }\right) ^{\otimes n+1}\rightarrow \left( \mathbf{V}%
^{\otimes }\right) ^{\otimes n},  \label{f1}
\end{equation}
$i\in \left\{ 1,...,n\right\} $, define canonical bijections $m_{i}^{R}:%
\mathbf{V}^{\otimes R}\backsimeq \mathbf{V}^{\otimes D_{i}^{n}\left(
R\right) }$. The latter give rise to group homomorphisms $\delta _{i}^{n}:%
\mathsf{C}^{n}\rightarrow \mathsf{C}^{n+1}$ by the assignment 
\begin{equation}
\psi \mapsto \delta _{i}^{n}\psi =\bigoplus\nolimits_{R\in \mathbb{N}%
_{0}^{\times n+1}}\left[ \delta _{i}^{n}\psi \right] _{R};\;\;\left[ \delta
_{i}^{n}\psi \right] _{R}=\left( m_{i}^{R}\right) ^{-1}\,\psi
_{D_{i}^{n}\left( R\right) }\,m_{i}^{R}.  \label{cl}
\end{equation}
The elements $\delta _{i}^{n}\psi $ are uniquely determined by 
\begin{equation}
m_{i}^{n}\,\delta _{i}^{n}\psi =\psi \,m_{i}^{n},  \label{im}
\end{equation}
since the restriction of above equation to $\mathbf{V}^{\otimes R}$ gives
precisely the last part of Eq. $\left( \ref{cl}\right) $. Using the
identification $\mathbf{V}^{\otimes R}\thickapprox \mathbf{V}^{\otimes
D_{i}^{n}\left( R\right) }$, we can write in compact form 
\begin{equation}
\delta _{i}^{n}\psi \thickapprox \bigoplus\nolimits_{R\in \mathbb{N}%
_{0}^{\times n+1}}\psi _{D_{i}^{n}\left( R\right) }.  \label{iden}
\end{equation}
Let us also define the homomorphisms $\delta _{0,n+1}^{n}:\mathsf{C}%
^{n}\rightarrow \mathsf{C}^{n+1}$, $n\neq 0$, as $\delta _{0}^{n}\psi =%
\mathbb{I}\otimes \psi $ and $\delta _{n+1}^{n}\psi =\psi \otimes \mathbb{I}$%
. For $n=0$ any element $\chi \in $ $\mathsf{C}^{0}$ is a map $1\mapsto
\lambda \in \Bbbk ^{\times }\equiv \Bbbk -\left\{ 0\right\} $. In other
words $\mathsf{C}^{0}\backsimeq \Bbbk ^{\times }$, the multiplicative group
of the field $\Bbbk $. So we can define $\delta _{0,1}^{0}\chi $ by the
assignment $v\in \mathbf{V}^{\otimes }\mapsto \lambda \cdot v\in \mathbf{V}%
^{\otimes }$. In resume, we have built up a set of group homomorphisms $%
\delta _{i}^{n}$, $i\in \left\{ 0,1,...,n+1\right\} $, $n\in \mathbb{N}_{0}$%
. We shall call them \emph{coface operators}.

Consider now, for $i\in \left\{ 0,1,...,n\right\} $, the linear
transformations 
\begin{equation}
\eta _{i}^{n}=\underset{i\;factors}{\underbrace{I\otimes ...\otimes I}}%
\otimes \underset{n+1-i\;factors}{\underbrace{\eta _{\otimes }\otimes
...\otimes I}}:\left( \mathbf{V}^{\otimes }\right) ^{\otimes n}\rightarrow
\left( \mathbf{V}^{\otimes }\right) ^{\otimes n+1}  \label{f2}
\end{equation}
related to the unit map $\eta _{\otimes }:\Bbbk \rightarrow \mathbf{V}%
^{\otimes }$, and the functions 
\begin{equation*}
S_{i}^{n}:\mathbb{N}_{0}^{\times n}\rightarrow \mathbb{N}_{0}^{\times
n+1}\;\;\;/\;\;\;S_{i}^{n}\left( r_{1},...,r_{i},r_{i+1},...,r_{n}\right)
=\left( r_{1},...,r_{i},0,r_{i+1},...,r_{n}\right) .
\end{equation*}
Since $\eta _{\otimes }$ as a map over $\mathbf{V}^{\otimes 0}=\Bbbk $ is
the identity, given $R\in \mathbb{N}_{0}^{\times n}$ the map $\eta _{i}^{n}$
restricted to $\mathbf{V}^{\otimes R}$ defines a canonical bijection $\eta
_{i}^{R}:\mathbf{V}^{\otimes R}\backsimeq \mathbf{V}^{\otimes
S_{i}^{n}\left( R\right) }$. From $\eta _{i}^{n}$ and $S_{i}^{n}$, the group
homomorphisms $\sigma _{i}^{n}:\mathsf{C}^{n+1}\rightarrow \mathsf{C}^{n}$,
with 
\begin{equation}
\psi \mapsto \sigma _{i}^{n}\psi =\bigoplus\nolimits_{R\in \mathbb{N}%
_{0}^{\times n}}\left[ \sigma _{i}^{n}\psi \right] _{R};\;\;\;\left[ \sigma
_{i}^{n}\psi \right] _{R}=\left( \eta _{i}^{R}\right) ^{-1}\,\psi
_{S_{i}^{n}\left( R\right) }\,\eta _{i}^{R},  \label{dege}
\end{equation}
can be defined. The equation $\eta _{i}^{n}\,\sigma _{i}^{n}\psi =\psi
\,\eta _{i}^{n}$ determinates $\sigma _{i}^{n}\psi $ completely. From the
identification $\mathbf{V}^{\otimes R}\thickapprox \mathbf{V}^{\otimes
S_{i}^{n}\left( R\right) }$ we can write 
\begin{equation*}
\sigma _{i}^{n}\psi \thickapprox \bigoplus\nolimits_{R\in \mathbb{N}%
_{0}^{\times n}}\psi _{S_{i}^{n}\left( R\right) }.
\end{equation*}
They will be called \emph{codegeneracies}.

\begin{theorem}
The map $\left[ n+1\right] \mapsto \mathsf{C}^{n}\left[ \mathbf{V}\right] $
defines a \textbf{cosimplicial object} in $\mathrm{Grp}$; i.e. it can be
extended to a functor $\mathbf{\Delta }\rightarrow \mathrm{Grp}$, being $%
\mathbf{\Delta }$ the simplicial category.
\end{theorem}

\begin{proof}
We must show that 
\begin{equation}
\begin{array}{ll}
\delta _{i}^{n+1}\,\delta _{j}^{n}=\delta _{j+1}^{n+1}\,\delta _{i}^{n}, & 
i\leq j; \\ 
&  \\ 
\sigma _{j}^{n-1}\,\sigma _{i}^{n}=\sigma _{i\,}^{n-1}\sigma _{j+1}^{n}, & 
i\leq j; \\ 
&  \\ 
\sigma _{j}^{n+1}\,\delta _{i}^{n+1}=\delta _{i\,}^{n}\,\sigma _{j-1}^{n}, & 
i<j; \\ 
&  \\ 
\sigma _{j}^{n}\,\delta _{i}^{n}=id, & i=j,i=j+1; \\ 
&  \\ 
\sigma _{j}^{n+1}\,\delta _{i}^{n+1}=\delta _{i-1}^{n}\,\sigma _{j}^{n}, & 
i>j+1.
\end{array}
\label{cos}
\end{equation}
In what follows, we will omit the index $n$ on maps $\delta _{i}^{n}$, $%
\sigma _{i}^{n}$, etc. It is easy to see that the opposite equations hold
for $\delta =D$, $\sigma =S$ and for $\delta =m$, $\sigma =\eta $. For
instance, for $i\leq j$ 
\begin{equation}
m_{j}\,m_{i}=m_{i}\,m_{j+1}\ \;and\;\;D_{j}\,D_{i}=D_{i}\,D_{j+1}.  \label{R}
\end{equation}
These are direct consequences of the associative monoidal structures defined
by $\left( D,S\right) $ in $\mathbb{N}_{0}$ and by $\left( m,\eta \right) $
in $\mathbf{V}^{\otimes }$. Using this fact, let us show the first one of
the above equalities. Taking $\psi \in \mathsf{C}^{n}$ we have for all $R\in 
\mathbb{N}_{0}^{\times n+2}$%
\begin{equation*}
\begin{array}{l}
\left[ \delta _{i}\,\delta _{j}\psi \right] _{R}=\left[ \delta _{i}\left(
\delta _{j}\psi \right) \right] _{R}=\left( m_{i}^{R}\right) ^{-1}\,\left[
\delta _{j}\psi \right] _{D_{i}^{n+1}\left( R\right) }\;m_{i}^{R} \\ 
\\ 
=\left( m_{i}^{R}\right) ^{-1}\,\left( m_{j}^{D_{i}^{n+1}\left( R\right)
}\right) ^{-1}\,\psi _{D_{j}^{n}\,D_{i}^{n+1}\left( R\right)
}\;m_{j}^{D_{i}^{n+1}\left( R\right) }\,m_{i}^{R} \\ 
\\ 
=\left( m_{j}^{D_{i}^{n+1}\left( R\right) }\,m_{i}^{R}\right) ^{-1}\,\psi
_{D_{j}^{n}\,D_{i}^{n+1}\left( R\right) }\;m_{j}^{D_{i}^{n+1}\left( R\right)
}\,m_{i}^{R}.
\end{array}
\end{equation*}
Note that, for instance, $m_{j}^{D_{i}^{n+1}\left( R\right) }$ is $%
m_{j}^{n+1}$ restricted to $\mathbf{V}^{\otimes D_{i}^{n+1}\left( R\right) }$%
. Restricting Eq. $\left( \ref{R}\right) $ to $\mathbf{V}^{\otimes R}$ and
taking into account the domains where $m$'s are applied, we arrive at 
\begin{equation*}
\begin{array}{l}
\left[ \delta _{i}\,\delta _{j}\psi \right] _{R}=\left(
m_{i}^{D_{j+1}^{n+1}\left( R\right) }\,m_{j+1}^{R}\right) ^{-1}\,\psi
_{D_{i}^{n}\,D_{j+1}^{n+1}\left( R\right) }\;m_{i}^{D_{j+1}^{n+1}\left(
R\right) }\,m_{j+1}^{R} \\ 
\\ 
=\left( m_{j+1}^{R}\right) ^{-1}\,\left[ \delta _{i}\psi \right]
_{D_{j+1}^{n+1}\left( R\right) }\;m_{j+1}^{R}=\left[ \delta _{j+1}\left(
\delta _{i}\psi \right) \right] _{R}=\left[ \delta _{j+1}\delta _{i}\psi 
\right] _{R}.
\end{array}
\end{equation*}
The others equalities can be shown in an analogous way.\bigskip
\end{proof}

Now, we are going to endow the family of groups $\left\{ \mathsf{C}^{n}\left[
\mathbf{V}\right] \right\} _{n\in \mathbb{N}_{0}}$ with a multiplicative
quasicomplex structure, in a very close fashion to the cochain quasicomplex
for twisting of bialgebras.

\subsection{The multiplicative quasicomplex $\mathsf{C}^{\bullet }$}

Unless a confusion may arise, the index $n$ on maps $\delta _{i}^{n}$ will
be omitted. Consider the map $\partial $ that assigns to each $\psi \in $ $%
\mathsf{C}^{n}$, for every $n\in \mathbb{N}_{0}$, an element of $\mathsf{C}%
^{n+1}$ given by the equation 
\begin{equation}
\partial \psi =\left( \underset{i\;odd}{\overset{\rightarrow }{\prod }}%
\delta _{i}\psi \right) \,\left( \underset{i\;even}{\overset{\leftarrow }{%
\prod }}\delta _{i}\psi \right) ^{-1},  \label{Aa}
\end{equation}
where $i\in \left\{ 0,...,n+1\right\} $ and being 
\begin{equation*}
\underset{\left\{ m_{1}<...<m_{k}\right\} }{\overset{\rightarrow }{\prod }}%
\varphi _{i}\doteq \varphi _{m_{1}}\,\varphi _{m_{2}}\,...\,\varphi
_{m_{k}},\;\;\;\underset{\left\{ m_{1}<...<m_{k}\right\} }{\overset{%
\leftarrow }{\prod }}\varphi _{i}\doteq \varphi _{m_{k}}\,\varphi
_{m_{k-1}}\,...\,\varphi _{m_{1}}.
\end{equation*}
The notation $\partial \psi =\partial _{-}\psi $\thinspace $\left( \partial
_{+}\psi \right) ^{-1}$, 
\begin{equation}
\partial _{-}\psi =\underset{i\;odd}{\overset{\rightarrow }{\prod }}\delta
_{i}\psi ,\;\;\partial _{+}\psi =\underset{i\;even}{\overset{\leftarrow }{%
\prod }}\delta _{i}\psi ,  \label{paq}
\end{equation}
will be very useful. Note that $\partial :\mathsf{C}^{n}\rightarrow \mathsf{C%
}^{n+1}$ is not a group homomorphism for every $n$, but satisfies $\partial 
\mathbb{I}^{\otimes n}=\mathbb{I}^{\otimes n+1}$. The same holds for maps $%
\partial _{\pm }$.

\begin{definition}
We define the \textbf{multiplicative cosimplicial quasicomplex of} $\mathbf{V%
}$ \emph{(}or its tensor algebra\emph{)}, denoted by the pair $\left( 
\mathsf{C}^{\bullet }\left[ \mathbf{V}\right] ,\partial \right) $, as the
sequence 
\begin{equation*}
\begin{diagram}[midshaft] \QTR{sf}{C}^{0}\left[ \QTR{bf}{V}\right] &
\rTo^{\partial} & \QTR{sf}{C}^{1}\left[ \QTR{bf}{V}\right] & \rTo^{\partial}
& \QTR{sf}{C}^{2}\left[ \QTR{bf}{V}\right] & \rTo^{\partial} &
\QTR{sf}{C}^{3}\left[ \QTR{bf}{V}\right]& \rTo^{\partial} & ...\\
\end{diagram}\;\;.\;\;\;\blacksquare
\end{equation*}
\end{definition}

The prefix \emph{quasi}\textbf{\ }is used because in general $\partial
^{2}\left( \mathsf{C}^{n}\left[ \mathbf{V}\right] \right) \neq \left\{ 
\mathbb{I}^{\otimes n+2}\right\} $. However, we shall see for $n=0,1$ the
equality holds. Sometimes we shall omit for brevity this prefix.

Let $\mathrm{Grp}_{\ast }$ be the full subcategory of $\mathrm{Set}_{\ast }$
whose objects are groups based on their respective unit elements (i.e. its
morphisms are unit preserving functions), and denote by $\mathrm{Grp}_{\ast
q}$ the category of quasicomplexes in $\mathrm{Grp}_{\ast }\subset \mathrm{%
Set}_{\ast }$. Its arrows $f^{\bullet }:\left( \mathsf{C}^{\bullet
},\partial _{\mathsf{C}}\right) \rightarrow \left( \mathsf{D}^{\bullet
},\partial _{\mathsf{D}}\right) $ are collections of unit preserving
functions $f^{n}:\mathsf{C}^{n}\rightarrow \mathsf{D}^{n}$ such that
equality $f^{n+1}\,\partial _{\mathsf{C}}=\partial _{\mathsf{D}}\,f^{n}$
holds for all $n$. Then, the pair $\left( \mathsf{C}^{\bullet }\left[ 
\mathbf{V}\right] ,\partial \right) $ is an object of $\mathrm{Grp}_{\ast q}$%
, and the triple $\left( \mathsf{C}^{\bullet }\left[ \mathbf{V}\right]
,\partial _{+},\partial _{-}\right) $ is a multiplicative \emph{parity}
quasicomplex in $\mathrm{Grp}_{\ast }$ (see \cite{ion} and references
therein).

\bigskip

The notions of cocycles, coboundaries, cohomology relation and counitality
for $\left( \mathsf{C}^{\bullet },\partial \right) $ are the following.

$\bullet $ An \textbf{n-cocycle }is a cochain $\chi \in \mathsf{C}^{n}$
satisfying $\partial \chi =\mathbb{I}^{\otimes n+1}$. We indicate $\mathsf{Z}%
^{n}=\mathsf{Z}^{n}\left[ \mathbf{V}\right] $ the set of $n$-cocycles in $%
\mathsf{C}^{n}$.

$\bullet $ Two $n$-cochains $\chi $ and $\chi ^{\prime }$ are called \textbf{%
cohomologous} if there exists $\theta \in \mathsf{C}^{n-1}$, $n\in \mathbb{N}
$, such that the equation $\partial _{-}\theta \,\chi \,\left( \partial
_{+}\theta \right) ^{-1}=\chi ^{\prime }$ holds. We denote $\chi \backsim
_{\theta }\chi ^{\prime }$ if $\chi $ and $\chi ^{\prime }$ are cohomologous
through $\theta $, or simply $\chi \backsim \chi ^{\prime }$, and call $%
\mathsf{Coh}^{n}=\mathsf{Coh}^{n}\left[ \mathbf{V}\right] $ the subset of $%
\mathsf{C}^{n}\times \mathsf{C}^{n}$ of cohomologous pairs.

$\bullet $ An \textbf{n-coboundary }is an element $\omega \in \mathsf{C}^{n}$
with $\mathbb{I}^{\otimes n}\backsim \omega $, i.e. such that $\exists
\theta \in \mathsf{C}^{n-1}$ satisfying $\partial _{-}\theta \,\left(
\partial _{+}\theta \right) ^{-1}=\omega $. We denote $\mathsf{B}^{n}=%
\mathsf{B}^{n}\left[ \mathbf{V}\right] $ the $n$-coboundaries of $\mathsf{C}%
^{n}$. In particular, $\left\{ \mathbb{I}^{\otimes n}\right\} \times \mathsf{%
B}^{n}\subset \mathsf{Coh}^{n}$, for all $n\in \mathbb{N}.$

$\bullet $ An \textbf{n-counital }cochain is a map $\chi \in \mathsf{C}^{n}$
satisfying for $0\leq i\leq n-1$, $\chi \,\eta _{i}^{n-1}=\eta _{i}^{n-1}$,
i.e. $\sigma _{i}^{n-1}\chi =\mathbb{I}^{\otimes n-1}$. Each $\frak{C}^{n}=%
\frak{C}^{n}\left[ \mathbf{V}\right] $ indicates the subgroup of counital
cochains in $\mathsf{C}^{n}$.

\bigskip

We must mention that the cohomology relation is not an equivalence relation
for all $n$. In general it is only reflexive, because $\chi \backsim \chi $
through $\theta =\mathbb{I}^{\otimes n}$.

\bigskip

As we have said at the beginning of this chapter, we are interested in
cochains satisfying Eq. $\left( \ref{couny}\right) $, i.e. $\chi \,\left(
\eta \otimes \mathbb{I}\right) =\left( \eta \otimes \mathbb{I}\right) $ and $%
\chi \,\left( \mathbb{I}\otimes \eta \right) =\left( \mathbb{I}\otimes \eta
\right) $. Since $\eta _{0}^{1}=\left( \mathbb{I}\otimes \eta \right) $ and $%
\eta _{1}^{1}=\left( \eta \otimes \mathbb{I}\right) $, such conditions can
be rewritten $\chi \,\eta _{i}^{1}=\eta _{i}^{1}$, $i=0,1$. These are
exactly the counital cochains. It follows from $\left( \ref{cos}\right) $
that, if $\chi \in \frak{C}^{n}$, then $\delta _{i}^{n}\chi $ is not
necessarily in $\frak{C}^{n+1}$. For example, provided $\sigma
_{i}^{n}\,\delta _{i}^{n}=id$, then $\sigma _{i}^{n}\,\delta _{i}^{n}\chi
=\chi $. Hence, the groups $\frak{C}^{n}$ do not define a cosimplicial
subobject of $\mathsf{C}^{\bullet }$. Nevertheless, for the quasicomplex
structure we have the following result.

\begin{proposition}
The restriction $\partial _{\frak{C}}$ of $\partial $ to the groups $\frak{C}%
^{n}$ defines a quasicomplex $\left( \frak{C}^{\bullet },\partial _{\frak{C}%
}\right) \in \mathrm{Grp}_{\ast q}$, such that the inclusions $\frak{u}:%
\frak{C}^{n}\hookrightarrow \mathsf{C}^{n}$ make $\frak{C}^{\bullet }$ a
subobject of $\mathsf{C}^{\bullet }$.
\end{proposition}

\begin{proof}
We must show that $\partial \left( \frak{C}^{n}\right) \subset \frak{C}%
^{n+1} $. So, let us calculate $\sigma _{j}^{n}\,\partial \psi $ for $\psi $
such that $\sigma _{j}^{n-1}\psi =\mathbb{I}^{\otimes n-1}$, with $%
j=0,...,n-1$. Since $\sigma _{i}^{n}$ is a group homomorphism, 
\begin{equation*}
\sigma _{j}^{n}\,\partial \psi =\left( \underset{i\;odd}{\overset{%
\rightarrow }{\prod }}\sigma _{j}^{n}\,\delta _{i}^{n}\psi \right) \,\left( 
\underset{i\;even}{\overset{\leftarrow }{\prod }}\sigma _{j}^{n}\,\delta
_{i}^{n}\psi \right) ^{-1}.
\end{equation*}
On the other hand, from $\left( \ref{cos}\right) $ (remember that $%
i=0,1,...,n+1$) 
\begin{equation*}
\begin{array}{ll}
\sigma _{j}^{n}\,\delta _{i}^{n}\psi =\,\delta _{i}^{n-1}\sigma
_{j-1}^{n-1}\psi =\,\delta _{i}^{n-1}\mathbb{I}^{\otimes n-1}=\mathbb{I}%
^{\otimes n}, & if\;\;i<j; \\ 
&  \\ 
\sigma _{j}^{n}\,\delta _{i}^{n}\psi =\,\delta _{i-1}^{n-1}\sigma
_{j}^{n-1}\psi =\,\delta _{i-1}^{n-1}\mathbb{I}^{\otimes n-1}=\mathbb{I}%
^{\otimes n}, & if\;\;i>j+1; \\ 
&  \\ 
\sigma _{j}^{n}\,\delta _{i}^{n}\psi =\psi , & if\;\;i=j,\;i=j+1;
\end{array}
\end{equation*}
hence, if $j$ is odd, $\sigma _{j}^{n}\,\partial \psi =\sigma
_{j}^{n}\,\delta _{j}^{n}\psi \,\left( \sigma _{j}^{n}\,\delta
_{j+1}^{n}\psi \right) ^{-1}=\psi \,\psi ^{-1}=\mathbb{I}^{\otimes n}$. The
same can be done for $j$ even. Therefore $\sigma _{j}^{n}\,\partial \psi =%
\mathbb{I}^{\otimes n}$ for $j=0,...,n-1$; that is to say, $\partial \psi
\in \frak{C}^{n+1}$ if $\psi \in \frak{C}^{n}$.
\end{proof}

From now on, $\partial $ will denote the coboundary operator for both $\frak{%
C}^{\bullet }$ and $\mathsf{C}^{\bullet }$. It is worth mentioning $\left( 
\frak{C}^{\bullet },\partial _{+},\partial _{-}\right) $ is \textbf{not} a
parity subquasicomplex of $\left( \mathsf{C}^{\bullet },\partial
_{+},\partial _{-}\right) $.

The subsets of cocycles, coboundaries and cohomologous cochains in $\frak{C}%
^{\bullet }$ will be indicated by $\frak{Z}^{n}\subset \mathsf{Z}^{n}$, $%
\frak{B}^{n}\subset \mathsf{B}^{n}$ and $\frak{Coh}^{n}\subset \mathsf{Coh}%
^{n}$, respectively. Note that $\frak{Z}^{n}=\mathsf{Z}^{n}\cap \frak{C}^{n}$%
, but for $\frak{B}^{n}$ and $\frak{Coh}^{n}$, in general, 
\begin{equation}
\frak{B}^{n}\subset \mathsf{B}^{n}\cap \frak{C}^{n}\;\;and\;\;\frak{Coh}%
^{n}\subset \mathsf{Coh}^{n}\cap \frak{C}^{n}\times \frak{C}^{n}.
\label{inq}
\end{equation}

\subsubsection{Elementary calculations in low \emph{geometric dimensions}}

$\smallskip $

$\mathbf{n=0:}$ Given $\chi :\Bbbk \backsimeq \Bbbk :1\mapsto \lambda $ (so, 
$\chi ^{-1}:1\mapsto \lambda ^{-1}$), it follows that $\delta _{0,1}\chi
:v\mapsto \lambda \cdot v$, and 
\begin{equation*}
\partial \chi =\delta _{1}\chi \,\left( \delta _{0}\chi \right) ^{-1}=%
\mathbb{I}:\mathbf{V}^{\otimes }\rightarrow \mathbf{V}^{\otimes }.
\end{equation*}
Hence $\mathsf{B}^{1}=\left\{ \mathbb{I}\right\} $, and as a consequence $%
\mathsf{B}^{1}=\frak{B}^{1}$. In particular, $\partial $ restricted to $%
\mathsf{C}^{0}$ is the trivial group morphism, and accordingly $\mathsf{Z}%
^{0}=\left\{ \chi \in \mathsf{C}^{0}:\partial \chi =\mathbb{I}\right\} =%
\mathsf{C}^{0}=\Bbbk ^{\times }$. In particular, $\partial ^{2}\chi =\mathbb{%
I}^{\otimes 2}$ for every $\chi \in $ $\mathsf{C}^{0}$. On the other hand,
the unique counital 0-cochain is the identity map $I_{\Bbbk }=\mathbb{I}_{0}$%
. Thus, $\frak{C}^{0}=\left\{ \mathbb{I}_{0}\right\} \subset \mathsf{C}^{0}$.

\bigskip

$\mathbf{n=1:}$ For any $\alpha :\mathbf{V}^{\otimes }\backsimeq \mathbf{V}%
^{\otimes }$ in $\mathsf{C}^{1}$, 
\begin{equation}
\partial \alpha =\delta _{1}\alpha \,\left( \delta _{0}\alpha \right)
^{-1}\,\left( \delta _{2}\alpha \right) ^{-1}=\left( \delta _{1}\alpha
\right) \,\left( \alpha \otimes \alpha \right) ^{-1}.  \label{1c}
\end{equation}
Then, a 1-cocycle is an element $\alpha \in $ $\mathsf{C}^{1}$ such that $%
\delta _{1}\alpha =\alpha \otimes \alpha $, or $\alpha _{r+s}\thickapprox
\alpha _{r}\otimes \alpha _{s}$. In particular, $\alpha _{s}\thickapprox
\alpha _{0}\otimes \alpha _{s}$, so $\alpha _{0}=I_{\Bbbk }$. This means
that $\alpha $ must be counital, i.e. the inclusion $\mathsf{Z}^{1}\subset 
\frak{Z}^{1}$ holds. Moreover, from Eq. $\left( \ref{im}\right) $ the
cocycle condition is equivalent to $\alpha \,m_{\otimes }=m_{\otimes
}\,\delta _{1}\alpha =m_{\otimes }\,\left( \alpha \otimes \alpha \right) $,
hence $\alpha $ is a 1-cocycle \emph{iff }is an algebra automorphism of $%
\mathbf{V}^{\otimes }$ obeying $\alpha \left( \mathbf{V}\right) \subset 
\mathbf{V}$. Then, such algebra automorphisms are in bijection with linear
automorphisms $Aut_{\mathrm{Vct}}\left[ \mathbf{V}\right] =GL\left( n\right) 
$, $n=\dim \mathbf{V}$, i.e. $GL\left( n\right) \backsimeq \mathsf{Z}^{1}=%
\frak{Z}^{1}.$With respect to the square of $\partial $ on $\mathsf{C}^{1}$,
we have 
\begin{equation*}
\partial ^{2}\alpha =\delta _{1}\left( \partial \alpha \right) \,\delta
_{3}\left( \partial \alpha \right) \,\delta _{0}\left( \partial \alpha
\right) ^{-1}\,\delta _{2}\left( \partial \alpha \right) ^{-1}=\delta
_{1}\left( \partial \alpha \right) \,\left( \partial \alpha \otimes \mathbb{I%
}\right) \,\left( \mathbb{I}\otimes \partial \alpha \right) ^{-1}\,\delta
_{2}\left( \partial \alpha \right) ^{-1},
\end{equation*}
and using that $\delta _{1}\left( \alpha \otimes \alpha \right) =\delta
_{1}\alpha \otimes \alpha $, $\delta _{2}\left( \alpha \otimes \alpha
\right) =\alpha \otimes \delta _{1}\alpha $, and that the $\delta _{i}^{n}$%
's are group homomorphisms, it follows from $\left( \ref{1c}\right) $ 
\begin{equation}
\partial ^{2}\alpha =\delta _{1}\left( \delta _{1}\alpha \right) \,\delta
_{2}\left( \delta _{1}\alpha \right) ^{-1}.  \label{18}
\end{equation}
Further, restricting $\delta _{1,2}\left( \delta _{1}\alpha \right) $ to
some $\mathbf{V}^{\otimes r}\otimes \mathbf{V}^{\otimes s}\otimes \mathbf{V}%
^{\otimes t}$, and using the identifications $\left( \ref{iden}\right) $, 
\begin{equation*}
\left[ \delta _{1}\left( \delta _{1}\alpha \right) \right]
_{r,s,t}\thickapprox \left[ \delta _{1}\alpha \right] _{r+s,t}\thickapprox
\alpha _{r+s+t},\;\;\left[ \delta _{2}\left( \delta _{1}\alpha \right) %
\right] _{r,s,t}\thickapprox \left[ \delta _{1}\alpha \right]
_{r,s+t}\thickapprox \alpha _{r+s+t},
\end{equation*}
and $\left[ \partial ^{2}\alpha \right] _{r,s,t}\thickapprox \alpha _{r+s+t}$%
\thinspace $\alpha _{r+s+t}^{-1}=\mathbb{I}_{r,s,t}$ follows, i.e. $\partial
^{2}\alpha =\mathbb{I}^{\otimes 3}$. That means, together with the $n=0$
case, 
\begin{equation}
\func{Im}\left. \partial ^{2}\right| _{\mathsf{C}^{0,1}}=\mathbb{I}^{\otimes
2,3},\;\;or\;\;\mathsf{B}^{1,2}\subset \mathsf{Z}^{1,2}.  \label{cua}
\end{equation}
But $\partial $ is not a morphism of groups for $n=1$.

Given $\alpha ,\beta \in \mathsf{C}^{1}$, they are cohomologous \emph{iff }%
there exists $\lambda \in \mathsf{C}^{0}=\Bbbk ^{\times }$ such that $%
\lambda \cdot \alpha \cdot \lambda ^{-1}=\beta $ \emph{iff} $\alpha =\beta $%
. In particular $\mathsf{B}^{1}=\left\{ \mathbb{I}\right\} =\frak{B}^{1}$,
as we have previously seen. Restricting ourself to counital cochains we have 
$\left( \alpha ,\beta \right) \in \frak{Coh}^{1}$ \emph{iff }there exists $%
\lambda \in \frak{C}^{0}$ such that $\lambda \cdot \alpha \cdot \lambda
^{-1}=\beta $. But $\frak{C}^{0}=\left\{ \mathbb{I}_{0}\right\} $, so $%
\left( \alpha ,\beta \right) \in \frak{Coh}^{1}$ \emph{iff }$\alpha =\beta $%
, and consequently $\frak{Coh}^{1}=\mathsf{Coh}^{1}\cap \frak{C}^{1}\times 
\frak{C}^{1}$.

\bigskip

$\mathbf{n=2:}$ A 2-cocycle is an element $\psi \in $ $\mathsf{C}^{2}$
satisfying $\delta _{1}\psi \,\left( \psi \otimes \mathbb{I}\right) =\delta
_{2}\psi \,\left( \mathbb{I}\otimes \psi \right) $, or 
\begin{equation}
\psi _{r+s,t}\,\left( \psi _{r,s}\otimes \mathbb{I}_{t}\right) \thickapprox
\psi _{r,s+t}\,\left( \mathbb{I}_{r}\otimes \psi _{s,t}\right) .  \label{2cc}
\end{equation}
A counital 2-cochain fulfills the equations $\psi \left( 1\otimes a\right)
=1\otimes a$ and $\psi \left( a\otimes 1\right) =a\otimes 1$, that is to
say, $\psi _{0,s}\thickapprox \mathbb{I}_{s}$ and $\psi _{r,0}\thickapprox 
\mathbb{I}_{r}$. From $\left( \ref{2cc}\right) $ is easy to see that 
\begin{equation}
\psi _{0,s}\thickapprox \psi _{0,0}\cdot \mathbb{I}_{s}\;\;\;and\;\;\;\psi
_{r,0}\thickapprox \psi _{0,0}\cdot \mathbb{I}_{r},  \label{2cc2}
\end{equation}
regarding $\psi _{0,0}$ as an element of $\Bbbk ^{\times }$. That means, 
\begin{equation}
\psi \in \frak{Z}^{2}\;\;\;iff\;\;\;\psi \in \mathsf{Z}^{2}\;and\;\psi
_{0,0}=1.  \label{c2c}
\end{equation}

Suppose $\left( \psi ,\varphi \right) \in \mathsf{Coh}^{2}$, i.e.\emph{\ }$%
\partial _{-}\theta \,\psi =\varphi \,\partial _{+}\theta $ for some $\theta
\in \mathsf{C}^{1}$, or in other terms, $\theta _{r+s}\,\psi
_{r,s}\thickapprox \varphi _{r,s}\,\left( \theta _{r}\otimes \theta
_{s}\right) $. If $\left( \psi ,\varphi \right) \in \mathsf{Coh}^{2}\cap 
\frak{C}^{2}\times \frak{C}^{2}$, then $\psi _{0,0}=\varphi _{0,0}=1$ and $%
\theta _{0}\thickapprox \theta _{0}\,\psi _{0,0}\thickapprox \varphi
_{0,0}\,\left( \theta _{0}\otimes \theta _{0}\right) \thickapprox \left(
\theta _{0}\otimes \theta _{0}\right) $, and in consequence $\theta _{0}\,=1$%
, i.e. $\theta \in \frak{C}^{1}$. Therefore, $\left( \psi ,\varphi \right)
\in \frak{Coh}^{2}$. Summing up (compare with $\left( \ref{inq}\right) $), 
\begin{equation}
\frak{Coh}^{i}=\mathsf{Coh}^{i}\cap \frak{C}^{i}\times \frak{C}%
^{i}\;\;\;and\;\;\;\frak{B}^{i}=\mathsf{B}^{i}\cap \frak{C}%
^{i},\;\;for\;i=1,2.  \label{eqy}
\end{equation}

\subsubsection{\emph{(}Anti\emph{)}bicharacter\ like\ cochains}

Using the terminology of linear forms over a bialgebra \cite{maj}, we shall
say $\psi \in $ $\mathsf{C}^{2}$ is a \emph{bicharacter} if 
\begin{equation}
\delta _{1}\psi =\psi _{23}\,\psi _{13},\;\;\;\;\;\delta _{2}\psi =\psi
_{12}\,\psi _{13},  \label{bica}
\end{equation}
and \emph{anti-bicharacter} (see \cite{cor} for examples in the bialgebra
case) if 
\begin{equation}
\delta _{1}\psi =\psi _{13}\,\psi _{23},\;\;\;\;\;\delta _{2}\psi =\psi
_{13}\,\psi _{12},  \label{abica}
\end{equation}
where $\psi _{12}=\psi \otimes \mathbb{I}$, $\psi _{23}=\mathbb{I}\otimes
\psi $ and $\psi _{13}=\left( \mathbb{I}\otimes f_{\otimes }^{-1}\right)
\,\psi _{12}\,\left( \mathbb{I}\otimes f_{\otimes }\right) $. We are
denoting by $f_{\otimes }$ the flipping operator on $\mathbf{V}^{\otimes
}\otimes \mathbf{V}^{\otimes }$. Using the identifications $\left( \ref{iden}%
\right) $ and following obvious notation, Eq. $\left( \ref{bica}\right) $
can be written 
\begin{equation}
\begin{array}{l}
\psi _{r+s,t}\thickapprox \left( \mathbb{I}_{r}\otimes \psi _{s,t}\right)
\,\left( \mathbb{I}_{r}\otimes f_{s,t}^{-1}\right) \,\left( \psi
_{r,t}\otimes \mathbb{I}_{s}\right) \,\left( \mathbb{I}_{r}\otimes
f_{s,t}\right) , \\ 
\\ 
\psi _{r,s+t}\thickapprox \left( \psi _{r,s}\otimes \mathbb{I}_{t}\right)
\,\left( \mathbb{I}_{r}\otimes f_{s,t}^{-1}\right) \,\left( \psi
_{r,t}\otimes \mathbb{I}_{s}\right) \,\left( \mathbb{I}_{r}\otimes
f_{s,t}\right) .
\end{array}
\label{rst}
\end{equation}
Note that every (anti)bicharacter $\psi $ is always counital; in fact, 
\begin{equation*}
\psi _{0,t}=\psi _{0+0,t}\thickapprox \left( \mathbb{I}_{0}\otimes \psi
_{0,t}\right) \,\left( \mathbb{I}_{0}\otimes f_{0,t}^{-1}\right) \,\left(
\psi _{0,t}\otimes \mathbb{I}_{0}\right) \,\left( \mathbb{I}_{0}\otimes
f_{0,t}\right) \thickapprox \psi _{0,t}^{2},
\end{equation*}
thus $\psi _{0,t}\thickapprox \mathbb{I}_{t}$. Also, $\psi $ is completely
define by $\psi _{1,1}$.

If $\psi $ is a bicharacter, it will be a 2-cocycle \emph{iff} 
\begin{equation}
\psi _{12}\,\psi _{13}\,\psi _{23}=\psi _{23}\,\psi _{13}\,\psi _{12},
\label{ybe}
\end{equation}
i.e. $\psi $ satisfies the \emph{Yang-Baxter} (YB) equation. For an
anti-bicharacter, $\psi $ is a 2-cocycle \emph{iff} 
\begin{equation}
\psi _{12}\,\psi _{23}=\psi _{23}\,\psi _{12}.  \label{te}
\end{equation}
It can be shown it is enough for $\psi _{1,1}$ to satisfy some of the above
equations (depending on the case) in order to ensure $\psi $ is a 2-cocycle.
Then, given $M\in GL\left( n^{2}\right) $ satisfying $\left( \ref{ybe}%
\right) $ (resp. $\left( \ref{te}\right) $), we can define a counital
2-cocycle $\psi $ related to an $n$-dimensional vector space $\mathbf{V}$,
just taking $\psi _{1,1}\left( a_{i}\otimes a_{j}\right)
=M_{ij}^{kl}\,a_{k}\otimes a_{l}$ and extending $\psi $ to all of $\mathbf{V}%
^{\otimes }\otimes \mathbf{V}^{\otimes }$ using Eq. $\left( \ref{bica}%
\right) $ (resp. $\left( \ref{abica}\right) $).

\subsubsection{First and second cohomology spaces}

We have just seen that the cohomology relation in $\mathsf{C}^{1}\left[ 
\mathbf{V}\right] $ is the identity relation, so is an equivalence one. In
particular, it can be used to define the \emph{first cohomology space} $%
H^{1} $ as the quotient of $\mathsf{Z}^{1}$, i.e. $H^{1}\doteq \mathsf{Z}%
^{1}=\frak{Z}^{1}\backsimeq GL\left[ \mathbf{V}\right] $. The same can be
done for $n=2$.

\begin{proposition}
The cohomology relation for $n=2$ is an equivalence relation, and every
equivalence class with some element in $\mathsf{Z}^{2}$ is contained there.
This is also valid for $\frak{C}^{\bullet }$.
\end{proposition}

\begin{proof}
We need to show that the relation is symmetric and transitive. Suppose $%
\left( \varphi ,\psi \right) \in \mathsf{Coh}^{2}$, i.e. there exists $%
\alpha \in \mathsf{C}^{1}$ such that $\partial _{-}\alpha \,\varphi \,\left(
\partial _{+}\alpha \right) ^{-1}=\psi $. Let us first note that $\partial
_{\pm }$ are group homomorphisms when restricted to $\mathsf{C}^{1}$.
Indeed, $\partial _{-}\left( \alpha \right) =\delta _{1}\left( \alpha
\right) $ and $\partial _{+}\left( \alpha \right) =\alpha \otimes \alpha $
for all $\alpha \in \mathsf{C}^{1}$. In particular $\partial _{-}\left(
\alpha ^{-1}\right) =\delta _{1}\left( \alpha ^{-1}\right) =\left( \delta
_{1}\alpha \right) ^{-1}=\left( \partial _{-}\alpha \right) ^{-1}$ and $%
\partial _{+}\left( \alpha ^{-1}\right) =\alpha ^{-1}\otimes \alpha
^{-1}=\left( \partial _{+}\alpha \right) ^{-1}$. Thus 
\begin{equation*}
\partial _{-}\left( \alpha ^{-1}\right) \,\psi \,\partial _{+}\left( \alpha
^{-1}\right) ^{-1}=\varphi ,\;\;\;i.e.\;\left( \psi ,\varphi \right) \in 
\mathsf{Coh}^{2},
\end{equation*}
and accordingly the relation is symmetric. To show the transitivity,
consider again the cohomologous maps $\varphi $ and $\psi $, and a 2-cochain 
$\phi \in \mathsf{C}^{2}$ such that there exists $\beta $ satisfying $%
\partial _{-}\beta \,\psi \,\left( \partial _{+}\beta \right) ^{-1}=\phi $,
i.e. $\psi \backsim \phi $. Then, 
\begin{equation*}
\partial _{-}\left( \beta \,\alpha \right) \,\varphi \,\left( \partial
_{+}\left( \alpha \,\beta \right) \right) ^{-1}=\partial _{-}\beta \,\left(
\partial _{-}\alpha \,\varphi \,\partial _{+}\left( \alpha \right)
^{-1}\right) \,\left( \partial _{+}\beta \right) ^{-1}=\phi .
\end{equation*}
Accordingly, $\varphi \backsim \phi $ through $\beta \alpha $. It rests to
prove each equivalence class containing an element of $\mathsf{Z}^{2}$ is a
subset of $\mathsf{Z}^{2}$, that is to say, given a 2-cocycle $\varphi \in 
\mathsf{C}^{2}$ and a 1-cochain $\alpha \in \mathsf{C}^{1}$, $\psi =\delta
_{1}\alpha \,\varphi \,\left( \alpha \otimes \alpha \right) ^{-1}$ is a
2-cocycle too. But, 
\begin{equation*}
\begin{array}{l}
\partial _{-}\psi =\delta _{1}\left( \delta _{1}\alpha \right) \,\delta
_{1}\varphi \,\left( \delta _{1}\left( \alpha \right) ^{-1}\otimes \alpha
^{-1}\right) \,\left( \delta _{1}\alpha \,\varphi \,\left( \alpha \otimes
\alpha \right) ^{-1}\otimes \mathbb{I}\right) \\ 
\\ 
=\delta _{1}\left( \delta _{1}\alpha \right) \,\delta _{1}\varphi \,\left(
\varphi \otimes \mathbb{I}\right) \,\left( \alpha \otimes \alpha \otimes
\alpha \right) ^{-1}
\end{array}
\end{equation*}
and 
\begin{equation*}
\begin{array}{l}
\left( \partial _{+}\psi \right) ^{-1}=\left( \mathbb{I}\otimes \left(
\alpha \otimes \alpha \right) \,\varphi ^{-1}\,\delta _{1}\left( \alpha
\right) ^{-1}\right) \,\left( \alpha \otimes \delta _{1}\alpha \right)
\,\delta _{2}\left( \varphi \right) ^{-1}\,\delta _{2}\left( \delta
_{1}\alpha \right) ^{-1} \\ 
\\ 
=\left( \alpha \otimes \alpha \otimes \alpha \right) \,\left( \mathbb{I}%
\otimes \varphi ^{-1}\right) \,\delta _{2}\left( \varphi \right)
^{-1}\,\delta _{2}\left( \delta _{1}\alpha \right) ^{-1}.
\end{array}
\end{equation*}
Hence, from the last equations and Eq. $\left( \ref{18}\right) $%
\begin{equation*}
\partial \psi =\delta _{1}\left( \delta _{1}\alpha \right) \,\partial
\varphi \,\delta _{2}\left( \delta _{1}\alpha \right) ^{-1}=\delta
_{1}\left( \delta _{1}\alpha \right) \,\delta _{2}\left( \delta _{1}\alpha
\right) ^{-1}=\partial ^{2}\alpha =\mathbb{I}^{\otimes 3},
\end{equation*}
as we wanted to show.\ 

For $\frak{C}^{\bullet }$ we just have to use the monic $\frak{u}:\frak{C}%
^{\bullet }\hookrightarrow \mathsf{C}^{\bullet }$, which is given by group
monomorphisms.
\end{proof}

Then we can define the second cohomology space $H^{2}\left[ \mathbf{V}\right]
$ as the set of cohomologous classes in $\mathsf{Z}^{2}$, i.e. $H^{2}=\left. 
\mathsf{Z}^{2}\right/ \backsim _{\mathsf{Coh}}$. We are going to see in the
next chapter that any 2-cocycle is a 2-coboundary. Therefore, every $\psi
\in \mathsf{Z}^{2}$ is cohomologous to $\mathbb{I}^{\otimes 2}$, and
accordingly $H^{2}\left[ \mathbf{V}\right] =\left\{ \mathbb{I}^{\otimes
2}\right\} $. Nevertheless, this does not mean all deformations on a given
quantum space generated by $\mathbf{V}$ are trivial, because the maps that
define twisted algebras isomorphic to the original one are those
cohomologous to the identity through certain class of 1-cochains (what we
shall call \emph{admissible}), and not through any of them.

\subsection{Functorial properties of $\mathsf{C}^{\bullet }\left[ \mathbf{V}%
\right] $}

Let us consider the cosimplicial objects $\mathsf{C}_{\mathbf{V}}^{\bullet }:%
\mathbf{\Delta }\rightarrow \mathrm{Grp}$ and $\mathsf{C}_{\mathbf{W}%
}^{\bullet }:\mathbf{\Delta }\rightarrow \mathrm{Grp}$ with $\mathbf{V}$ and 
$\mathbf{W}$ in $\mathrm{Vct}$.

\begin{theorem}
The functors $\mathsf{C}_{\mathbf{V}}^{\bullet }$ and $\mathsf{C}_{\mathbf{W}%
}^{\bullet }$ are naturally equivalent if $\mathbf{V}$ and $\mathbf{W}$ are
isomorphic vector spaces.
\end{theorem}

\begin{proof}
We must show that there exists a family $\left\{ \alpha _{\left[ n+1\right]
}\right\} _{n\in \mathbb{N}_{0}}$ of group isomorphisms 
\begin{equation*}
\mathsf{C}_{\mathbf{V}}^{\bullet }\left[ n+1\right] =\mathsf{C}^{n}\left[ 
\mathbf{V}\right] \backsimeq \mathsf{C}^{n}\left[ \mathbf{W}\right] =\mathsf{%
C}_{\mathbf{W}}^{\bullet }\left[ n+1\right]
\end{equation*}
such that 
\begin{equation*}
\alpha _{\left[ n+2\right] }\,\delta _{i,\mathbf{V}}^{n}=\delta _{i,\mathbf{W%
}}^{n}\,\alpha _{\left[ n+1\right] }\;\;\;and\;\;\;\alpha _{\left[ n+1\right]
}\,\sigma _{i,\mathbf{V}}^{n}=\sigma _{i,\mathbf{W}}^{n}\,\alpha _{\left[ n+2%
\right] }.
\end{equation*}
Let $f$ be an isomorphism between $\mathbf{V}$ and $\mathbf{W}$. Then, the
map $f^{\otimes }:\mathbf{V}^{\otimes }\rightarrow \mathbf{W}^{\otimes }$,
the unique extension of $f$ to $\mathbf{V}^{\otimes }$ as an algebra
homomorphism, defines for every $n\in \mathbb{N}_{0}$ a group isomorphism $%
\mathsf{C}^{n}\left[ \mathbf{V}\right] \backsimeq \mathsf{C}^{n}\left[ 
\mathbf{W}\right] $ given by 
\begin{equation}
\psi \mapsto \psi ^{f}\doteq \left( f^{\otimes }\right) ^{\otimes n}\,\psi
\,\left( \left( f^{\otimes }\right) ^{\otimes n}\right) ^{-1}.  \label{fiefe}
\end{equation}
In particular, since 
\begin{equation*}
\left( f^{\otimes }\right) ^{\otimes n}\left( \mathbf{V}^{\otimes
r_{1}}\otimes ...\otimes \mathbf{V}^{\otimes r_{n}}\right) \subset \mathbf{W}%
^{\otimes r_{1}}\otimes ...\otimes \mathbf{W}^{\otimes r_{n}},
\end{equation*}
or in compact form $\left( f^{\otimes }\right) ^{\otimes n}\left( \mathbf{V}%
^{\otimes R}\right) \subset \mathbf{W}^{\otimes R}$, we can write 
\begin{equation*}
\left( f^{\otimes }\right) ^{\otimes n}=\bigoplus\nolimits_{R\in \mathbb{N}%
_{0}^{\times n}}f^{\otimes R}\;\;\;\;and\;\;\;\;\psi
^{f}=\bigoplus\nolimits_{R\in \mathbb{N}_{0}^{\times n}}f^{\otimes R}\,\psi
_{R}\,\left( f^{\otimes R}\right) ^{-1}.
\end{equation*}
Let us call $f^{n}$ such group isomorphisms, and show that 
\begin{equation}
f^{n+1}\,\delta _{i,\mathbf{V}}^{n}=\delta _{i,\mathbf{W}}^{n}\,f^{n}\;\;%
\;and\;\;\;f^{n+1}\,\sigma _{i,\mathbf{V}}^{n}=\sigma _{i,\mathbf{W}%
}^{n}\,f^{n}.  \label{mag}
\end{equation}
Since $f^{\otimes }$ is an algebra isomorphism, equations 
\begin{equation*}
\left( f^{\otimes }\right) ^{\otimes n}\,m_{i,\mathbf{V}}^{n}=m_{i,\mathbf{W}%
}^{n}\,\left( f^{\otimes }\right) ^{\otimes n+1}\;\;\;and\;\;\;\left(
f^{\otimes }\right) ^{\otimes n+1}\,\eta _{i,\mathbf{V}}^{n}=\eta _{i,%
\mathbf{W}}^{n},
\end{equation*}
and analogous equations to its inverse, hold. Restricting to some $R\in 
\mathbb{N}_{0}^{\times n+1}$ they translate into 
\begin{equation}
f^{\otimes D_{i}^{n}\left( R\right) }\,m_{i,\mathbf{V}}^{R}=m_{i,\mathbf{W}%
}^{R}\,f^{\otimes R}\;\;\;and\;\;\;f^{\otimes S_{i}^{n+1}\left( R\right)
}\,\eta _{i,\mathbf{V}}^{R}=\eta _{i,\mathbf{W}}^{R}.  \label{morf}
\end{equation}
For the first equality of $\left( \ref{mag}\right) $, given $\psi \in 
\mathsf{C}^{n}\left[ \mathbf{V}\right] $ and $R\in \mathbb{N}_{0}^{\times
n+1}$, we have 
\begin{equation*}
\left[ f^{n+1}\,\delta _{i,\mathbf{V}}^{n}\psi \right] _{R}=f^{\otimes R}\,%
\left[ \delta _{i,\mathbf{V}}^{n}\psi \right] _{R}\,\left( f^{\otimes
R}\right) ^{-1}=\,f^{\otimes R}\,\left( \left( m_{i,\mathbf{V}}^{R}\right)
^{-1}\,\psi _{D_{i}^{n}\left( R\right) }\,m_{i,\mathbf{V}}^{R}\right)
\,\left( f^{\otimes R}\right) ^{-1},
\end{equation*}
and using $\left( \ref{morf}\right) $, 
\begin{equation*}
\begin{array}{l}
\left[ f^{n+1}\,\delta _{i,\mathbf{V}}^{n}\psi \right] _{R}=\left( m_{i,%
\mathbf{W}}^{R}\right) ^{-1}\,f^{\otimes D_{i}^{n}\left( R\right) }\,\psi
_{D_{i}^{n}\left( R\right) }\,\,\left( f^{\otimes D_{i}^{n}\left( R\right)
}\right) ^{-1}\,m_{i,\mathbf{W}}^{R} \\ 
\\ 
=\left( m_{i,\mathbf{W}}^{R}\right) ^{-1}\,\psi _{D_{i}^{n}\left( R\right)
}^{f}\,m_{i,\mathbf{W}}^{R}=\left[ \delta _{i,\mathbf{W}}^{n}\psi ^{f}\right]
_{R}=\left[ \delta _{i,\mathbf{W}}^{n}\,f^{n}\psi \right] _{R},
\end{array}
\end{equation*}
as we wanted to show. The second equality of $\left( \ref{mag}\right) $ (and
the $i=0,n+1$ cases) can be shown similarly. Thus, taking $\alpha _{\left[
n+1\right] }=f^{n}$ we have the natural equivalence $\mathsf{C}_{\mathbf{V}%
}^{\bullet }\backsimeq \mathsf{C}_{\mathbf{W}}^{\bullet }$ we are looking
for.
\end{proof}

This equivalence can be translated to the corresponding multiplicative
quasicomplex structures, as follows. Let us call $\mathcal{G}\left[ \mathrm{%
Vct}\right] \subset \mathrm{Vct}$ the groupoid associated to $\mathrm{Vct}$,
i.e. the subcategory of $\mathrm{Vct}$ whose morphisms are isomorphisms.

\begin{theorem}
The map $\mathbf{V}\mapsto \left( \mathsf{C}^{\bullet }\left[ \mathbf{V}%
\right] ,\partial _{\mathbf{V}}\right) $ defines a functor $\mathcal{G}\left[
\mathrm{Vct}\right] \rightarrow \mathrm{Grp}_{\ast q}$.
\end{theorem}

\begin{proof}
Let us show the map $f\mapsto f^{\bullet }$, $f:\mathbf{V}\backsimeq \mathbf{%
W}$ and $f^{n}\psi =\psi ^{f}$ for all $\psi $ in $\mathsf{C}^{n}\left[ 
\mathbf{V}\right] $ (see Eq. $\left( \ref{fiefe}\right) $), extends $\mathbf{%
V}\mapsto \left( \mathsf{C}^{\bullet }\left[ \mathbf{V}\right] ,\partial _{%
\mathbf{V}}\right) $ to a functor. That is to say, $f^{\bullet }$ is such
that $f^{n+1}\,\partial _{\mathbf{V}}=\partial _{\mathbf{W}}\,f^{n}$. Using
Eq. $\left( \ref{mag}\right) $%
\begin{eqnarray}
f^{n+1}\,\partial _{\mathbf{V}}\psi =\left( \underset{i\;odd}{\overset{%
\rightarrow }{\prod }}\left( f^{n+1}\,\delta _{i,\mathbf{V}}\psi \right)
\right) \,\left( \underset{i\;even}{\overset{\leftarrow }{\prod }}\left(
f^{n+1}\,\delta _{i,\mathbf{V}}\psi \right) \right) ^{-1}  \notag \\
&&  \label{conm} \\
=\left( \underset{i\;odd}{\overset{\rightarrow }{\prod }}\delta _{i,\mathbf{W%
}}\,\left( f^{n}\psi \right) \right) \,\left( \underset{i\;even}{\overset{%
\leftarrow }{\prod }}\delta _{i,\mathbf{W}}\,\left( f^{n}\psi \right)
\right) ^{-1}=\partial _{\mathbf{W}}\,f^{n}\psi ,  \notag
\end{eqnarray}
since each $f^{n}$ is a group homomorphism. Hence, the claim of the
proposition follows.
\end{proof}

Note the functor $\mathcal{G}\left[ \mathrm{Vct}\right] \rightarrow \mathrm{%
Grp}_{\ast q}$ defines an homomorphism of groupoids 
\begin{equation*}
\mathcal{G}\left[ \mathrm{Vct}\right] \rightarrow \mathcal{G}\left[ \mathrm{%
Grp}_{\ast q}\right]
\end{equation*}
that assigns to each linear isomorphism$\mathbf{\ }f:\mathbf{V}\backsimeq 
\mathbf{W}$ a quasicomplex isomorphism (see $\left( \ref{fiefe}\right) $) 
\begin{equation*}
f^{\bullet }:\left( \mathsf{C}^{\bullet }\left[ \mathbf{V}\right] ,\partial
_{\mathbf{V}}\right) \backsimeq \left( \mathsf{C}^{\bullet }\left[ \mathbf{W}%
\right] ,\partial _{\mathbf{W}}\right) :\psi \mapsto \psi ^{f}.
\end{equation*}
From this, we have the bijections 
\begin{equation*}
\mathsf{Z}^{n}\left[ \mathbf{V}\right] \backsimeq \mathsf{Z}^{n}\left[ 
\mathbf{W}\right] ,\;\mathsf{B}^{n}\left[ \mathbf{V}\right] \backsimeq 
\mathsf{B}^{n}\left[ \mathbf{W}\right] \;\;and\;\;\mathsf{Coh}^{n}\left[ 
\mathbf{V}\right] \backsimeq \mathsf{Coh}^{n}\left[ \mathbf{W}\right]
\end{equation*}
for all $n\in \mathbb{N}_{0}$. Clearly, all these results are also true for
the subquasicomplexes $\frak{C}^{\bullet }\left[ \mathbf{V}\right] $ of
counital cochains.

\bigskip

In the following subsections we present some structural results. They relate
the quasicomplex of a vector space $\mathbf{W}$ with one of its subspaces $%
\mathbf{V}\subset \mathbf{W}$, and to its dual $\mathbf{W}^{\ast }$. Also
relate $\mathsf{C}^{\bullet }\left[ \mathbf{U}\right] $ and $\mathsf{C}%
^{\bullet }\left[ \mathbf{V}\right] $ with the quasicomplexes $\mathsf{C}%
^{\bullet }\left[ \mathbf{U}\otimes \mathbf{V}\right] $. From these results
we are able to define a full subcategory of $\mathrm{Grp}_{\ast q}$, with
monoidal structure given by the direct product $\times $ of groups.

All constructions and results will be based on $\mathsf{C}^{\bullet }$, but
are also valid, at the quasicomplex level, for the related subquasicomplex $%
\frak{C}^{\bullet }$.

\subsubsection{Linear subspaces and subobjects in $\mathrm{Grp}_{\ast q}$}

Consider the inclusion of vector spaces $\mathbf{V}\subset \mathbf{W}$. The
cochains $\psi \in \mathsf{C}^{n}\left[ \mathbf{W}\right] $ such that $\psi
\left( \mathbf{V}^{\otimes R}\right) \subset \mathbf{V}^{\otimes R}$ holds $%
\forall R\in \mathbb{N}_{0}^{\times n}$, form a subgroup of $\mathsf{C}^{n}%
\left[ \mathbf{W}\right] $ that we shall denote $\mathsf{C}^{n}\left[ 
\mathbf{V}\subset \mathbf{W}\right] $. Given $\psi \in \mathsf{C}^{n}\left[ 
\mathbf{V}\subset \mathbf{W}\right] $, let us indicate by $\left. \psi
\right| _{\mathbf{V}}$ and $\left[ \left. \psi \right| _{\mathbf{V}}\right]
_{R}$ the restrictions of $\psi $ and $\psi _{R}$ to $\left( \mathbf{V}%
^{\otimes }\right) ^{\otimes n}$ and $\mathbf{V}^{\otimes R}$, respectively.
Of course, $\left. \psi \right| _{\mathbf{V}}$ is an element of $\mathsf{C}%
^{n}\left[ \mathbf{V}\right] $ and any element there can be obtained in that
way. We are interested on the relation between the cohomological properties
of the elements $\psi \in \mathsf{C}^{n}\left[ \mathbf{V}\subset \mathbf{W}%
\right] $, and the ones of their associated elements $\left. \psi \right| _{%
\mathbf{V}}$ in $\mathsf{C}^{n}\left[ \mathbf{V}\right] $. To this end,
consider the following result.

\begin{proposition}
The subgroups $\mathsf{C}^{n}\left[ \mathbf{V}\subset \mathbf{W}\right] $
define a multiplicative cosimplicial quasicomplex $\mathsf{C}^{\bullet }%
\left[ \mathbf{V}\subset \mathbf{W}\right] $, such that the canonical
inclusions and projections 
\begin{equation*}
\mathsf{C}^{n}\left[ \mathbf{V}\subset \mathbf{W}\right] \hookrightarrow 
\mathsf{C}^{n}\left[ \mathbf{W}\right] ;\;\;\mathsf{C}^{n}\left[ \mathbf{V}%
\subset \mathbf{W}\right] \twoheadrightarrow \mathsf{C}^{n}\left[ \mathbf{V}%
\right] :\psi \mapsto \left. \psi \right| _{\mathbf{V}},
\end{equation*}
give rise to natural transformations of the respective cosimplicial objects,
and to monic and epi arrows, respectively, in $\mathrm{Grp}_{\ast q}$.
\end{proposition}

\begin{proof}
Note that if $\mathbf{V}\subset \mathbf{W}$ as vector spaces, then $\mathbf{V%
}^{\otimes }\subset \mathbf{W}^{\otimes }$ as algebras. Hence (see Eq. $%
\left( \ref{cl}\right) $), for every $\psi \in \mathsf{C}^{n}\left[ \mathbf{V%
}\subset \mathbf{W}\right] $ and $R\in \mathbb{N}_{0}^{\times n+1}$ we have 
\begin{eqnarray}
\left. \left( \delta _{i,\mathbf{W}}\psi \right) \right| _{\mathbf{V}%
^{\otimes R}}=\left. \left( m_{i,\mathbf{W}}^{R}\right) ^{-1}\,\psi
_{D_{i}^{n}\left( R\right) }\,m_{i,\mathbf{W}}^{R}\right| _{\mathbf{V}%
^{\otimes R}}  \notag \\
&&  \label{ep} \\
=\left( m_{i,\mathbf{V}}^{R}\right) ^{-1}\,\left[ \left. \psi \right| _{%
\mathbf{V}}\right] _{D_{i}^{n}\left( R\right) }\,m_{i,\mathbf{V}}^{R}=\left[
\delta _{i,\mathbf{V}}\left( \left. \psi \right| _{\mathbf{V}}\right) \right]
_{R}.  \notag
\end{eqnarray}
In particular, 
\begin{equation*}
\delta _{i,\mathbf{W}}\psi \left( \mathbf{V}^{\otimes R}\right) \subset 
\mathbf{V}^{\otimes R},\;\;\forall R\in \mathbb{N}_{0}^{\times n+1},
\end{equation*}
i.e. $\delta _{i,\mathbf{W}}\psi \in \mathsf{C}^{n}\left[ \mathbf{V}\subset 
\mathbf{W}\right] $ if $\psi \in \mathsf{C}^{n}\left[ \mathbf{V}\subset 
\mathbf{W}\right] $, and consequently Eq. $\left( \ref{ep}\right) $ can be
written 
\begin{equation}
\left. \left( \delta _{i,\mathbf{W}}\psi \right) \right| _{\mathbf{V}%
}=\delta _{i,\mathbf{V}}\left( \left. \psi \right| _{\mathbf{V}}\right) .
\label{epp}
\end{equation}
Furthermore, the restriction of $\delta _{i,\mathbf{W}}$ to $\mathsf{C}^{n}%
\left[ \mathbf{V}\subset \mathbf{W}\right] $ gives a map 
\begin{equation*}
\delta _{i,\mathbf{V}\subset \mathbf{W}}\doteq \left. \delta _{i,\mathbf{W}%
}\right| _{\mathsf{C}^{n}\left[ \mathbf{V}\subset \mathbf{W}\right] }:%
\mathsf{C}^{n}\left[ \mathbf{V}\subset \mathbf{W}\right] \rightarrow \mathsf{%
C}^{n+1}\left[ \mathbf{V}\subset \mathbf{W}\right] ,
\end{equation*}
such that 
\begin{equation}
\delta _{i,\mathbf{V}\subset \mathbf{W}}\psi =\delta _{i,\mathbf{W}}\psi
\label{epp1}
\end{equation}
and (using Eq. $\left( \ref{epp}\right) $) 
\begin{equation}
\left. \left( \delta _{i,\mathbf{V}\subset \mathbf{W}}\psi \right) \right| _{%
\mathbf{V}}=\delta _{i,\mathbf{V}}\left( \left. \psi \right| _{\mathbf{V}%
}\right) .  \label{epp2}
\end{equation}
The same is true for the maps $\sigma _{i,\mathbf{W}}$. In consequence, from 
$\delta _{i,\mathbf{V}\subset \mathbf{W}}$ and $\sigma _{i,\mathbf{V}\subset 
\mathbf{W}}$ a cosimplicial object 
\begin{equation*}
\mathsf{C}_{\mathbf{V}\subset \mathbf{W}}^{\bullet }:\mathbf{\Delta }%
\rightarrow \mathrm{Grp}:\left[ n+1\right] \mapsto \mathsf{C}^{n}\left[ 
\mathbf{V}\subset \mathbf{W}\right]
\end{equation*}
is defined and, since Eqs. $\left( \ref{epp1}\right) $ and $\left( \ref{epp2}%
\right) $, the proposition follows.
\end{proof}

We must mention the group injections $\mathsf{C}^{n}\left[ \mathbf{V}\right]
\hookrightarrow \mathsf{C}^{n}\left[ \mathbf{V}\subset \mathbf{W}\right] $
given by the map $\phi \mapsto \phi \oplus \mathbb{I}_{\mathbf{K}}$, where
we are decomposing the algebra $\left( \mathbf{W}^{\otimes }\right)
^{\otimes n}$ in the subalgebra $\left( \mathbf{V}^{\otimes }\right)
^{\otimes n}$ and a bilateral ideal $\mathbf{K}$, are not morphisms of
quasicomplexes (the problem is with the coface operators $\delta
_{0,n+1}^{n} $). Hence, the resulting group monomorphisms $\mathsf{C}^{n}%
\left[ \mathbf{V}\right] \hookrightarrow \mathsf{C}^{n}\left[ \mathbf{W}%
\right] $ do not define $\mathsf{C}^{\bullet }\left[ \mathbf{V}\right] $ as
a subobject of $\mathsf{C}^{\bullet }\left[ \mathbf{W}\right] $.

Now, to the wanted result. The cohomological properties of the elements $%
\psi \in \mathsf{C}^{\bullet }\left[ \mathbf{V}\subset \mathbf{W}\right] $,
determinate completely the ones of the restrictions $\left. \psi \right| _{%
\mathbf{V}}\in \mathsf{C}^{\bullet }\left[ \mathbf{V}\right] $. More
precisely,

\begin{proposition}
Let $\psi ,\varphi $ be elements of $\mathsf{C}^{n}\left[ \mathbf{V}\subset 
\mathbf{W}\right] $.

\thinspace

\thinspace

a\emph{)} $\psi \in \mathsf{Z}^{n}\left[ \mathbf{V}\subset \mathbf{W}\right] 
$ \emph{iff }$\psi \in \mathsf{Z}^{n}\left[ \mathbf{W}\right] $ \emph{iff }$%
\left. \psi \right| _{\mathbf{V}}\in \mathsf{Z}^{n}\left[ \mathbf{V}\right] $%
.

\thinspace \thinspace

\thinspace \thinspace

b\emph{)} If $\left( \psi ,\varphi \right) \in \mathsf{Coh}^{n}\left[ 
\mathbf{V}\subset \mathbf{W}\right] $, then $\left( \psi ,\varphi \right)
\in \mathsf{Coh}^{n}\left[ \mathbf{W}\right] $ and $\left( \left. \psi
\right| _{\mathbf{V}},\left. \varphi \right| _{\mathbf{V}}\right) \in 
\mathsf{Coh}^{n}\left[ \mathbf{V}\right] .$
\end{proposition}

\begin{proof}
a) The first part follows from the facts that $\partial _{\mathbf{V}\subset 
\mathbf{W}}\psi =\partial _{\mathbf{W}}\psi $, and 
\begin{equation*}
\left. \left( \partial _{\mathbf{V}\subset \mathbf{W}}\psi \right) \right| _{%
\mathbf{V}}=\left. \partial _{\mathbf{W}}\psi \right| _{\mathbf{V}}=\partial
_{\mathbf{V}}\left( \left. \psi \right| _{\mathbf{V}}\right)
\end{equation*}
(direct consequences of Eqs. $\left( \ref{epp1}\right) $ and $\left( \ref
{epp2}\right) $, resp.).

b) $\left( \psi ,\varphi \right) \in \mathsf{Coh}^{n}\left[ \mathbf{V}%
\subset \mathbf{W}\right] $ means that there exists $\theta \in \mathsf{C}%
^{n-1}\left[ \mathbf{V}\subset \mathbf{W}\right] $ such that 
\begin{equation*}
\left[ \partial _{\mathbf{V}\subset \mathbf{W}}\right] _{-}\theta \,\psi
=\varphi \,\left[ \partial _{\mathbf{V}\subset \mathbf{W}}\right] _{+}\theta
.
\end{equation*}
It is clear this implies $\left( \psi ,\varphi \right) \in \mathsf{Coh}^{n}%
\left[ \mathbf{W}\right] $. In addition, $\left. \theta \right| _{\mathbf{V}%
} $, $\left. \left[ \partial _{\mathbf{V}\subset \mathbf{W}}\right]
_{-}\theta \right| _{\mathbf{V}}$ and $\left. \left[ \partial _{\mathbf{V}%
\subset \mathbf{W}}\right] _{+}\theta \right| _{\mathbf{V}}$ define elements
of $\mathsf{C}^{\bullet }\left[ \mathbf{V}\right] $, thus 
\begin{equation*}
\partial _{\mathbf{V},-}\left( \left. \theta \right| _{\mathbf{V}}\right)
\,\left( \left. \psi \right| _{\mathbf{V}}\right) =\left. \left( \left[
\partial _{\mathbf{V}\subset \mathbf{W}}\right] _{-}\theta \right) \right| _{%
\mathbf{V}}\,\left. \psi \right| _{\mathbf{V}}=\left. \left[ \partial _{%
\mathbf{V}\subset \mathbf{W}}\right] _{-}\theta \,\psi \right| _{\mathbf{V}%
}=\left. \varphi \,\left[ \partial _{\mathbf{V}\subset \mathbf{W}}\right]
_{+}\theta \right| _{\mathbf{V}}=\left( \left. \varphi \right| _{\mathbf{V}%
}\right) \,\partial _{\mathbf{V},+}\left( \left. \theta \right| _{\mathbf{V}%
}\right) ,
\end{equation*}
i.e. $\left. \psi \right| _{\mathbf{V}}\backsim \left. \varphi \right| _{%
\mathbf{V}}$ through $\left. \theta \right| _{\mathbf{V}}$.
\end{proof}

\subsubsection{Coadjoints and products}

Let us restrict ourself to the category \textrm{Vct}$_{f}\subset $\textrm{Vct%
} of finite dimensional $\Bbbk $-vector spaces. In such a case, for each
isomorphism $\mathbf{V}\backsimeq \mathbf{V}^{\ast }$ there is a related
isomorphism of quasicomplexes $\left( \mathsf{C}^{\bullet }\left[ \mathbf{V}%
\right] ,\partial _{\mathbf{V}}\right) \backsimeq \left( \mathsf{C}^{\bullet
}\left[ \mathbf{V}^{\ast }\right] ,\partial _{\mathbf{V}^{\ast }}\right) $
and, in particular, group isomorphisms $\mathsf{C}^{n}\left[ \mathbf{V}%
\right] \backsimeq \mathsf{C}^{n}\left[ \mathbf{V}^{\ast }\right] $ for
every $n$. On the other hand, the natural pairing between $\left( \mathbf{V}%
^{\otimes }\right) ^{\otimes n}$ and $\left( \mathbf{V}^{\ast \otimes
}\right) ^{\otimes n}$, namely 
\begin{equation*}
\left\langle \cdot ,\cdot \right\rangle :\left( \mathbf{V}^{\ast \otimes
}\right) ^{\otimes n}\times \left( \mathbf{V}^{\otimes }\right) ^{\otimes
n}\rightarrow \Bbbk ,
\end{equation*}
gives rise to group anti-homomorphisms (the \emph{transposition maps}) 
\begin{equation*}
\ast _{l}:\mathsf{C}^{n}\left[ \mathbf{V}\right] \rightarrow \mathsf{C}^{n}%
\left[ \mathbf{V}^{\ast }\right] :\psi \mapsto \psi ^{\ast }
\end{equation*}
with $\left\langle \psi ^{\ast }v,w\right\rangle =\left\langle v,\psi
w\right\rangle $, and 
\begin{equation*}
\ast _{r}:\mathsf{C}^{n}\left[ \mathbf{V}^{\ast }\right] \rightarrow \mathsf{%
C}^{n}\left[ \mathbf{V}\right] :\psi \mapsto \psi ^{\ast },
\end{equation*}
being $\left\langle v,\psi ^{\ast }w\right\rangle =$ $\left\langle \psi
v,w\right\rangle $. In both cases we have $\left( \psi ^{\ast }\right)
_{R}=\left( \psi _{R}\right) ^{\ast }$, $\forall $ $R$, and consequently, 
\begin{equation}
\delta _{\mathbf{V}^{\ast },i}\psi ^{\ast }=\left( \delta _{i,\mathbf{V}%
}\psi \right) ^{\ast }\;\;and\;\;\sigma _{\mathbf{V}^{\ast },i}\psi ^{\ast
}=\left( \sigma _{i,\mathbf{V}}\psi \right) ^{\ast }.  \label{star}
\end{equation}
Of course, $\ast _{l}\,\ast _{r}=\ast _{r}\,\ast _{l}=id$; that is to say, $%
\ast _{l}$ and $\ast _{r}$ define mutually inverse natural equivalences of
the corresponding cosimplicial objects. But these maps do not define
morphisms between quasicomplexes $\left( \mathsf{C}^{\bullet }\left[ \mathbf{%
V}\right] ,\partial _{\mathbf{V}}\right) $ and $\left( \mathsf{C}^{\bullet }%
\left[ \mathbf{V}^{\ast }\right] ,\partial _{\mathbf{V}^{\ast }}\right) $,
unless they are composed with the inversion $\psi \mapsto \psi ^{-1}$.

\begin{proposition}
The maps $\psi \in \mathsf{C}^{n}\left[ \mathbf{V}\right] \mapsto \psi
^{\ast -1}=\psi ^{-1\ast }$, $n\in \mathbb{N}_{0}$, define an isomorphism in 
$\mathrm{Grp}_{\ast q}$. The same is true for $\psi \in \mathsf{C}^{n}\left[ 
\mathbf{V}^{\ast }\right] \mapsto \psi ^{\ast -1}$.
\end{proposition}

\begin{proof}
Since they are clearly bijective maps, we just must prove that 
\begin{equation*}
\left( \partial _{\mathbf{V}}\psi \right) ^{\ast -1}=\partial _{\mathbf{V}%
^{\ast }}\left( \psi ^{\ast -1}\right) ,\;\;\forall \psi \in \mathsf{C}^{n}%
\left[ \mathbf{V}\right] ,\;n\in \mathbb{N}_{0}.
\end{equation*}
Using Eq. $\left( \ref{star}\right) $, we have 
\begin{equation*}
\begin{array}{l}
\left( \partial _{\mathbf{V}}\psi \right) ^{\ast }=\left( \underset{i\;odd}{%
\overset{\rightarrow }{\prod }}\delta _{i,\mathbf{V}}\psi \,\underset{i\;even%
}{\overset{\rightarrow }{\prod }}\left( \delta _{i,\mathbf{V}}\psi \right)
^{-1}\right) ^{\ast }=\underset{i\;even}{\overset{\leftarrow }{\prod }}%
\left( \delta _{i,\mathbf{V}}\psi \right) ^{\ast -1}\,\underset{i\;odd}{%
\overset{\leftarrow }{\prod }}\left( \delta _{i,\mathbf{V}}\psi \right)
^{\ast } \\ 
\\ 
=\underset{i\;even}{\overset{\leftarrow }{\prod }}\left( \delta _{\mathbf{V}%
^{\ast },i}\psi ^{\ast }\right) ^{-1}\,\underset{i\;odd}{\overset{\leftarrow 
}{\prod }}\delta _{\mathbf{V}^{\ast },i}\psi ^{\ast }=\underset{i\;even}{%
\overset{\leftarrow }{\prod }}\delta _{\mathbf{V}^{\ast },i}\left( \psi
^{\ast -1}\right) \,\underset{i\;odd}{\overset{\leftarrow }{\prod }}\delta _{%
\mathbf{V}^{\ast },i}\psi ^{\ast },
\end{array}
\end{equation*}
where in the last step we use $\delta $'s are homomorphisms of groups.
Inverting the last member we arrive at the wanted result, i.e. 
\begin{equation*}
\left( \partial _{\mathbf{V}}\psi \right) ^{\ast -1}=\underset{i\;odd}{%
\overset{\rightarrow }{\prod }}\delta _{\mathbf{V}^{\ast },i}\left( \psi
^{\ast -1}\right) \,\underset{i\;even}{\overset{\rightarrow }{\prod }}\left(
\delta _{\mathbf{V}^{\ast },i}\left( \psi ^{\ast -1}\right) \right)
^{-1}=\partial _{\mathbf{V}^{\ast }}\left( \psi ^{\ast -1}\right) .
\end{equation*}
In the same way the another claim can be shown.
\end{proof}

For reason that will become clear later, we indicate $\psi ^{\ast -1}$ by $%
\psi ^{!}$. Following this notation, we can call $\mathsf{C}^{\bullet }\left[
\mathbf{V}^{\ast }\right] $ the \emph{coadjoint}\textbf{\ }quasicomplex to $%
\mathsf{C}^{\bullet }\left[ \mathbf{V}\right] $, and denote it $\mathsf{C}%
^{\bullet }\left[ \mathbf{V}\right] ^{!}$. Its coboundary is $\partial _{%
\mathbf{V}}^{!}\psi =\partial _{\mathbf{V}^{\ast }}\psi =\left( \partial _{%
\mathbf{V}}\left( \psi ^{!}\right) \right) ^{!}$. Given $g:\mathsf{C}%
^{\bullet }\left[ \mathbf{V}\right] \rightarrow \mathsf{C}^{\bullet }\left[ 
\mathbf{W}\right] $, its coadjoint arrow $g^{!}:\mathsf{C}^{\bullet }\left[ 
\mathbf{V}\right] ^{!}\rightarrow \mathsf{C}^{\bullet }\left[ \mathbf{W}%
\right] ^{!}$ can be defined as $g^{!}\left( \psi \right) =\left( g\left(
\psi ^{!}\right) \right) ^{!}$. Of course, $g^{!}$ is (by composition) a
morphism in $\mathrm{Grp}_{\ast q}$ and, $\mathsf{C}^{\bullet }\left[ 
\mathbf{V}\right] ^{!!}=\mathsf{C}^{\bullet }\left[ \mathbf{V}\right] $, $%
g^{!!}=g$.

In these terms, if $\psi \in \mathsf{Z}^{n}\left[ \mathbf{V}\right] $, we
have as a consequence of the last proposition that $\psi ^{!}\in \mathsf{Z}%
^{n}\left[ \mathbf{V}\right] ^{!}$. In addition, if $\psi $ is an
(anti)bicharacter in $\mathsf{C}^{2}\left[ \mathbf{V}\right] $, then $\psi
^{!}$ is an (anti)bicharacter in $\mathsf{C}^{2}\left[ \mathbf{V}\right]
^{!} $.

\bigskip\ 

We also can define a product between our quasicomplexes. Indeed, for every $%
n\in \mathbb{N}_{0}$, and $\mathbf{V},\mathbf{W}\in \mathrm{Vct}$, consider
the groups $\mathsf{C}^{n}\left[ \mathbf{V}\right] \times \mathsf{C}^{n}%
\left[ \mathbf{W}\right] $. We shall indicate $\mathsf{C}^{\bullet }\left[ 
\mathbf{V}\right] \times \mathsf{C}^{\bullet }\left[ \mathbf{W}\right] $ the 
\emph{product cosimplex}$\mathsf{\ }$of $\mathsf{C}^{\bullet }\left[ \mathbf{%
V}\right] $ and $\mathsf{C}^{\bullet }\left[ \mathbf{W}\right] $, with
object function $\left[ n+1\right] \mapsto \mathsf{C}^{n}\left[ \mathbf{V}%
\right] \times \mathsf{C}^{n}\left[ \mathbf{W}\right] $ and coface operators
and codegeneracies 
\begin{equation}
\delta _{i}^{n}\left( \phi ,\varphi \right) =\left( \delta _{i,\mathbf{V}%
}^{n}\phi ,\delta _{i,\mathbf{W}}^{n}\varphi \right) ,\;\;\;\sigma
_{i}^{n}\left( \phi ,\varphi \right) =\left( \sigma _{i,\mathbf{V}}^{n}\phi
,\sigma _{i,\mathbf{W}}^{n}\varphi \right) .  \label{fd}
\end{equation}
Naturally, the associated \emph{product quasicomplex} will have a coboundary
map such that $\partial \left( \phi ,\varphi \right) =\left( \partial _{%
\mathbf{V}}\phi ,\partial _{\mathbf{W}}\varphi \right) $. It is clear that
given a pair of morphisms $f^{\bullet }:\mathsf{C}^{\bullet }\left[ \mathbf{V%
}\right] \rightarrow \mathsf{C}^{\bullet }\left[ \mathbf{V}^{\prime }\right] 
$ and $g^{\bullet }:\mathsf{C}^{\bullet }\left[ \mathbf{W}\right]
\rightarrow \mathsf{C}^{\bullet }\left[ \mathbf{W}^{\prime }\right] $, the
related maps 
\begin{equation*}
f^{n}\times g^{n}:\mathsf{C}^{n}\left[ \mathbf{V}\right] \times \mathsf{C}%
^{n}\left[ \mathbf{W}\right] \rightarrow \mathsf{C}^{n}\left[ \mathbf{V}%
^{\prime }\right] \times \mathsf{C}^{n}\left[ \mathbf{W}^{\prime }\right]
\end{equation*}
define an arrow $\mathsf{C}^{\bullet }\left[ \mathbf{V}\right] \times 
\mathsf{C}^{\bullet }\left[ \mathbf{W}\right] \rightarrow \mathsf{C}%
^{\bullet }\left[ \mathbf{V}^{\prime }\right] \times \mathsf{C}^{\bullet }%
\left[ \mathbf{W}^{\prime }\right] $.

\begin{proposition}
$\mathsf{C}^{\bullet }\left[ \mathbf{V}\right] \times \mathsf{C}^{\bullet }%
\left[ \mathbf{W}\right] $ is a subobject of $\mathsf{C}^{\bullet }\left[ 
\mathbf{V}\otimes \mathbf{W}\right] $ in $\mathrm{Grp}_{\ast q}$.
\end{proposition}

\begin{proof}
Consider the canonical isomorphisms $t_{R}:\left( \mathbf{V}\otimes \mathbf{W%
}\right) ^{\otimes R}\backsimeq \mathbf{V}^{\otimes R}\otimes \mathbf{W}%
^{\otimes R}$. We shall show the group monomorphisms $\frak{j}:\mathsf{C}^{n}%
\left[ \mathbf{V}\right] \times \mathsf{C}^{n}\left[ \mathbf{W}\right]
\hookrightarrow \mathsf{C}^{n}\left[ \mathbf{V}\otimes \mathbf{W}\right] $, 
\begin{equation*}
\frak{j}\left( \phi ,\varphi \right) =\bigoplus\nolimits_{R\in \mathbb{N}%
_{0}^{\times n}}\left[ \frak{j}\left( \phi ,\varphi \right) \right]
_{R}=\bigoplus\nolimits_{R\in \mathbb{N}_{0}^{\times n}}t_{R}^{-1}\,\left(
\phi _{R}\otimes \varphi _{R}\right) \,t_{R},
\end{equation*}
define a morphism in $\mathrm{Grp}_{\ast q}$. To see that, it is enough to
show $\frak{j}\,\delta _{i}^{n}\left( \phi ,\varphi \right) =\delta _{i,%
\mathbf{V}\otimes \mathbf{W}}^{n}\,\frak{j}\left( \phi ,\varphi \right) $.
Given $R\in \mathbb{N}_{0}^{\times n+1}$, by definition of $\frak{j}$, 
\begin{equation*}
\begin{array}{l}
\left[ \frak{j}\,\delta _{i}^{n}\left( \phi ,\varphi \right) \right]
_{R}=t_{R}^{-1}\,\left( \left[ \delta _{i,\mathbf{V}}^{n}\phi \right]
_{R}\otimes \left[ \delta _{i,\mathbf{W}}^{n}\varphi \right] _{R}\right)
\,t_{R} \\ 
\\ 
=t_{R}^{-1}\,\left( \left( m_{i,\mathbf{V}}^{R}\right) ^{-1}\,\phi
_{D_{i}^{n}\left( R\right) }\,m_{i,\mathbf{V}}^{R}\otimes \left( m_{i,%
\mathbf{W}}^{R}\right) ^{-1}\,\varphi _{D_{i}^{n}\left( R\right) }\,m_{i,%
\mathbf{W}}^{R}\right) \,t_{R} \\ 
\\ 
=t_{R}^{-1}\,\left( m_{i,\mathbf{V}}^{R}\otimes m_{i,\mathbf{W}}^{R}\right)
^{-1}\,\left( \phi _{D_{i}^{n}\left( R\right) }\otimes \varphi
_{D_{i}^{n}\left( R\right) }\right) \,\left( m_{i,\mathbf{V}}^{R}\otimes
m_{i,\mathbf{W}}^{R}\right) \,t_{R}.
\end{array}
\end{equation*}
Noting that $m_{i,\mathbf{V}\otimes \mathbf{W}}^{R}=t_{D_{i}^{n}\left(
R\right) }^{-1}\,\left( m_{i,\mathbf{V}}^{R}\otimes m_{i,\mathbf{W}%
}^{R}\right) \,t_{R}\,$, the last member is precisely 
\begin{eqnarray*}
&& \\
\left( m_{i,\mathbf{V}\otimes \mathbf{W}}^{R}\right) ^{-1}\,\left[ \frak{j}%
\left( \phi ,\varphi \right) \right] _{D_{i}^{n}\left( R\right) }\,m_{i,%
\mathbf{V}\otimes \mathbf{W}}^{R}=\left[ \delta _{i,\mathbf{V}\otimes 
\mathbf{W}}^{n}\frak{j}\left( \phi ,\varphi \right) \right] _{R}. \\
&&
\end{eqnarray*}
The same can be done for $\delta _{0,n+1}^{n}$. This concludes our proof.
\end{proof}

In forthcoming sections, we shall identify $\left( \mathbf{V}\otimes \mathbf{%
W}\right) ^{\otimes R}$ with $\mathbf{V}^{\otimes R}\otimes \mathbf{W}%
^{\otimes R}$, enabling us to describe each $\frak{j}\left( \phi ,\varphi
\right) $ as a map 
\begin{equation}
\frak{j}\left( \phi ,\varphi \right) =\bigoplus\nolimits_{R\in \mathbb{N}%
_{0}^{\times n}}\left[ \frak{j}\left( \phi ,\varphi \right) \right]
_{R}=\bigoplus\nolimits_{R\in \mathbb{N}_{0}^{\times n}}\phi _{R}\otimes
\varphi _{R},  \label{moni}
\end{equation}
or equivalently, as the restriction of $\phi \otimes \varphi $ to 
\begin{equation*}
\bigoplus\nolimits_{R\in \mathbb{N}_{0}^{\times n}}\mathbf{V}^{\otimes
R}\otimes \mathbf{W}^{\otimes R}\subset \left( \mathbf{V}^{\otimes }\right)
^{\otimes n}\otimes \left( \mathbf{W}^{\otimes }\right) ^{\otimes n},
\end{equation*}
i.e. $\frak{j}\left( \phi ,\varphi \right) \subset \phi \otimes \varphi $.

\bigskip

Let us denote $\mathsf{C}^{\bullet }\left[ \mathrm{Vct}\right] $ and $%
\mathsf{C}^{\bullet }\left[ \mathrm{Vct}_{f}\right] $ the full subcategories
of $\mathrm{Grp}_{\ast q}$ whose objects are of the form 
\begin{equation*}
\mathsf{C}^{\bullet }\left\{ \mathbf{V}_{i}\right\} \doteq \mathsf{C}%
^{\bullet }\left[ \mathbf{V}_{1}\right] \times ...\times \mathsf{C}^{\bullet
}\left[ \mathbf{V}_{m}\right] ;\;\;\mathbf{V}_{i}\in \mathrm{Vct},
\end{equation*}
and define $\mathsf{C}^{\bullet }\left\{ \mathbf{V}_{i}\right\} ^{!}\doteq 
\mathsf{C}^{\bullet }\left[ \mathbf{V}_{1}\right] ^{!}\times ...\times 
\mathsf{C}^{\bullet }\left[ \mathbf{V}_{m}\right] ^{!}$ for $\mathbf{V}%
_{i}\in \mathrm{Vct}_{f}$. The following result resumes all we have
discussed.

\begin{theorem}
Each quasicomplex $\mathsf{C}^{\bullet }\left\{ \mathbf{V}_{i}\right\} $ is
a subobject of $\mathsf{C}^{\bullet }\left[ \otimes _{i}\mathbf{V}_{i}\right]
$. The map 
\begin{equation*}
\times \,:\left( \mathsf{C}^{\bullet }\left\{ \mathbf{V}_{i}\right\} ,%
\mathsf{C}^{\bullet }\left\{ \mathbf{W}_{i}\right\} \right) \mapsto \mathsf{C%
}^{\bullet }\left\{ \mathbf{V}_{i}\right\} \times \mathsf{C}^{\bullet
}\left\{ \mathbf{W}_{i}\right\}
\end{equation*}
defines a monoidal structure in $\mathsf{C}^{\bullet }\left[ \mathrm{Vct}%
\right] $ with unit $\mathsf{C}^{\bullet }\left[ \Bbbk \right] $; while $!:%
\mathsf{C}^{\bullet }\left\{ \mathbf{V}_{i}\right\} \mapsto \mathsf{C}%
^{\bullet }\left\{ \mathbf{V}_{i}\right\} ^{!}$ gives rise to a covariant
monoidal functor in $\mathsf{C}^{\bullet }\left[ \mathrm{Vct}_{f}\right] $.\
\ \ $\blacksquare $
\end{theorem}

It is important to note, among other things, the set inclusions (under
natural identifications) 
\begin{equation}
\mathsf{Z}^{n}\left[ \mathbf{V}_{1}\right] \times ...\times \mathsf{Z}^{n}%
\left[ \mathbf{V}_{m}\right] \subset \mathsf{Z}^{n}\left[ \mathbf{V}%
_{1}\otimes ...\otimes \mathbf{V}_{m}\right]  \label{sinc}
\end{equation}
hold. This is a consequence of the following more general fact. Since the
image of $\mathcal{G}\left[ \mathrm{Vct}\right] \rightarrow \mathrm{Grp}%
_{\ast q}$ (see \S \textbf{2.3}), is contained in\textrm{\ }$\mathsf{C}%
^{\bullet }\left[ \mathrm{Vct}\right] $, it can be regarded as a functor $%
\mathsf{C}^{\bullet }:\mathcal{G}\left[ \mathrm{Vct}\right] \rightarrow 
\mathsf{C}^{\bullet }\left[ \mathrm{Vct}\right] $. Then, the monics $\frak{j}%
:\mathsf{C}^{\bullet }\left[ \mathbf{V}\right] \times \mathsf{C}^{\bullet }%
\left[ \mathbf{W}\right] \hookrightarrow \mathsf{C}^{\bullet }\left[ \mathbf{%
V}\otimes \mathbf{W}\right] $ define natural transformations $\mathsf{C}%
^{\bullet }\times \mathsf{C}^{\bullet }\rightarrow \mathsf{C}^{\bullet }%
\mathsf{\,}\left( \cdot \otimes \cdot \right) $. The same is true for the
full subfunctor $\mathsf{C}_{f}^{\bullet }:\mathcal{G}\left[ \mathrm{Vct}_{f}%
\right] \rightarrow \mathsf{C}^{\bullet }\left[ \mathrm{Vct}_{f}\right] $.
On the other hand, for each $\mathbf{V}\in \mathrm{Vct}_{f}$, $\mathsf{C}%
^{\bullet }\left[ \mathbf{V}\right] \backsimeq \mathsf{C}^{\bullet }\left[ 
\mathbf{V}^{\ast }\right] $ is a functorial isomorphism for the natural
equivalences $\mathsf{C}_{f}^{\bullet }\,\backsimeq \mathsf{C}_{f}^{\bullet
}\,\ast $ and $\mathsf{C}_{f}^{\bullet }\,\backsimeq !\mathsf{\,C}%
_{f}^{\bullet }$ .

\subsection{Comparison with the quasicomplex for bialgebra twisting}

Let $\mathbf{A}$ be a bialgebra, i.e. $\mathbf{A}\in \Bbbk -\mathrm{Bialg}=%
\mathrm{Bialg}$, with coalgebra structure $\left( \Delta ,\varepsilon
\right) $. Consider the groups $\mathsf{G}^{n}\left[ \mathbf{A}\right] $, $%
n\in $ $\mathbb{N}_{0}$, formed out by the invertible elements of $Lin\left[ 
\mathbf{A}^{\otimes n},\Bbbk \right] $ under the convolution product $\ast $%
. Such a product is given by 
\begin{equation*}
\psi \ast \varphi =\left( \psi \otimes \varphi \right) \,\Delta ^{\left(
n\right) },\;\;\psi ,\varphi \in Lin\left[ \mathbf{A}^{\otimes n},\Bbbk %
\right] ,
\end{equation*}
being $\Delta ^{\left( n\right) }$ the usual coproduct on $\mathbf{A}%
^{\otimes n}$ (in particular, $\Delta ^{\left( 1\right) }=\Delta $ and $%
\Delta ^{\left( 0\right) }=I_{\Bbbk }$). The unit of $\ast $ is $\varepsilon
^{\otimes n}$, the usual counit of $\mathbf{A}^{\otimes n}$.

Denote by the pair $\left( m,\eta \right) $ the algebra structure of $%
\mathbf{A}$, and define (as we have done for a tensor algebra $\mathbf{V}%
^{\otimes }$) the maps $m_{i}^{n}$ and $\eta _{i}^{n}$ as in Eqs. $\left( 
\ref{f1}\right) $ and $\left( \ref{f2}\right) $. Then, the group
homomorphisms $d_{i}^{n}:\mathsf{G}^{n}\left[ \mathbf{A}\right] \rightarrow 
\mathsf{G}^{n+1}\left[ \mathbf{A}\right] $ and $s_{i}^{n}:\mathsf{G}^{n+1}%
\left[ \mathbf{A}\right] \rightarrow \mathsf{G}^{n}\left[ \mathbf{A}\right] $%
, with 
\begin{equation*}
d_{i}^{n}:\psi \mapsto \left\{ 
\begin{array}{l}
\varepsilon \otimes \psi ,\;if\;i=0, \\ 
\\ 
\psi \,m_{i}^{n},\;if\;i\in \left\{ 1,...,n\right\} , \\ 
\\ 
\psi \otimes \varepsilon ,\;if\;i=n+1,
\end{array}
\right.
\end{equation*}
and 
\begin{equation*}
s_{i}^{n}:\psi \mapsto \psi \,\eta _{i}^{n},\;\;i\in \left\{
0,1,...,n\right\} ,
\end{equation*}
define a cosimplicial object $\mathsf{G}_{\mathbf{A}}^{\bullet }:\left[ n+1%
\right] \mapsto \mathsf{G}^{n}\left[ \mathbf{A}\right] $ in $\mathrm{Grp}$ 
\cite{drin}\cite{maj}. Moreover, using the convolution product, the map 
\begin{equation}
d:\psi \mapsto \left( \underset{i\;even}{\overset{\rightarrow }{\prod }}%
d_{i}\psi \right) ^{-1}\ast \left( \underset{i\;odd}{\overset{\leftarrow }{%
\prod }}d_{i}\psi \right) ,  \label{Bb}
\end{equation}
where $i\in \left\{ 0,...,n+1\right\} $, supplies the groups $\mathsf{G}^{n}%
\left[ \mathbf{A}\right] $ with a (cosimplicial multiplicative) cochain
quasicomplex structure.\footnote{%
We reorder the factor, with respect to \cite{maj}, in a way that is
convenient for our work.} Note that the order is reversed w.r.t. the
coboundary $\partial $ of the quasicomplexes $\mathsf{C}^{\bullet }$.

We recall that $\psi \in \mathsf{G}^{n}\left[ \mathbf{A}\right] $ is a
counital $n$-cochain if satisfies $s_{i}^{n}\psi =\psi \,\eta
_{i}^{n}=\varepsilon ^{\otimes n}$ for all $i$. Also recall, a bicharacter
is a map $\psi \in \mathsf{G}^{2}\left[ \mathbf{A}\right] $ such that 
\begin{equation*}
d_{1}\psi =\psi \,\left( m\otimes I\right) =\psi _{13}\ast \psi
_{23},\;\;d_{2}\psi =\psi \,\left( I\otimes m\right) =\psi _{13}\ast \psi
_{12},
\end{equation*}
while an anti-bicharacter \cite{cor} fulfill the opposite equations: $%
d_{1}\psi =\psi _{23}\ast \psi _{13}$ and $d_{2}\psi =\psi _{12}\ast \psi
_{13}$.

\bigskip

Since $d$ is a unit preserving map, the pair $\left( \mathsf{G}^{\bullet }%
\left[ \mathbf{A}\right] ,d\right) $ is a quasicomplex in $\mathrm{Grp}%
_{\ast }$, and the factorization of $d$ in $d_{+},d_{-}$, namely $d\psi
=\left( d_{+}\psi \right) ^{-1}$\thinspace $d_{-}\psi $, makes the triple $%
\left( \mathsf{G}^{\bullet }\left[ \mathbf{A}\right] ,d_{+},d_{-}\right) $ a
parity quasicomplex there. For later convenience, we mention the following
known facts without proof.

\begin{proposition}
If a pair of bialgebras $\mathbf{A}$ and $\mathbf{B}$ are homomorphic by
means of a map $\alpha :\mathbf{A}\rightarrow \mathbf{B}$, then there exists
a natural transformation $\mathsf{G}_{\mathbf{B}}^{\bullet }\rightarrow 
\mathsf{G}_{\mathbf{A}}^{\bullet }$ of cosimplicial objects given by
functions 
\begin{equation}
\psi \in \mathsf{G}^{n}\left[ \mathbf{B}\right] \mapsto \psi \,\alpha
^{\otimes n}\in \mathsf{G}^{n}\left[ \mathbf{A}\right] .  \label{nt}
\end{equation}
Moreover, the same functions define a morphism $\left( \mathsf{G}^{\bullet }%
\left[ \mathbf{B}\right] ,d_{\mathbf{B}}\right) \rightarrow \left( \mathsf{G}%
^{\bullet }\left[ \mathbf{A}\right] ,d_{\mathbf{A}}\right) $ of
quasicomplexes, becoming the assignment $\mathbf{A}\mapsto \mathsf{G}%
^{\bullet }\left[ \mathbf{A}\right] $ into a contravariant functor $\mathrm{%
Bialg}\rightarrow \mathrm{Grp}_{\ast q}$.\ \ \ $\blacksquare $
\end{proposition}

When the bialgebra $\mathbf{A}$ is a tensor algebra $\mathbf{C}^{\otimes }$,%
\footnote{%
Every tensor algebra $\mathbf{C}^{\otimes }$ has, fixing a basis $\left\{
v_{i}\right\} $ on $\mathbf{C}$, a bialgebra structure given by the
extensions of $\Delta :v_{i}\mapsto v_{i}\otimes v_{i}$ and $\varepsilon
:v_{i}\mapsto 1$ to algebra homomorphisms. This is cocommutative.} we can
ask which is the relation between $\mathsf{G}^{\bullet }\left[ \mathbf{C}%
^{\otimes }\right] $ and $\mathsf{C}^{\bullet }\left[ \mathbf{C}\right] $.
In order to find this relation, given $\psi \in \mathsf{G}^{n}\left[ \mathbf{%
C}^{\otimes }\right] $, consider the endomorphism of $\left( \mathbf{C}%
^{\otimes }\right) ^{\otimes n}$ of the form\footnote{%
Given a bialgebra $\mathbf{A}$ and two linear maps $\alpha ,\beta :\mathbf{A}%
\rightarrow \mathbf{B}_{\alpha ,\beta }$, $\alpha \ast \beta $ will always
mean $\left( \alpha \otimes \beta \right) \,\Delta :\mathbf{A}\rightarrow 
\mathbf{B}_{\alpha }\otimes \mathbf{B}_{\beta }$.} 
\begin{equation*}
\psi \ast \mathbb{I}^{\otimes n}\ast \psi ^{-1}=\left( \psi \otimes \mathbb{I%
}^{\otimes n}\otimes \psi ^{-1}\right) \,\left( \mathbb{I}^{\otimes
n}\otimes \Delta ^{\left( n\right) }\right) \,\Delta ^{\left( n\right) }.
\end{equation*}
Given another $\varphi \in \mathsf{G}^{n}\left[ \mathbf{C}^{\otimes }\right] 
$, we have 
\begin{equation*}
\left( \psi \ast \varphi \right) \ast \mathbb{I}^{\otimes n}\ast \left( \psi
\ast \varphi \right) ^{-1}=\psi \ast \varphi \ast \mathbb{I}^{\otimes n}\ast
\varphi ^{-1}\ast \psi ^{-1}=\left( \varphi \ast \mathbb{I}^{\otimes n}\ast
\varphi ^{-1}\right) \,\left( \psi \ast \mathbb{I}^{\otimes n}\ast \psi
^{-1}\right) ,
\end{equation*}
and obviously for $\varepsilon ^{\otimes n}$, $\varepsilon ^{\otimes n}\ast 
\mathbb{I}^{\otimes n}\ast \varepsilon ^{\otimes n}=\mathbb{I}^{\otimes n}$.
It follows that $\psi ^{-1}\ast \mathbb{I}^{\otimes n}\ast \psi $ is the
inverse of $\psi \ast \mathbb{I}^{\otimes n}\ast \psi ^{-1}$, hence the
latter is a linear automorphism of $\left( \mathbf{C}^{\otimes }\right)
^{\otimes n}$. Therefore, the function $\digamma ^{n}:\psi \mapsto \psi
^{-1}\ast \mathbb{I}^{\otimes n}\ast \psi $ is a group \textbf{anti}%
-homomorphism, i.e. $\digamma ^{n}\left( \varepsilon ^{\otimes n}\right) =%
\mathbb{I}^{\otimes n}$ and $\digamma ^{n}\left( \psi \ast \varphi \right)
=\digamma ^{n}\varphi \,\digamma ^{n}\psi $. Suppose from now on the
comultiplication on $\mathbf{C}^{\otimes }$ is such that the inclusion $%
\Delta \left( \mathbf{C}\right) \subset \mathbf{C}\otimes \mathbf{C}$ holds.
In other words, we are saying $\mathbf{C}$ is a coalgebra. Then, since $%
\Delta ^{\left( n\right) }$ is an algebra morphism, for every $R\in \mathbb{N%
}_{0}^{\times n}$ 
\begin{equation*}
\func{Im}\Delta _{R}\subset \mathbf{C}^{\otimes R}\otimes \mathbf{C}%
^{\otimes R},\;\;\Delta _{R}\doteq \left. \Delta ^{\left( n\right) }\right|
_{\mathbf{C}^{\otimes R}},
\end{equation*}
and accordingly each map $\digamma ^{n}\psi $ is homogeneous, thus an $n$%
-cochain in $\mathsf{C}^{n}\left[ \mathbf{C}\right] $. Moreover,

\begin{theorem}
Let $\mathbf{C}$ be a coalgebra. Then there exist group anti-homomorphisms $%
\digamma ^{n}:\mathsf{G}^{n}\left[ \mathbf{C}^{\otimes }\right] \rightarrow 
\mathsf{C}^{n}\left[ \mathbf{C}\right] $, $n\in \mathbb{N}_{0}$, that define:

$\bullet $ a natural transformation of the functors $\mathsf{G}_{\mathbf{C}%
^{\otimes }}^{\bullet }\mathbf{\ }$and $\mathsf{C}_{\mathbf{C}}^{\bullet }$,
and

$\bullet $ a morphism of cochain quasicomplexes $\left( \mathsf{G}^{\bullet }%
\left[ \mathbf{C}^{\otimes }\right] ,d\right) \rightarrow \left( \mathsf{C}%
^{\bullet }\left[ \mathbf{C}\right] ,\partial \right) $ in $\mathrm{Grp}%
_{\ast q}$.
\end{theorem}

\begin{proof}
For the first statement we must show that $\digamma ^{n+1}\,d_{i}^{n}=\delta
_{i}^{n}\,\digamma ^{n}$ and $\digamma ^{n+1}\,s_{i}^{n}=\sigma
_{i}^{n}\,\digamma ^{n+2}$. Let us apply the first member of the first
equality to $\psi \in \mathsf{G}^{n}\left[ \mathbf{C}^{\otimes }\right] $.
We have, for $i\neq 0,n+1$, the element of $\mathsf{C}^{n+1}\left[ \mathbf{C}%
\right] $%
\begin{eqnarray*}
\digamma ^{n+1}\,d_{i}^{n}\psi =\digamma ^{n+1}\,\left( \psi
\,m_{i}^{n}\right) =\left( \psi \,m_{i}^{n}\right) \ast \mathbb{I}^{\otimes
n+1}\ast \left( \psi ^{-1}\,m_{i}^{n}\right) \\
&& \\
=\left( \psi \,m_{i}^{n}\otimes \mathbb{I}^{\otimes n+1}\otimes \psi
^{-1}\,m_{i}^{n}\right) \,\left( \mathbb{I}^{\otimes n+1}\otimes \Delta
^{\left( n+1\right) }\right) \,\Delta ^{\left( n+1\right) },
\end{eqnarray*}
which restricted to $\mathbf{C}^{\otimes R}$, for some $R\in \mathbb{N}%
_{0}^{\times n+1}$, gives 
\begin{equation*}
\left[ \digamma ^{n+1}\,d_{i}^{n}\psi \right] _{R}=\left( \psi
\,m_{i}^{R}\otimes \mathbb{I}_{R}\otimes \psi ^{-1}\,m_{i}^{R}\right)
\,\left( \mathbb{I}_{R}\otimes \Delta _{R}\right) \,\Delta _{R}.
\end{equation*}
Writing $\mathbb{I}_{R}=\left( m_{i}^{R}\right) ^{-1}\,m_{i}^{R}$ and using 
\begin{equation}
\left( m_{i}^{R}\otimes m_{i}^{R}\right) \,\Delta _{R}=\Delta
_{D_{i}^{n+1}\left( R\right) }\,m_{i}^{R},  \label{md}
\end{equation}
it follows that $\left[ \digamma ^{n+1}\,d_{i}^{n}\psi \right] _{R}$ is
equal to 
\begin{equation*}
\begin{array}{l}
\left( m_{i}^{R}\right) ^{-1}\,\left[ \left( \psi \otimes \mathbb{I}%
_{D_{i}^{n+1}\left( R\right) }\otimes \psi ^{-1}\right) \,\left( \mathbb{I}%
_{D_{i}^{n+1}\left( R\right) }\otimes \Delta _{D_{i}^{n+1}\left( R\right)
}\right) \,\Delta _{D_{i}^{n+1}\left( R\right) }\right] \,m_{i}^{R} \\ 
\\ 
=\left( m_{i}^{R}\right) ^{-1}\,\left[ \psi \ast \mathbb{I}^{\otimes n}\ast
\psi ^{-1}\right] _{D_{i}^{n+1}\left( R\right) }\,m_{i}^{R}=\left[ \delta
_{i}^{n}\,\digamma ^{n}\psi \right] _{R}.
\end{array}
\end{equation*}
If for instance $i=0$, then\footnote{%
We are using that given a pair of coalgebras $\left( \mathbf{A},\Delta _{%
\mathbf{A}}\right) $ and $\left( \mathbf{B},\Delta _{\mathbf{B}}\right) $,
and linear maps $\alpha ,\beta $ and $\gamma ,\delta $ with domain $\mathbf{A%
}$ and $\mathbf{B}$, resp., the convolution product related to the usual
coalgebra in $\mathbf{A}\otimes \mathbf{B}$ satisfies $\left( \alpha \otimes
\gamma \right) \ast \left( \beta \otimes \delta \right) =\left( \alpha \ast
\beta \right) \otimes \left( \gamma \ast \delta \right) $.} 
\begin{eqnarray*}
\digamma ^{n+1}\,\left( \varepsilon \otimes \psi \right) =\left( \varepsilon
\otimes \psi \right) \ast \mathbb{I}^{\otimes n+1}\ast \left( \varepsilon
\otimes \psi ^{-1}\right) =\left( \varepsilon \otimes \psi \right) \ast
\left( \mathbb{I}\otimes \mathbb{I}^{\otimes n}\right) \ast \left(
\varepsilon \otimes \psi ^{-1}\right) \\
&& \\
=\left( \varepsilon \ast \mathbb{I}\ast \varepsilon \right) \otimes \left(
\psi \ast \mathbb{I}^{\otimes n}\ast \psi ^{-1}\right) =\mathbb{I}\otimes
\digamma ^{n}\psi .
\end{eqnarray*}
The other equalities can be shown analogously. We must use for $\digamma
^{n+1}\,s_{i}^{n}=\sigma _{i}^{n}\,\digamma ^{n+2}$ that 
\begin{equation}
\varepsilon _{R}=\varepsilon _{D_{i}^{n+1}\left( R\right) }\,m_{i}^{R},
\label{me}
\end{equation}
what follows from the algebra map character of $\varepsilon $.

Let us prove the second statement. $\digamma ^{\bullet }:\mathsf{G}^{\bullet
}\left[ \mathbf{C}^{\otimes }\right] \rightarrow \mathsf{C}^{\bullet }\left[ 
\mathbf{C}\right] $ is a morphism in $\mathrm{Grp}_{\ast q}$ if and only if
equation $\digamma ^{n+1}\,d=\partial \,\digamma ^{n}$ is satisfied. Using $%
\digamma ^{n+1}\,d_{i}^{n}=\delta _{i}^{n}\,\digamma ^{n}$ and the fact that
each $\digamma ^{n}$ is an anti-homomorphism of groups, the last equation is
immediately fulfilled (compare Eqs. $\left( \ref{Aa}\right) $ with $\left( 
\ref{Bb}\right) $). That concludes our proof.
\end{proof}

The above theorem describes completely the relation between $\mathsf{G}%
^{\bullet }\left[ \mathbf{C}^{\otimes }\right] $ and $\mathsf{C}^{\bullet }%
\left[ \mathbf{C}\right] $ (provided $\mathbf{C}$ is a coalgebra) as
multiplicative cosimplicial quasicomplexes. Nevertheless, let us look a
little closer at the maps $\digamma ^{n}$. Observe that 
\begin{equation*}
\digamma ^{n}\psi =\mathbb{I}^{\otimes n}\;\;iff\;\;\psi \ast \mathbb{I}%
^{\otimes n}=\mathbb{I}^{\otimes n}\ast \psi
\end{equation*}
\emph{iff } $\psi \ast \varphi =\varphi \ast \psi $ for all $\varphi \in 
\mathsf{G}^{n}\left[ \mathbf{C}^{\otimes }\right] $. In other words, $%
\digamma ^{n}\psi =\mathbb{I}^{\otimes n}$ \emph{iff }$\psi \in \mathcal{Z}%
^{n}$, the center of $\mathsf{G}^{n}\left[ \mathbf{C}^{\otimes }\right] $.
In particular, if the coproduct $\Delta $ of $\mathbf{C}^{\otimes }$ is
cocommutative, then $\digamma ^{n}\psi =\mathbb{I}^{\otimes n}$ for all $%
\psi $ in $\mathsf{G}^{n}\left[ \mathbf{C}^{\otimes }\right] $, for all $n$.
In any case, we may consider the quotient groups $\mathsf{G}_{\mathcal{Z}%
}^{n}\left[ \mathbf{C}^{\otimes }\right] \doteq \left. \mathsf{G}^{n}\left[ 
\mathbf{C}^{\otimes }\right] \right/ \mathcal{Z}^{n}$ and the corresponding
injections\footnote{%
As we will see later, from the point of view of twisting of bialgebras, the
groups of interest are precisely $\mathsf{G}_{\mathcal{Z}}^{n}$, instead of $%
\mathsf{G}^{n}$.} $i^{n}:\mathsf{G}_{\mathcal{Z}}^{n}\left[ \mathbf{C}%
^{\otimes }\right] \hookrightarrow \mathsf{C}^{n}\left[ \mathbf{C}\right] $
induced by each $\digamma ^{n}$. Thus, roughly speaking, the elements of $%
\mathsf{G}^{n}$ form, up to the center, a subset of $\mathsf{C}^{n}$.
Moreover, we have:

\begin{theorem}
The groups $\mathsf{G}_{\mathcal{Z}}^{n}\left[ \mathbf{C}^{\otimes }\right] $
define a multiplicative cosimplicial quasicomplex in such a way that the
commutative diagrams of group \emph{(}anti\emph{)}homomorphisms 
\begin{equation*}
\begin{diagram}[midshaft] & &\QTR{sf}{G}_{\mathcal{Z}}^{n}\left[
\QTR{bf}{C}^{\otimes }\right] & & \\ &\ruTo^{p ^{n}} & & \rdTo^{i^{n}} & \\
\QTR{sf}{G}^{n}\left[ \QTR{bf}{C}^{\otimes }\right] & & \rTo^{\digamma ^{n}}
& &\QTR{sf}{C}^{n}\left[ \QTR{bf}{C}\right] \\ \end{diagram}
\end{equation*}
where $p^{n}:\mathsf{G}^{n}\left[ \mathbf{C}^{\otimes }\right]
\twoheadrightarrow \mathsf{G}_{\mathcal{Z}}^{n}\left[ \mathbf{C}^{\otimes }%
\right] $ is the canonical projection, give rise to a commutative diagram of
natural transformations among their related cosimplicial objects, and a
commutative diagram of morphisms in $\mathrm{Grp}_{\ast q}$ of their
corresponding cochain quasicomplex.\ \ \ $\blacksquare $
\end{theorem}

The theorem follows from the lemmas below. To enunciate the first one, let
us note that given $R\in \mathbb{N}_{0}^{\times n}$, the pair $\left( \Delta
_{R},\varepsilon _{R}\right) $ defines a coalgebra structure on $\mathbf{C}%
^{\otimes R}$. Then, the set $Lin\left[ \mathbf{C}^{\otimes }\right] $ is
supplied with a convolution product $\ast $ which enable us to define the
groups $\mathsf{G}^{R}=\mathsf{G}^{R}\left[ \mathbf{C}^{\otimes }\right] $,
given by the invertible linear forms $\mathbf{C}^{\otimes R}\rightarrow
\Bbbk $. Of course, the unit is $\varepsilon _{R}$. On the other hand,
because of Eq. $\left( \ref{tee}\right) $, every element $\psi $ in $\mathsf{%
G}^{n}=\mathsf{G}^{n}\left[ \mathbf{C}^{\otimes }\right] $ is defined by a
family of linear forms $\psi _{R}\in \mathsf{G}^{R}$, one for each $R\in 
\mathbb{N}_{0}^{\times n}$, being $\psi _{R}=\left. \psi \right| _{\mathbf{C}%
^{\otimes R}}$. In other words, the elements of $\mathsf{G}^{n}$ can be
characterized as functions $R\in \mathbb{N}_{0}^{\times n}\mapsto \psi
_{R}\in \mathsf{G}^{R}$. Furthermore, the restrictions $\psi \mapsto \psi
_{R}$ are group epimorphisms, since $\left[ \psi \ast \varphi \right]
_{R}=\psi _{R}\ast \varphi _{R}$.

\begin{lemma}
The following statements are equivalent:

\qquad i\emph{)} $\psi \in \mathcal{Z}^{n}\subset \mathsf{G}^{n}$.

\qquad ii\emph{)} The equation $\psi _{R}\ast \gamma \,=\gamma \ast \psi
_{R} $ holds $\forall \gamma \in \mathsf{G}^{R}$, $\forall R\in \mathbb{N}%
_{0}^{\times n}$.
\end{lemma}

\begin{proof}
A linear form $\psi $ is in the center of $\mathsf{G}^{n}$ \emph{iff }$\psi
_{R}\ast \varphi _{R}=\varphi _{R}\ast \psi _{R}$ for all $R\in \mathbb{N}%
_{0}^{\times n}$, for all $\varphi \in \mathsf{G}^{n}$. But, to give an
arbitrary element $\varphi \in \mathsf{G}^{n}$ is the same as giving for
each $R$ an arbitrary element of $\mathsf{G}^{R}$, by the characterization
above. Then, $\psi \in \mathcal{Z}^{n}$ \emph{iff }$\psi _{R}\ast \gamma
\,=\gamma \ast \psi _{R}$ for all $\gamma \in \mathsf{G}^{R}$ and $R\in 
\mathbb{N}_{0}^{\times n}$.
\end{proof}

\begin{lemma}
For all $n\in \mathbb{N}_{0}$, $d_{i}^{n}\left( \mathcal{Z}^{n}\right)
\subset \mathcal{Z}^{n+1}$ and $s_{i}^{n}\left( \mathcal{Z}^{n+1}\right)
\subset \mathcal{Z}^{n}$.
\end{lemma}

\begin{proof}
We shall only prove the inclusions $d_{i}^{n}\left( \mathcal{Z}^{n}\right)
\subset \mathcal{Z}^{n+1}$, since the another ones can be proven in a
similar way. Let $\psi $ and $\varphi $ be arbitrary linear forms in $%
\mathsf{G}^{n}$ and $\mathsf{G}^{n+1}$, respectively. Then, for $i\neq 0,n+1$%
, \textsf{\ } 
\begin{equation*}
\left( d_{i}^{n}\psi \right) \ast \varphi =\left( \psi \,m_{i}^{n}\otimes
\varphi \right) \,\Delta ^{\left( n+1\right) }
\end{equation*}
restricted to some $\mathbf{C}^{\otimes R}$ ($R\in \mathbb{N}_{0}^{\times
n+1}$) is equal to 
\begin{eqnarray*}
\left( \psi \,m_{i}^{R}\otimes \varphi \right) \,\Delta _{R}=\left( \psi
\otimes \varphi \,\left( m_{i}^{R}\right) ^{-1}\right) \,\left(
m_{i}^{R}\otimes m_{i}^{R}\right) \,\Delta _{R} \\
&& \\
=\left( \left( \psi \otimes \gamma \right) \,\Delta _{D_{i}^{n+1}\left(
R\right) }\right) \,m_{i}^{R}=\left( \psi _{D_{i}^{n+1}\left( R\right) }\ast
\gamma \right) \,\,m_{i}^{R},
\end{eqnarray*}
where $\gamma \doteq \varphi \,\left( m_{i}^{R}\right) ^{-1}$ is a linear
form $\mathbf{C}^{\otimes D_{i}^{n+1}\left( R\right) }\rightarrow \Bbbk $.
This linear form is invertible w.r.t. the convolution product, and $\gamma
^{-1}\doteq \varphi ^{-1}\,\left( m_{i}^{R}\right) ^{-1}$, because (see Eqs. 
$\left( \ref{md}\right) $ and $\left( \ref{me}\right) $) 
\begin{equation*}
\gamma ^{-1}\ast \gamma =\left( \varphi ^{-1}\,\left( m_{i}^{R}\right)
^{-1}\otimes \varphi \,\left( m_{i}^{R}\right) ^{-1}\right) \,\Delta
_{D_{i}^{n+1}\left( R\right) }=\left( \varphi ^{-1}\otimes \varphi \right)
\,\Delta _{R}\,\left( m_{i}^{R}\right) ^{-1}=\varepsilon _{R}\,\left(
m_{i}^{R}\right) ^{-1}=\varepsilon _{D_{i}^{n+1}\left( R\right) }.
\end{equation*}
Therefore, $\gamma \in \mathsf{G}^{D_{i}^{n+1}\left( R\right) }$. Now, if $%
\psi \in \mathcal{Z}^{n}$, from the last lemma 
\begin{equation*}
\left( \psi _{D_{i}^{n+1}\left( R\right) }\ast \gamma \right)
\,\,m_{i}^{R}=\left( \gamma \ast \psi _{D_{i}^{n+1}\left( R\right) }\right)
\,\,m_{i}^{R}=\left( \varphi \otimes \psi \,m_{i}^{R}\right) \,\Delta _{R},
\end{equation*}
for all $R$, hence $\left( d_{i}^{n}\psi \right) \ast \varphi =\varphi \ast
\left( d_{i}^{n}\psi \right) $ for all $\varphi \in \mathsf{G}^{n+1}$. In
the $i=0$ case, 
\begin{equation*}
\left( \varepsilon \otimes \psi \right) \ast \varphi =\varphi \,\left[
\left( \varepsilon \otimes \psi \right) \ast \left( \mathbb{I}\otimes 
\mathbb{I}^{\otimes n}\right) \right] =\varphi \,\left[ \left( \varepsilon
\ast \mathbb{I}\right) \otimes \left( \psi \ast \mathbb{I}^{\otimes
n}\right) \right] =\varphi \,\left[ \mathbb{I}\otimes \left( \psi \ast 
\mathbb{I}^{\otimes n}\right) \right] ,
\end{equation*}
and similarly 
\begin{equation*}
\varphi \ast \left( \varepsilon \otimes \psi \right) =\varphi \,\left[
\left( \mathbb{I}\ast \varepsilon \right) \otimes \left( \mathbb{I}^{\otimes
n}\ast \psi \right) \right] =\varphi \,\left[ \mathbb{I}\otimes \left( 
\mathbb{I}^{\otimes n}\ast \psi \right) \right] .
\end{equation*}
In addition, we have supposed $\psi \in \mathcal{Z}^{n}$, then $\psi \ast 
\mathbb{I}^{\otimes n}=\mathbb{I}^{\otimes n}\ast \psi $ and $d_{0}^{n}\psi
\ast \varphi =\varphi \ast d_{0}^{n}\psi $ holds. And analogously for $i=n+1$%
. Hence, $d_{i}^{n}\left( \mathcal{Z}^{n}\right) \subset \mathcal{Z}^{n+1}$
for all $n,i$.\ 
\end{proof}

As a direct consequence of the last theorem, the cocycles,
(anti)bicharacters, coboundaries and every cohomologous class of $\mathsf{G}%
_{\mathcal{Z}}^{\bullet }$ can be seen as subsets of the cocycles, the
coboundaries and corresponding cohomologous class of $\mathsf{C}^{\bullet }$%
, respectively. Furthermore, the counital cochains of $\mathsf{G}_{\mathcal{Z%
}}^{\bullet }$ define another quasicomplex, namely $\frak{G}_{\mathcal{Z}%
}^{\bullet }$, such that there exists a monic $\frak{G}_{\mathcal{Z}%
}^{\bullet }\hookrightarrow \frak{C}^{\bullet }$ given by group
anti-monomorphisms.

\medskip

By last, consider a vector space $\mathbf{V}$ and the free bialgebra
generated by a \emph{multiplicative matrix} $t_{i}^{j}=v^{j}\otimes v_{i}\in 
\mathbf{V}^{\ast }\otimes \mathbf{V}$. That is to say, the algebra is $\left[
\mathbf{V}^{\ast }\otimes \mathbf{V}\right] ^{\otimes }$ and the coalgebra
structure is given by assignments $\Delta :t_{i}^{j}\mapsto t_{i}^{k}\otimes
t_{k}^{j}$ and $\varepsilon :t_{i}^{j}\mapsto \delta _{i}^{j}$. It can be
shown the bijection $End\left[ \mathbf{V}\right] =Lin\left[ \mathbf{V},%
\mathbf{V}\right] \backsimeq Lin\left[ \mathbf{V}^{\ast }\otimes \mathbf{V}%
,\Bbbk \right] $ gives rise to an isomorphism of quasicomplexes 
\begin{equation}
\mathsf{C}^{\bullet }\left[ \mathbf{V}\right] \backsimeq \mathsf{G}^{\bullet
}\left[ \left[ \mathbf{V}^{\ast }\otimes \mathbf{V}\right] ^{\otimes }\right]
.  \label{cg}
\end{equation}

\section{Twisting of quantum linear spaces}

As a previous step toward the quantum space twist transformations we need
the concept of \emph{admissible} cochain. Consider a pair $\mathcal{A}%
=\left( \mathbf{A}_{1},\mathbf{A}\right) $ and a 2-cochain $\psi \in \mathsf{%
C}^{2}\left[ \mathbf{A}_{1}\right] $. The map $\psi $ defines a linear map
in $\mathbf{A}\otimes \mathbf{A}$ \emph{iff} 
\begin{equation}
\psi \left( \mathbf{A}_{1}^{\otimes }\otimes \ker \Pi +\ker \Pi \otimes 
\mathbf{A}_{1}^{\otimes }\right) \subset \mathbf{A}_{1}^{\otimes }\otimes
\ker \Pi +\ker \Pi \otimes \mathbf{A}_{1}^{\otimes },  \label{nuc}
\end{equation}
where $\Pi $ is the canonical epimorphism $\mathbf{A}_{1}^{\otimes
}\twoheadrightarrow \mathbf{A}$. Such a map, namely $\Psi :\mathbf{A}\otimes 
\mathbf{A}\rightarrow \mathbf{A}\otimes \mathbf{A}$, would be given by the
equation $\Psi \,\Pi ^{\otimes 2}=\Pi ^{\otimes 2}\,\psi $.

\begin{definition}
We shall call $\mathcal{A}$\textbf{-admissible}, or \textbf{admissible} for $%
\mathcal{A}$, the $n$-cochains satisfying 
\begin{equation*}
\psi \left( \ker \Pi ^{\otimes n}\right) \subset \ker \Pi ^{\otimes
n}\;\;\;and\;\;\;\psi ^{-1}\left( \ker \Pi ^{\otimes n}\right) \subset \ker
\Pi ^{\otimes n},
\end{equation*}
or equivalently, $\psi \left( \ker \Pi ^{\otimes n}\right) =\ker \Pi
^{\otimes n}$.\ \ \ $\blacksquare $
\end{definition}

If $\mathcal{A}$ is conic, admissibility condition for $\psi \in \mathsf{C}%
^{2}\left[ \mathbf{A}_{1}\right] $ reduces to 
\begin{equation}
\psi \left( \mathbf{A}_{1}^{\otimes r}\otimes \mathbf{I}_{s}+\mathbf{I}%
_{r}\otimes \mathbf{A}_{1}^{\otimes s}\right) \subset \mathbf{A}%
_{1}^{\otimes r}\otimes \mathbf{I}_{s}+\mathbf{I}_{r}\otimes \mathbf{A}%
_{1}^{\otimes s},  \label{ap}
\end{equation}
for all $r,s\in \mathbb{N}_{0}$, where $\ker \Pi =\bigoplus_{n\geq 2}\mathbf{%
I}_{n}$. The other inclusions (and therefore the equalities) follow from the
fact that the involved vector spaces are finite dimensional.

The $\mathcal{A}$-admissible $n$-cochains form a subgroup of $\mathsf{C}^{n}%
\left[ \mathbf{A}_{1}\right] $. It can be shown the quotient of this
subgroup by\emph{\ }the identification 
\begin{equation}
\psi \backsim \psi ^{\prime }\;\;\Leftrightarrow \;\;\Pi ^{\otimes n}\left(
\psi -\psi ^{\prime }\right) =0,  \label{quo}
\end{equation}
is such that the map $\left[ \psi \right] \mapsto \Psi $, from admissible
equivalence classes to $Aut_{\mathrm{Vct}}\left[ \mathbf{A}^{\otimes n}%
\right] $, is an injective group homomorphism.

It is worth remarking admissibility of $\psi $ does not imply the one of $%
\partial \psi $. For instance, consider the algebra generated by $\mathbf{A}%
_{1}=span\left[ a,b\right] $ and quotient by the relation $ab=0$. Define the
1-cochain $\theta :\mathbf{A}_{1}^{\otimes }\backsimeq \mathbf{A}%
_{1}^{\otimes }$ such that $\theta _{n}$ is the identity unless for the
elements $abx$ and $xab$, where $\theta \left( abx\right) =xab$ and $\theta
\left( xab\right) =abx$, being $x$ equal to $a$ or $b$. Hence, $\theta $ is
clearly admissible. Then, for instance, 
\begin{equation*}
\delta _{1}\theta \,\left( ab\otimes b\right) =\left( m_{1}^{R}\right)
^{-1}\,\theta _{3}\,m_{1}^{R}\left( ab\otimes b\right) =\left(
m_{1}^{R}\right) ^{-1}\left( bab\right) =ba\otimes b,
\end{equation*}
with $R=\left( 2,1\right) $, what implies $\delta _{1}\theta \left( \ker \Pi
^{\otimes 2}\right) \neq \ker \Pi ^{\otimes 2}$, i.e. $\partial \theta $ is
not admissible. Therefore, unfortunately the subgroups of admissible
cochains do not define a subcomplex of $\mathsf{C}^{\bullet }\left[ \mathbf{A%
}_{1}\right] $. Nevertheless, in \S \textbf{4.1} we consider for each conic
quantum space $\mathcal{A}=\left( \mathbf{A}_{1},\mathbf{A}\right) $ a
stronger condition than admissibility, which enable us to define a
subcomplex $\mathsf{C}^{\bullet }\left[ \mathcal{A}\right] \subset \mathsf{C}%
^{\bullet }\left[ \mathbf{A}_{1}\right] $.

\subsection{Counital 2-cocycle twisting}

Consider a quantum space $\mathcal{A}=\left( \mathbf{A}_{1},\mathbf{A}%
\right) $ with related unital associative algebra structure $\left( m,\eta
\right) $, and an admissible counital 2-cochain $\psi \in \frak{C}^{2}\left[ 
\mathbf{A}_{1}\right] $. Since $\psi $ gives rise to a unique automorphism $%
\Psi $ of $\mathbf{A}\otimes \mathbf{A}$ by the equation $\Psi \,\Pi
^{\otimes 2}=\Pi ^{\otimes 2}\,\psi $, we can define, as at the beginning of
\S \textbf{2}, a new algebra structure over $\mathbf{A}$, namely $\mathbf{A}%
_{\psi }\doteq \left( \mathbf{A},m_{\psi }\right) $, with $m_{\psi }=m$%
\thinspace $\Psi $. In these terms we have the following notion of twisting
process on quantum spaces.

\begin{definition}
Given a quantum space $\mathcal{A}=\left( \mathbf{A}_{1},\mathbf{A}\right)
\in \mathrm{FGA}$ \emph{(}resp. $\mathrm{CA}$\emph{) }and a related
admissible $\psi \in \frak{Z}^{2}\left[ \mathbf{A}_{1}\right] $, we define
the \textbf{twisting of} $\mathcal{A}$ \textbf{by} $\psi $ as the pair $%
\mathcal{A}_{\psi }=\left( \mathbf{A}_{1},\mathbf{A}_{\psi }\right) \in 
\mathrm{FGA}$ \emph{(}resp. $\mathrm{CA}$\emph{)}.\ \ \ $\blacksquare $
\end{definition}

The following theorem justifies above definition.

\begin{theorem}
Let $\Psi $ be a map defined by an $\mathcal{A}$-\textbf{admissible counital
2-cocycle} $\psi \in \frak{Z}^{2}\left[ \mathbf{A}_{1}\right] $. Then $%
m_{\psi }=m$\thinspace $\Psi $ is an associative product with unit $\eta $,
and the associative unital algebra $\mathbf{A}_{\psi }\doteq \left( \mathbf{A%
},m_{\psi },\eta \right) $ is also generated by $\mathbf{A}_{1}$. Moreover, $%
m_{\psi }$ defines a filtration that coincides with the one associated to $m$
\emph{(}see Eq. $\left( \ref{filt}\right) $\emph{)}. In particular, for the
conic case, $m_{\psi }$ defines the same gradation as $m$.
\end{theorem}

\begin{proof}
Since $\Pi $ is a morphism of unital algebras $\Pi \,m_{\otimes }=m\,\Pi $
and $\Pi \left( \lambda \right) =\eta \left( \lambda \right) $, $\forall
\lambda \in \mathbf{A}_{1}^{\otimes 0}=\Bbbk $. Then, if $\psi $ is counital
we have 
\begin{equation*}
\Psi \left( \eta \left( \lambda \right) \otimes \left[ a\right] \right)
=\Psi \,\Pi ^{\otimes 2}\left( \lambda \otimes a\right) =\Pi ^{\otimes
2}\,\psi \left( \lambda \otimes a\right) =\eta \left( \lambda \right)
\otimes \left[ a\right] ,
\end{equation*}
being $\left[ a\right] =\Pi \left( a\right) $. So, $\Psi $\thinspace $\left(
\eta \otimes I\right) =\left( \eta \otimes I\right) $ and 
\begin{equation*}
m_{\psi }\,\left( \eta \otimes I\right) =m\,\Psi \left( \eta \otimes
I\right) \,=m\,\left( \eta \otimes I\right) =I.
\end{equation*}
The same holds for $I\otimes \eta $; hence $\eta $ is a unit map for the
product $m_{\psi }$. Let us see that $m_{\psi }$ is associative. It will be
associative \emph{iff} 
\begin{equation*}
m_{\psi }\,\left( m_{\psi }\otimes I\right) \,\Pi ^{\otimes 3}=m_{\psi
}\,\left( I\otimes m_{\psi }\right) \,\Pi ^{\otimes 3},
\end{equation*}
\emph{iff} 
\begin{equation*}
\Pi \,m_{\otimes }\,\psi \,\left( m_{\otimes }\otimes I\right) \,\left( \psi
\otimes I\right) =\Pi \,m_{\otimes }\,\psi \,\left( I\otimes m_{\otimes
}\right) \,\left( I\otimes \psi \right) .
\end{equation*}
From Equation $\left( \ref{im}\right) $ for $n=2$, $\psi $\thinspace $\left(
m_{\otimes }\otimes I\right) =\left( m_{\otimes }\otimes I\right) $%
\thinspace $\delta _{1}\psi $ and $\psi $\thinspace $\left( I\otimes
m_{\otimes }\right) =\left( I\otimes m_{\otimes }\right) $\thinspace $\delta
_{2}\psi $, so $m_{\psi }$ is associative \emph{iff} 
\begin{equation*}
\Pi \,m_{\otimes }\,\left( I\otimes m_{\otimes }\right) \,\left( \partial
\psi -\mathbb{I}^{\otimes 3}\right) =0,
\end{equation*}
where we have used the associativity of $m_{\otimes }$. But the last eq.
holds trivially since $\psi $ is a 2-cocycle.

It rests to show that $\mathbf{A}_{1}$ generates $\mathbf{A}$ through $%
m_{\psi }$. We know that, given a basis $\left\{ a_{i}\right\} $ of $\mathbf{%
A}_{1}$, it is possible to express any element of $\mathbf{A}$ as a linear
combination of words written with the \emph{letters} $a_{i}$'s and \emph{%
glued }with the product $m$. Let us show the same is true for $m_{\psi }$.
We will proceed by induction on the $\emph{length}$ of the words.

Denoting $m$ by a blank space and $m_{\psi }$ by $``\cdot "$, the map $\Psi $
for $r,s\in \mathbb{N}_{0}$ is given by a $\left( \dim \mathbf{A}_{1}\right)
^{r+s}$ square invertible matrix representing $\psi _{r,s}$, namely
(remember that $\Pi \left( a_{i}\right) =a_{i}$) 
\begin{equation*}
\Psi \left( a_{i_{1}}...a_{i_{r}}\otimes a_{j_{1}}...a_{j_{s}}\right) =\psi
_{i_{1}...i_{r},j_{1}...,j_{s}}^{k_{1}...k_{r},l_{1}...l_{s}}%
\;a_{k_{1}}...a_{k_{r}}\otimes a_{l_{1}}...a_{l_{s}}.
\end{equation*}
Since 
\begin{equation*}
a_{i}\cdot a_{j}=m_{\psi }\left( a_{i}\otimes a_{j}\right) =\psi
_{i,j}^{k,l}\;m\left( a_{k}\otimes a_{l}\right) =\psi
_{i,j}^{k,l}\;a_{k}\,a_{l},
\end{equation*}
we have 
\begin{equation*}
a_{k}\,a_{l}=\left( \psi ^{-1}\right) _{k,l}^{i,j}\;a_{i}\cdot a_{j}.
\end{equation*}
Then, for words of length $2$ there exists an invertible linear map $\theta
_{2}:\mathbf{A}_{1}^{\otimes 2}\backsimeq \mathbf{A}_{1}^{\otimes 2}$, with
matrix coefficients $\left( \theta _{2}\right) _{ij}^{kl}=\psi _{i,j}^{k,l}$%
, such that 
\begin{equation*}
a_{k}\,a_{l}=\left( \theta _{2}^{-1}\right) _{ij}^{kl}\;a_{i}\cdot a_{j}.
\end{equation*}
Suppose that for words of $n$ letters there exists an invertible linear map $%
\theta _{n}$ such that 
\begin{equation}
a_{i_{1}}...a_{i_{n}}=\left( \theta _{n}^{-1}\right)
_{i_{1}...i_{n}}^{j_{1}...j_{n}}\;a_{j_{1}}\cdot ...\cdot a_{j_{n}}.
\label{pitita}
\end{equation}
Then, for $n+1$%
\begin{equation*}
\begin{array}{l}
\left( a_{j_{1}}\cdot ...\cdot a_{j_{n}}\right) \cdot a_{j_{n+1}}=\left(
\theta _{n}\right) _{j_{1}...j_{n}}^{i_{1}...i_{n}}\;\left(
a_{i_{1}}...a_{i_{n}}\right) \cdot a_{j_{n+1}} \\ 
\\ 
=\left( \theta _{n}\right) _{j_{1}...j_{n}}^{i_{1}...i_{n}}\;\psi
_{i_{1}...i_{n},j_{n+1}}^{k_{1}...k_{r},k_{n+1}}\;a_{k_{1}}...a_{k_{n+1}}%
\doteq \left( \theta _{n+1}\right)
_{j_{1}...j_{n+1}}^{k_{1}...k_{n+1}}\;a_{k_{1}}...a_{k_{n+1}},
\end{array}
\end{equation*}
i.e. 
\begin{equation*}
a_{k_{1}}...a_{k_{n+1}}=\left( \theta _{n+1}^{-1}\right)
_{k_{1}...k_{n+1}}^{j_{1}...j_{n+1}}\;a_{j_{1}}\cdot ...\cdot a_{j_{n+1}},
\end{equation*}
being 
\begin{equation}
\left( \theta _{n+1}^{-1}\right) _{k_{1}...k_{n+1}}^{j_{1}...j_{n+1}}\doteq
\left( \psi ^{-1}\right)
_{k_{1}...k_{n},k_{n+1}}^{i_{1}...i_{n},j_{n+1}}\;\left( \theta
_{n}^{-1}\right) _{i_{1}...i_{n}}^{j_{1}...j_{n}}.  \label{ind}
\end{equation}
Hence, by induction, any element of $\mathbf{A}$ can be expressed as linear
combinations of words constructed with elements of $\mathbf{A}_{1}$ and
glued with $m_{\psi }$. Both words have the same length, in the sense that
can be obtained as the image by $\Pi $ of homogeneous elements with same
number of factors. From that it is clear $m_{\psi }$ defines the same
filtration as $m$, and for the conic case the same gradation, as we wanted
to show.
\end{proof}

In other words, the last theorem says the pair $\left( \mathbf{A}_{1},%
\mathbf{A}_{\psi }\right) $ is a quantum space, and if $\left( \mathbf{A}%
_{1},\mathbf{A}\right) $ is conic, so is its twisting. In the following, we
shall describe some interesting examples.

\subsubsection{The quantum plane as a 2-cocycle twisting}

Let $\mathcal{A}^{2|0}=\left( \mathbf{A}_{1},\mathbf{A}\right) $ be the
quantum space with $\mathbf{A}_{1}=span\left[ a,b\right] $ and $\mathbf{A}%
=\Bbbk \left[ a,b\right] $, the commutative algebra freely generated by the
symbols $a$ and $b$. In other words, $\mathcal{A}^{2|0}$ gives the
coordinate ring of the two dimensional affine space. Consider the basis 
\begin{equation*}
\left\{ a^{f_{1}}b^{g_{1}}...a^{f_{k}}b^{g_{k}}:\;f,g\in \mathbb{N}%
_{0},\;k\in \mathbb{N}\right\}
\end{equation*}
of $\mathbf{A}_{1}^{\otimes }$, and for every element $v$ in that basis the
positive numbers $m^{a}\left( v\right) \doteq \sum_{i=1}^{k}\,f_{i}$ and $%
m^{b}\left( v\right) \doteq \sum_{i=1}^{k}\,g_{i}$. Now, given $\hslash \in
\Bbbk $ let us define the 2-cochain 
\begin{equation*}
\psi _{\hslash }:v\otimes w\mapsto \exp \left[ m^{a}\left( v\right)
\,m^{b}\left( w\right) \,\hslash \right] \,v\otimes w.
\end{equation*}
The counitality is immediate, since taking $v=w=1$, equation $m^{a}\left(
v\right) \,m^{b}\left( w\right) =0$ holds. This cochain is $\mathcal{A}$%
-admissible, in fact, $\psi _{\hslash }$ gives rise to a linear automorphism 
$\Psi _{\hslash }$ of $\mathbf{A}^{\otimes 2}$ such that 
\begin{equation}
\Psi _{\hslash }\left( a^{n}b^{m}\otimes a^{r}b^{s}\right) =\exp \left(
ns\hslash \right) \;a^{n}b^{m}\otimes a^{r}b^{s}.  \label{fig}
\end{equation}
The formula above defines $\Psi _{\hslash }$ completely, because the set $%
\left\{ a^{n}b^{m}\right\} _{n,m\in \mathbb{N}_{0}}$ is a basis for $\mathbf{%
A}$. Let us see that $\psi _{\hslash }$ is a 2-cocycle. Straightforward
calculations show that 
\begin{equation*}
\begin{array}{l}
\delta _{1}\psi _{\hslash }\left( u\otimes v\otimes w\right) =\exp \left[
\left( m^{a}\left( u\right) +m^{a}\left( v\right) \right) \,m^{b}\left(
w\right) \,\hslash \right] \,u\otimes v\otimes w, \\ 
\\ 
\delta _{2}\psi _{\hslash }\left( u\otimes v\otimes w\right) =\exp \left[
m^{a}\left( u\right) \,\left( m^{b}\left( v\right) +m^{b}\left( w\right)
\right) \,\hslash \right] \,u\otimes v\otimes w,
\end{array}
\end{equation*}
for any elements $u$, $v$ and $w$ of the basis of $\mathbf{A}_{1}^{\otimes }$%
. Now, applying consecutively $\psi _{\hslash }\otimes \mathbb{I}$ and $%
\,\delta _{1}\psi _{\hslash }$ on such elements we have 
\begin{equation*}
\begin{array}{l}
u\otimes v\otimes w\mapsto \exp \left[ m^{a}\left( u\right) \,m^{b}\left(
v\right) \,\hslash \right] \,u\otimes v\otimes w \\ 
\\ 
\mapsto \exp \left[ \left( m^{a}\left( u\right) +m^{a}\left( v\right)
\right) \,m^{b}\left( w\right) \,\hslash \right] \exp \left[ m^{a}\left(
u\right) \,m^{b}\left( v\right) \,\hslash \right] \,u\otimes v\otimes w \\ 
\\ 
=\exp \left[ \left( m^{a}\left( u\right) \,m^{b}\left( w\right) +m^{a}\left(
v\right) \,m^{b}\left( w\right) +m^{a}\left( u\right) \,m^{b}\left( v\right)
\right) \,\hslash \right] \,u\otimes v\otimes w,
\end{array}
\end{equation*}
and doing the same for $\,\mathbb{I}\otimes \psi _{\hslash }$ and $\delta
_{2}\psi _{\hslash }$, 
\begin{equation*}
\begin{array}{l}
u\otimes v\otimes w\mapsto \exp \left[ m^{a}\left( v\right) \,m^{b}\left(
w\right) \,\hslash \right] \,u\otimes v\otimes w \\ 
\\ 
\mapsto \exp \left[ m^{a}\left( u\right) \,\left( m^{b}\left( v\right)
+m^{b}\left( w\right) \right) \,\hslash \right] \exp \left[ m^{a}\left(
v\right) \,m^{b}\left( w\right) \,\hslash \right] \,u\otimes v\otimes w \\ 
\\ 
=\exp \left[ \left( m^{a}\left( u\right) \,m^{b}\left( v\right) +m^{a}\left(
u\right) \,m^{b}\left( w\right) +m^{a}\left( v\right) \,m^{b}\left( w\right)
\right) \,\hslash \right] \,u\otimes v\otimes w.
\end{array}
\end{equation*}
Thus $\delta _{1}\psi _{\hslash }\,\left( \psi _{\hslash }\otimes \mathbb{I}%
\right) =\delta _{2}\psi _{\hslash }\,\left( \mathbb{I}\otimes \psi
_{\hslash }\right) $. Moreover, it can be seen that $\psi _{\hslash }$ is a
bicharacter and also an anti-bicharacter. Accordingly, the product $m_{\psi
_{\hslash }}=m\,\Psi _{\hslash }=\star $ is associative, and from Eq. $%
\left( \ref{fig}\right) $ for $n=s=0$ and $m=r=1$, $a\star b=\exp \hslash
\;b\star a$ follows. Thus, $\mathcal{A}_{\psi _{\hslash }}^{2|0}$ is
isomorphic to the quantum plane 
\begin{equation*}
\mathcal{A}_{\hslash }^{2|0}\doteq \left. \Bbbk \left\langle
x,y\right\rangle \right/ \left\{ xy-\exp \hslash \;yx\right\}
\end{equation*}
under the extension of $x\mapsto a$, $y\mapsto b$ to an algebra map. The map 
$\psi _{\hslash }$ also defines an admissible 2-cocycle for the superplane $%
\mathcal{A}^{0|2}$, given by a Grassmann algebra in two variables. Writing
again $m_{\psi _{\hslash }}=\star $, we have the relation $a\star b=-\exp
\hslash \;b\star a$, thus $\mathcal{A}_{\psi _{\hslash }}^{0|2}$ is
isomorphic to the quantum $0|2$-dimensional superplane 
\begin{equation*}
\mathcal{A}_{-\hslash }^{0|2}=\left. \Bbbk \left\langle \eta ,\xi
\right\rangle \right/ \left\{ \eta \xi +\exp \hslash \;\xi \eta ;\;\eta
^{2};\;\xi ^{2}\right\} .
\end{equation*}
This result can easily be generalized to any quantum space $\mathcal{A}%
_{\hslash }^{n|0}$ and $\mathcal{A}_{\hslash }^{0|n}$. In general, given $%
M\in GL\left( n^{2}\right) $ fulfilling the YB equation or such that $%
M_{12}\,M_{23}=M_{23}\,M_{12}$, we can define, as at the end of \S \textbf{%
2.2.2}, a bicharacter or an anti-bicharacter $\psi $ related to an $n$%
-dimensional vector space $\mathbf{A}_{1}$. Such a cochain is always
admissible for commutative and anticommutative quantum spaces $\mathcal{A}$
generated by $\mathbf{A}_{1}$. The ideal related to $\mathcal{A}_{\psi }$ is
generated by 
\begin{equation}
\left( N_{ij}^{kl}\mp N_{ji}^{kl}\right) \,a_{k}\otimes a_{l};\;\;\;N=M^{-1},
\label{dr}
\end{equation}
the $\mp $ corresponding to the $\mathcal{A}$ commutative or $\mathcal{A}$
anticommutative cases, respectively. (We need the results of \S \textbf{%
3.2.2 }to show that.) Note relations in $\left( \ref{dr}\right) $ can be
written $\left( R_{ij}^{kl}\mp P_{ij}^{kl}\right) \,a_{k}\otimes a_{l}$,
being $R=M_{21}N$ and $P$ the permutation matrix. It can be shown $R$ is a
YB operator, in fact is a \emph{triangular} one, i.e. $R_{21}R=I$. Thus,
twisting $\mathcal{A}$ by (anti)bicharacters define quantum spaces with
relations given by triangular YB operators of the form $M_{21}M^{-1}$. As a
particular case, to obtain $\mathcal{A}_{\hslash }^{n|0}$ or $\mathcal{A}%
_{-\hslash }^{0|n}$ from $\mathcal{A}$, $\psi $ can be defined by 
\begin{equation*}
M_{ij}^{kl}=\exp \left\{ \left. \left[ 1-\func{sg}\left( i-j\right) \right]
\,\hslash \right/ 2\right\} \,\delta _{i}^{k}\delta
_{j}^{l};\;\;i,j,k,l=1...n;
\end{equation*}
being $\func{sg}\left( 0\right) =1$. The extension of $M$ to a bicharacter
and to an anti-bicharacter gives rise to the same map. Since $M$ is a YB
operator, the resulting cochain $\psi $ is effectively a 2-cocycle.

\subsubsection{Symmetric twisted tensor products}

Given a couple of quantum spaces $\mathcal{A}$ and $\mathcal{B}$, we recall
that a symmetric twisted tensor product (STTP) of them is (essentially) a
quantum space $\mathcal{A}\circ _{\tau }\mathcal{B}=\left( \mathbf{A}%
_{1}\otimes \mathbf{B}_{1},\mathbf{A}\circ _{\tau }\mathbf{B}\right) $,
being $\mathbf{A}\circ _{\tau }\mathbf{B}$ the subalgebra of $\mathbf{A}%
\otimes _{\tau }\mathbf{B}$ \cite{cap} generated by $\mathbf{A}_{1}\otimes 
\mathbf{B}_{1}$, and where the related (symmetric) twisting map $\tau $ is a
linear bijection $\mathbf{B}\otimes \mathbf{A\backsimeq A}\otimes \mathbf{B}$
that defines, by restriction, an isomorphism $\mathbf{B}_{1}\otimes \mathbf{A%
}_{1}\backsimeq \mathbf{A}_{1}\otimes \mathbf{B}_{1}$. Thus, $\tau $ is
completely defined by the last isomorphism and by the properties of a
twisting map, namely 
\begin{equation}
\begin{array}{l}
\tau \left( m_{B}\otimes I_{A}\right) =\left( I_{A}\otimes m_{B}\right)
\,\left( \tau \otimes I_{B}\right) \,\left( I_{B}\otimes \tau \right) , \\ 
\tau \left( I_{B}\otimes m_{A}\right) =\left( m_{A}\otimes I_{B}\right)
\,\left( I_{A}\otimes \tau \right) \,\left( \tau \otimes I_{A}\right) ,
\end{array}
\label{mu}
\end{equation}
and 
\begin{equation}
\tau \left( I_{B}\otimes \eta _{A}\right) =\eta _{A}\otimes I_{B};\;\;\tau
\left( \eta _{B}\otimes I_{A}\right) =I_{A}\otimes \eta _{B}.  \label{un}
\end{equation}
Remember that the algebra structure of $\mathbf{A}\otimes _{\tau }\mathbf{B}$
is given by the maps 
\begin{equation*}
m_{A\otimes B}^{\tau }=\left( m_{A}\otimes m_{B}\right) \,\left(
I_{A}\otimes \tau \otimes I_{B}\right) \;\;and\;\;\eta _{A\otimes B}=\eta
_{A}\otimes \eta _{B}.
\end{equation*}
We will show that any STTP can be seen as a twisting of $\mathcal{A}\circ 
\mathcal{B}$ by an element of $\frak{Z}^{2}\left[ \mathbf{A}_{1}\otimes 
\mathbf{B}_{1}\right] $.

\bigskip

The usual product in $\mathbf{A}\otimes \mathbf{B}$ is $m_{\mathbf{A}\otimes 
\mathbf{B}}=\left( m_{\mathbf{A}}\otimes m_{\mathbf{B}}\right) \,\left(
I\otimes f\otimes I\right) $, being $f$ the canonical flipping map. Then,
given $\tau $ we have that $m_{\mathbf{A}\otimes \mathbf{B}}^{\tau }=m_{%
\mathbf{A}\otimes \mathbf{B}}\,\Omega $, where 
\begin{equation*}
\Omega \doteq \left( I\otimes f^{-1}\,\tau \otimes I\right) :\left( \mathbf{A%
}\otimes \mathbf{B}\right) ^{\otimes 2}\rightarrow \left( \mathbf{A}\otimes 
\mathbf{B}\right) ^{\otimes 2}.
\end{equation*}
And since $\Omega \left( \left( \mathbf{A}\circ \mathbf{B}\right) ^{\otimes
2}\right) \subset \left( \mathbf{A}\circ \mathbf{B}\right) ^{\otimes 2}$
(which follows from property $\left( \ref{mu}\right) $), we also have $m_{%
\mathbf{A}\circ \mathbf{B}}^{\tau }=m_{\mathbf{A}\circ \mathbf{B}}\,\Omega $%
, by making the corresponding restriction. Now, let us extend the
isomorphism 
\begin{equation*}
\left. \tau \right| _{\mathbf{B}_{1}\otimes \mathbf{A}_{1}}:\mathbf{B}%
_{1}\otimes \mathbf{A}_{1}\backsimeq \mathbf{A}_{1}\otimes \mathbf{B}_{1}
\end{equation*}
to all of $\mathbf{B}_{1}^{\otimes }\otimes \mathbf{A}_{1}^{\otimes }$ using
Eq. $\left( \ref{mu}\right) $ and $\left( \ref{un}\right) $, and call $\tau
_{\otimes }$ the resulting map. By construction, $\tau _{\otimes }$ defines
isomorphisms 
\begin{equation}
\mathbf{B}_{1}^{\otimes r}\otimes \mathbf{A}_{1}^{\otimes s}\backsimeq 
\mathbf{A}_{1}^{\otimes s}\otimes \mathbf{B}_{1}^{\otimes r},\;\;\forall
r,s\in \mathbb{N}_{0}.  \label{if}
\end{equation}
Then, using the canonical flipping map $f_{\otimes }$ on $\mathbf{B}%
_{1}^{\otimes }\otimes \mathbf{A}_{1}^{\otimes }$ (which obviously defines
similar isomorphisms to the ones given in equation above) we have a linear
automorphism 
\begin{equation*}
\widehat{\omega }\doteq \left( \mathbb{I}\otimes f_{\otimes }^{-1}\,\tau
_{\otimes }\otimes \mathbb{I}\right) :\left( \mathbf{A}_{1}^{\otimes
}\otimes \mathbf{B}_{1}^{\otimes }\right) ^{\otimes 2}\rightarrow \left( 
\mathbf{A}_{1}^{\otimes }\otimes \mathbf{B}_{1}^{\otimes }\right) ^{\otimes
2}
\end{equation*}
that is homogeneous (since Eq. $\left( \ref{if}\right) $) and that
restricted to 
\begin{equation*}
\left( \mathbf{A}_{1}^{\otimes }\circ \mathbf{B}_{1}^{\otimes }\right)
^{\otimes 2}\subset \left( \mathbf{A}_{1}^{\otimes }\otimes \mathbf{B}%
_{1}^{\otimes }\right) ^{\otimes 2}
\end{equation*}
satisfies $\Omega \,\Pi _{\circ }^{\otimes 2}=\Pi _{\circ }^{\otimes 2}\,%
\widehat{\omega }$. (Recall the identification $\mathbf{A}_{1}^{\otimes
}\circ \mathbf{B}_{1}^{\otimes }\backsimeq \left[ \mathbf{A}_{1}\otimes 
\mathbf{B}_{1}\right] ^{\otimes }$.) Thus, $\widehat{\omega }$ gives rise
(by restriction) to a 2-cochain in $\mathsf{C}^{2}\left[ \mathbf{A}%
_{1}\otimes \mathbf{B}_{1}\right] $, which is automatically $\mathcal{A}%
\circ \mathcal{B}$-admissible. Let us call $\omega $ such a cochain. Because 
$\tau _{\otimes }$ satisfies Eq. $\left( \ref{un}\right) $, $\omega $ is
counital, i.e. $\omega \in \frak{C}^{2}\left[ \mathbf{A}_{1}\otimes \mathbf{B%
}_{1}\right] $. It remains to see that $\omega $ is a 2-cocycle in order to
prove our claim. Eq. $\left( \ref{mu}\right) $ for $\omega $ translates into
equations $\delta _{1}\omega =\omega _{13}\,\omega _{23}$ and $\delta
_{2}\omega =\omega _{13}\,\omega _{12}$, with $\omega _{12}=\omega \otimes 
\mathbb{I}$, $\omega _{23}=\mathbb{I}\otimes \omega $ and $\omega
_{13}=\left( \mathbb{I}\otimes f_{\otimes }^{-1}\right) \,\omega
_{12}\,\left( \mathbb{I}\otimes f_{\otimes }\right) $. That is to say, $%
\omega $ is an anti-bicharacter. Then, since clearly $\omega _{23}\,\omega
_{12}=\omega _{12}\,\omega _{23}$, $\omega $ is effectively a 2-cocycle;
explicitly 
\begin{equation*}
\delta _{1}\omega \,\left( \omega \otimes \mathbb{I}\right) =\omega
_{13}\,\omega _{23}\,\omega _{12}=\omega _{13}\,\omega _{12}\,\omega
_{23}=\delta _{2}\omega \,\left( \mathbb{I}\otimes \omega \right) .
\end{equation*}

Let us resume these results in one proposition.

\begin{proposition}
Given an STTP $\mathcal{A}\circ _{\tau }\mathcal{B}$ of a couple $\mathcal{A}
$ and $\mathcal{B}$ of quantum spaces, there exists an admissible counital
2-cocycle $\omega $ for $\mathcal{A}\circ \mathcal{B}$ such that $\mathcal{A}%
\circ _{\tau }\mathcal{B}=\left( \mathcal{A}\circ \mathcal{B}\right)
_{\omega }$. Moreover, $\omega $ is an anti-bicharacter.\ \ \ $\blacksquare $
\end{proposition}

\subsubsection{Twisting of bialgebras}

Suppose a quantum space $\mathcal{A}$ is a bialgebra in \textrm{FGA }(resp. $%
\mathrm{CA}$), i.e. the coalgebra structure is given by arrows $\Delta :%
\mathcal{A}\rightarrow \mathcal{A}\circ \mathcal{A}$ and $\varepsilon :%
\mathcal{A}\rightarrow \mathcal{I}$ (resp. $\varepsilon :\mathcal{A}%
\rightarrow $ $\mathcal{K}$) there. In particular, $\mathbf{A}_{1}$ is a
coalgebra with coproduct $\Delta _{1}=\left. \Delta \right| _{\mathbf{A}%
_{1}} $ and counit $\varepsilon _{1}=\left. \varepsilon \right| _{\mathbf{A}%
_{1}}$. That is to say, a bialgebra in \textrm{FGA }is a bialgebra
algebraically generated by a finite dimensional coalgebra. The pair $\left(
\Delta _{1}^{\otimes },\varepsilon _{1}^{\otimes }\right) $ defines a
bialgebra structure on $\mathbf{A}_{1}^{\otimes }$, and the canonical
epimorphism $\Pi $ is a bialgebra homomorphism w.r.t. this structure.

Let $\chi :\mathbf{A}^{\otimes 2}\rightarrow \Bbbk $ be a 2-cocycle of $%
\mathsf{G}^{2}\left[ \mathbf{A}\right] $ (see \S \textbf{2.4}). From the
morphism of quasicomplexes $\mathsf{G}^{\bullet }\left[ \mathbf{A}\right]
\rightarrow \mathsf{G}^{\bullet }\left[ \mathbf{A}_{1}^{\otimes }\right] $
induced by the bialgebra map $\Pi :\mathbf{A}_{1}^{\otimes }\rightarrow 
\mathbf{A}$ (see Eq. $\left( \ref{nt}\right) $ in \textbf{Prop. 7}), the map 
$\chi \,\Pi ^{\otimes 2}:$ $\mathbf{A}_{1}^{\otimes }\otimes \mathbf{A}%
_{1}^{\otimes }\rightarrow \Bbbk $ is also a 2-cocycle. Since $\Pi $ is
surjective, the maps $\chi \mapsto \chi \,\Pi ^{\otimes n}$, $n\in \mathbb{N}%
_{0}$, give rise to a monic $\mathsf{G}^{\bullet }\left[ \mathbf{A}\right]
\hookrightarrow \mathsf{G}^{\bullet }\left[ \mathbf{A}_{1}^{\otimes }\right] 
$. Now, from $\digamma ^{\bullet }:\mathsf{G}^{\bullet }\left[ \mathbf{A}%
_{1}^{\otimes }\right] \rightarrow \mathsf{C}^{\bullet }\left[ \mathbf{A}_{1}%
\right] $, we have in addition the automorphism 
\begin{equation*}
\digamma _{\chi }\doteq \digamma ^{2}\left( \chi \,\Pi ^{\otimes 2}\right)
=\left( \chi \,\Pi ^{\otimes 2}\right) \ast \mathbb{I}^{\otimes 2}\ast
\left( \chi \,\Pi ^{\otimes 2}\right) ^{-1},
\end{equation*}
which, by \textbf{Theor. 6}, defines a 2-cocycle in $\mathsf{C}^{2}\left[ 
\mathbf{A}_{1}\right] $. The composition of $\digamma ^{\bullet }$ with the
above monic, quotient out the corresponding centers $\mathcal{Z}$, is
another monic $\mathsf{G}_{\mathcal{Z}}^{\bullet }\left[ \mathbf{A}\right]
\hookrightarrow \mathsf{C}^{\bullet }\left[ \mathbf{A}_{1}\right] $ (see 
\textbf{Theor. 7}), which enable us to think of $\mathsf{G}^{\bullet }\left[ 
\mathbf{A}\right] $ as a subcomplex of $\mathsf{C}^{\bullet }\left[ \mathbf{A%
}_{1}\right] $. Of course, all that is also true if $\mathsf{G}^{\bullet }$
and $\mathsf{C}^{\bullet }$ are replaced by their counital subgroups $\frak{G%
}^{\bullet }$ and $\frak{C}^{\bullet }$.

We recall that the twisting \cite{drin}\cite{maj} of a bialgebra $\mathbf{A}$
by $\chi $ is the bialgebra $\mathbf{A}_{\chi }$ with the same coalgebra
structure, and associative product $m_{\chi }=\chi \ast m\ast \chi ^{-1}$
with the same unit.\footnote{%
If a couple of elements $\chi ,\vartheta \in \mathsf{G}^{2}$ `differ' by an
element of the center $\mathcal{Z}$, both linear forms define the same
twisted bialgebra. This is why the interesting objects are the qoutients $%
\mathsf{G}_{\mathcal{Z}}^{n}$ instead of entire groups $\mathsf{G}^{n}$.}

$\digamma _{\chi }$ is admissible, in fact 
\begin{equation}
\Pi ^{\otimes 2}\,\digamma _{\chi }=\Pi ^{\otimes 2}\,\digamma ^{2}\left(
\chi \,\Pi ^{\otimes 2}\right) =\left( \chi \ast \mathbb{I}^{\otimes 2}\ast
\chi ^{-1}\right) \,\Pi ^{\otimes 2}.\,  \label{addm}
\end{equation}
Hence, given a counital 2-cocycle $\chi \in \mathsf{G}_{\mathcal{Z}%
}^{\bullet }\left[ \mathbf{A}\right] $, we have an admissible cocycle $%
\digamma _{\chi }$ in $\frak{Z}^{2}\left[ \mathbf{A}_{1}\right] $. Moreover,
the twisting $\mathcal{A}_{\chi }=\left( \mathbf{A}_{1},\mathbf{A}_{\chi
}\right) $ by $\chi \in \frak{G}^{2}\left[ \mathbf{A}\right] $ of the
bialgebra $\mathcal{A}$, coincides with the twisting $\mathcal{A}_{\digamma
_{\chi }}$ of $\mathcal{A}$ as a quantum space, because from Eq. $\left( \ref
{addm}\right) $, $m_{\digamma _{\chi }}=\chi \ast m\ast \chi ^{-1}=m_{\chi }$%
.

Let us consider the twisted internal coEnd objects $\underline{end}%
^{\Upsilon }\left[ \mathcal{A}\right] $, related to a symmetric twisting map
defined by $\widehat{\tau }_{\mathcal{A}}=id\otimes \sigma ^{!}\otimes
\sigma $, where $\sigma =\sigma _{\mathcal{A}}:\mathbf{A}_{1}\backsimeq 
\mathbf{A}_{1}$ and $\sigma ^{!}=\sigma ^{\ast -1}$, such that $\sigma _{%
\mathcal{A}}$ can be extended to quantum space automorphism $\mathcal{A}%
\backsimeq \mathcal{A}$. We mentioned in \S \textbf{1.2} that each one of
these bialgebras are isomorphic to the twisting $\left( \mathcal{A}%
\triangleright \mathcal{A}\right) _{\chi }$ (as bialgebras). Such $\chi $ is
a 2-cocycle in $\frak{G}^{2}\left[ \mathbf{A}\triangleright \mathbf{A}\right]
$ defined by 
\begin{equation*}
\chi \left( z_{k_{1}}^{j_{1}}...z_{k_{r}}^{j_{r}}\otimes
z_{k_{r+1}}^{j_{r+1}}...z_{k_{r+s}}^{j_{r+s}}\right) =\delta
_{k_{1}}^{j_{1}}...\delta _{k_{r}}^{j_{r}}\left( \sigma ^{-r}\right)
_{k_{r+1}}^{j_{r+1}}...\left( \sigma ^{-r}\right) _{k_{r+s}}^{j_{r+s}}\in
\Bbbk ,
\end{equation*}
being $z_{i}^{j}=a^{j}\otimes a_{i}$. Then, \emph{via} $\digamma $, $%
\underline{end}^{\Upsilon }\left[ \mathcal{A}\right] $ can be seen as a
twisting of $\underline{end}\left[ \mathcal{A}\right] $ by $\digamma _{\chi
}\in \frak{Z}^{2}\left[ \mathbf{A}_{1}^{\ast }\otimes \mathbf{A}_{1}\right] $%
, where 
\begin{equation}
\begin{array}{l}
\digamma _{\chi }\left( z_{k_{1}}^{j_{1}}...z_{k_{r}}^{j_{r}}\otimes
z_{k_{r+1}}^{j_{r+1}}...z_{k_{r+s}}^{j_{r+s}}\right) = \\ 
\\ 
=z_{k_{1}}^{j_{1}}...z_{k_{r}}^{j_{r}}\otimes \left( \left( \sigma
^{-r}\right) _{k_{r+1}}^{a_{1}}...\left( \sigma ^{-r}\right)
_{k_{r+s}}^{a_{s}}\right) \,z_{a_{1}}^{b_{1}}...z_{a_{s}}^{b_{s}}\,\left(
\left( \sigma ^{r}\right) _{b_{1}}^{j_{r+1}}...\left( \sigma ^{r}\right)
_{b_{s}}^{j_{r+s}}\right) .
\end{array}
\label{fgy}
\end{equation}
One can show $\chi $, and consequently $\digamma _{\chi }$, is bicharacter
and anti-bicharacter.

In the following we generalize this correspondence between $\underline{end}%
^{\Upsilon }\left[ \mathcal{A}\right] $ and $\underline{end}\left[ \mathcal{A%
}\right] $ to all of twisted internal coHom objects.

\subsubsection{The twisted coHom objects}

Consider the twisted coHom objects $\underline{hom}^{\Upsilon }\left[ 
\mathcal{B},\mathcal{A}\right] $ related to a twisting map defined by 
\begin{equation}
\widehat{\tau }_{\mathcal{A},\mathcal{B}}=id\otimes \rho \otimes \phi
,\;\phi =\sigma _{\mathcal{A}},\;\rho ^{-1}=\sigma _{\mathcal{B}}^{\ast },
\label{tau}
\end{equation}
where $\sigma _{\mathcal{A}}$ and $\sigma _{\mathcal{B}}$ can be extended to
automorphisms $\mathcal{A}\backsimeq \mathcal{A}$ and $\mathcal{B}\backsimeq 
\mathcal{B}$. Then, we know that $\underline{hom}^{\Upsilon }\left[ \mathcal{%
B},\mathcal{A}\right] $ is the algebra generated by $z_{i}^{j}=b^{j}\otimes
a_{i}\in \mathbf{B}_{1}^{\ast }\otimes \mathbf{A}_{1}$ and quotient by the
ideal algebraically generated by (see Eqs. $\left( \ref{I}\right) $ to $%
\left( \ref{VI}\right) $) 
\begin{equation*}
\left\{ \left\{ _{\sigma }R_{\lambda
_{n}}^{k_{1}...k_{n}}\;z_{k_{1}}^{j_{1}}\cdot
z_{k_{2}}^{j_{2}}\;...\;z_{k_{n}}^{j_{n}}\;\left( _{\sigma }S^{\bot }\right)
_{j_{1}...j_{n}}^{\omega _{n}}\right\} _{\omega _{n}\in \Omega
_{n}}^{\lambda _{n}\in \Lambda _{n}}\right\} _{n\in \mathbb{N}_{0}}.
\end{equation*}
Now, let us define for $\underline{hom}\left[ \mathcal{B},\mathcal{A}\right]
=\mathcal{B}\triangleright \mathcal{A}$ the counital 2-cochain $\varsigma
\in \frak{C}^{2}\left[ \mathbf{B}_{1}^{\ast }\otimes \mathbf{A}_{1}\right] $
given by 
\begin{equation}
\varsigma \left( z_{k_{1}}^{j_{1}}...z_{k_{r}}^{j_{r}}\otimes
z_{k_{r+1}}^{j_{r+1}}...z_{k_{r+s}}^{j_{r+s}}\right)
=z_{k_{1}}^{j_{1}}...z_{k_{r}}^{j_{r}}\otimes \xi ^{\left[ r\right] }\left(
z_{k_{r+1}}^{j_{r+1}}...z_{k_{r+s}}^{j_{r+s}}\right) ,  \label{sed}
\end{equation}
being $\xi ^{\left[ r\right] }\in \frak{C}^{1}\left[ \mathbf{B}_{1}^{\ast
}\otimes \mathbf{A}_{1}\right] $ such that 
\begin{equation}
\xi ^{\left[ r\right] }\left( z_{k_{1}}^{j_{1}}...z_{k_{s}}^{j_{s}}\right)
=\left( \left( \phi ^{-r}\right) _{k_{1}}^{a_{1}}...\left( \phi ^{-r}\right)
_{k_{s}}^{a_{s}}\right) \,z_{a_{1}}^{b_{1}}...z_{a_{s}}^{b_{s}}\,\left(
\left( \rho ^{-r}\right) _{b_{1}}^{j_{1}}...\left( \rho ^{-r}\right)
_{b_{s}}^{j_{s}}\right) .  \label{sedd}
\end{equation}
Defining $\xi ^{\left[ 0\right] }\doteq \mathbb{I}$, we can write $\varsigma
_{r,s}=\mathbb{I}_{r}\otimes \xi _{s}^{\left[ r\right] }$ for $r,s\in 
\mathbb{N}_{0}$. From the fact that $\phi $ and $\rho $ define algebra
automorphisms, it follows that $\varsigma $ is $\mathcal{B}\triangleright 
\mathcal{A}$-admissible. Let us see that $\varsigma $ is a 2-cocycle.
Straightforward calculations show that 
\begin{equation}
\xi _{s+t}^{\left[ r\right] }\thickapprox \xi _{s}^{\left[ r\right] }\otimes
\xi _{t}^{\left[ r\right] }\;\;\;and\;\;\;\xi _{t}^{\left[ r+s\right] }=\xi
_{t}^{\left[ r\right] }\,\xi _{t}^{\left[ s\right] }=\xi _{t}^{\left[ s%
\right] }\,\xi _{t}^{\left[ r\right] }.  \label{seda}
\end{equation}
Then, comparing 
\begin{equation*}
\varsigma _{r+s,t}\,\left( \varsigma _{r,s}\otimes \mathbb{I}_{t}\right)
=\left( \mathbb{I}_{r+s}\otimes \xi _{t}^{\left[ r+s\right] }\right)
\,\left( \mathbb{I}_{r}\otimes \xi _{s}^{\left[ r\right] }\otimes \mathbb{I}%
_{t}\right) =\mathbb{I}_{r}\otimes \xi _{s}^{\left[ r\right] }\otimes \xi
_{t}^{\left[ r+s\right] }
\end{equation*}
and 
\begin{equation*}
\varsigma _{r,s+t}\,\left( \mathbb{I}_{r}\otimes \varsigma _{s,t}\right)
=\left( \mathbb{I}_{r}\otimes \xi _{s+t}^{\left[ r\right] }\right) \,\left( 
\mathbb{I}_{r}\otimes \mathbb{I}_{s}\otimes \xi _{t}^{\left[ s\right]
}\right) =\mathbb{I}_{r}\otimes \xi _{s+t}^{\left[ r\right] }\,\left( 
\mathbb{I}_{s}\otimes \xi _{t}^{\left[ s\right] }\right) .
\end{equation*}
we have from Eq. $\left( \ref{seda}\right) $ that $\varsigma
_{r+s,t}\,\left( \varsigma _{r,s}\otimes \mathbb{I}_{t}\right) \thickapprox
\varsigma _{r,s+t}\,\left( \mathbb{I}_{r}\otimes \varsigma _{s,t}\right) $,
i.e. $\varsigma $ is an admissible element of $\frak{Z}^{2}\left[ \mathbf{B}%
_{1}^{\ast }\otimes \mathbf{A}_{1}\right] $. Moreover, using again Eq. $%
\left( \ref{seda}\right) $, we can see $\varsigma $ is a bicharacter and
anti-bicharacter. The twisting of $\mathcal{B}\triangleright \mathcal{A}$ by 
$\varsigma $ defines on the vector space $\mathbf{B}\triangleright \mathbf{A}
$ the product $m_{\varsigma }=m\,\Xi $, with $\Xi \,\Pi ^{\otimes 2}=\Pi
^{\otimes 2}\,\varsigma $. Then, the underlying algebra of $\left( \mathcal{B%
}\triangleright \mathcal{A}\right) _{\varsigma }$ is the one generated by $%
z_{i}^{j}=b^{j}\otimes a_{i}$ and quotient by the ideal algebraically
generated by 
\begin{equation*}
\left\{ \left\{ _{\sigma }R_{\lambda
_{n}}^{k_{1}...k_{n}}\;z_{k_{1}}^{j_{1}}\cdot _{\varsigma
}z_{k_{2}}^{j_{2}}\;...\;\cdot _{\varsigma }z_{k_{n}}^{j_{n}}\;\left(
_{\sigma }S^{\bot }\right) _{j_{1}...j_{n}}^{\omega _{n}}\right\} _{\omega
_{n}\in \Omega _{n}}^{\lambda _{n}\in \Lambda _{n}}\right\} _{n\in \mathbb{N}%
_{0}},
\end{equation*}
denoting by $\cdot _{\varsigma }$ the product $m_{\varsigma }$, and the
isomorphism $\underline{hom}^{\Upsilon }\left[ \mathcal{B},\mathcal{A}\right]
\backsimeq \underline{hom}\left[ \mathcal{B},\mathcal{A}\right] _{\varsigma
} $ follows. Therefore,

\begin{proposition}
The quantum spaces $\underline{hom}^{\Upsilon }\left[ \mathcal{B},\mathcal{A}%
\right] $ and $\underline{hom}\left[ \mathcal{B},\mathcal{A}\right] $ are
related by a twist transformation, i.e. 
\begin{equation*}
\underline{hom}^{\Upsilon }\left[ \mathcal{B},\mathcal{A}\right] \backsimeq 
\underline{hom}\left[ \mathcal{B},\mathcal{A}\right] _{\varsigma }
\end{equation*}
being $\varsigma $ a bicharacter in $\frak{Z}^{2}\left[ \mathbf{B}_{1}^{\ast
}\otimes \mathbf{A}_{1}\right] $.\ \ \ $\blacksquare $
\end{proposition}

For the coEnd objects the equation 
\begin{equation*}
\underline{end}\left[ \mathcal{A}\right] _{\chi }\backsimeq \underline{end}%
^{\Upsilon }\left[ \mathcal{A}\right] \backsimeq \underline{end}\left[ 
\mathcal{A}\right] _{\varsigma }
\end{equation*}
holds, provided $\varsigma =\digamma _{\chi }$ (see Eq. $\left( \ref{fgy}%
\right) $). Moreover, $\mathcal{A}_{\chi }$ and $\mathcal{B}_{\chi }$ can
alternatively be seen as twisting of the quantum spaces $\mathcal{A}=%
\underline{hom}\left[ \mathcal{K},\mathcal{A}\right] $ and $\mathcal{B}=%
\underline{hom}\left[ \mathcal{K},\mathcal{B}\right] $ by cocycles 
\begin{eqnarray*}
a_{k_{1}}...a_{k_{r}}\otimes a_{k_{r+1}}...a_{k_{r+s}} &\mapsto
&a_{k_{1}}...a_{k_{r}}\otimes \left( \phi ^{-r}\right)
_{k_{r+1}}^{j_{1}}...\left( \phi ^{-r}\right)
\,_{k_{r+s}}^{j_{s}}\,a_{j_{1}}...a_{j_{s}}, \\
&& \\
b_{k_{1}}...b_{k_{r}}\otimes b_{k_{r+1}}...b_{k_{r+s}} &\mapsto
&b_{k_{1}}...b_{k_{r}}\otimes \left( \rho ^{r}\right)
_{k_{r+1}}^{j_{1}}...\left( \rho ^{r}\right)
_{k_{r+s}}^{j_{s}}\,b_{j_{1}}...b_{j_{s}},
\end{eqnarray*}
in $\frak{Z}^{2}\left[ \mathbf{A}_{1}\right] $ and $\frak{Z}^{2}\left[ 
\mathbf{B}_{1}\right] $, respectively, which we also denote $\varsigma $. In
resume, 
\begin{equation*}
\underline{hom}^{\Upsilon }\left[ \mathcal{B},\mathcal{A}\right] \backsimeq 
\mathcal{B}_{\varsigma }\triangleright \mathcal{A}_{\varsigma }=\underline{%
hom}\left[ \mathcal{B}_{\varsigma },\mathcal{A}_{\varsigma }\right]
\backsimeq \underline{hom}\left[ \mathcal{B},\mathcal{A}\right] _{\varsigma
},
\end{equation*}
unifying in this way the correspondence between the objects $\underline{hom}%
^{\Upsilon }\left[ \mathcal{B},\mathcal{A}\right] $ and $\underline{hom}%
\left[ \mathcal{B},\mathcal{A}\right] $ (in particular, $\mathcal{A}%
^{\Upsilon }=\underline{hom}^{\Upsilon }\left[ \mathcal{K},\mathcal{A}\right]
$ and $\mathcal{A}=\underline{hom}\left[ \mathcal{K},\mathcal{A}\right] $),
just in terms of twisting of quantum spaces.

\subsection{Twistings and the cohomology relation}

We have seen that cocycle condition for a 2-cochain insure the associativity
of its related deformed product. But, when can it be insured these products
define a twisted quantum space isomorphic to the original one, or when two
twisted quantum spaces are isomorphic ? The following theorem gives the
answer to these questions.

\begin{theorem}
Given an $\mathcal{A}$-admissible cochain $\psi \in \frak{Z}^{2}$, $\mathcal{%
A}\backsimeq \mathcal{A}_{\psi }$ \emph{iff }$\psi =\partial \theta $, being 
$\theta $ an $\mathcal{A}$-\textbf{admissible} counital 1-cochain. More
generally, consider a pair of $\mathcal{A}$-admissible counital 2-cocycles $%
\varphi $ and $\psi $. Then, $\mathcal{A}_{\varphi }\backsimeq \mathcal{A}%
_{\psi }$ \emph{iff }$\varphi $ and $\psi $ are cohomologous through an $%
\mathcal{A}$-admissible counital 1-cochain $\theta $. In particular, if $%
\varphi $ and $\psi $ are related as in Eq. $\left( \ref{quo}\right) $, then 
$\varphi \backsim _{\theta }\psi $ through an admissible $\theta $.\ \ \ $%
\blacksquare $
\end{theorem}

We show below one of the implications. The proof of the other will be given
at the end of \S \textbf{3.2.2}.

\begin{proof}
\textbf{(part 1)} Suppose $\psi =\delta _{1}\theta \,\left( \theta \otimes
\theta \right) ^{-1}=\partial \theta $ with 
\begin{equation*}
\psi \left( \ker \Pi ^{\otimes 2}\right) =\ker \Pi ^{\otimes
2}\;\;and\;\;\theta \left( \ker \Pi \right) =\ker \Pi .
\end{equation*}
(Since $\partial ^{2}\theta =\mathbb{I}^{\otimes 3}$ for all $\theta \in 
\mathsf{C}^{1}$, $\psi $ is a 2-cocycle.) In particular $\psi $ defines a
linear automorphism $\Psi :\mathbf{A}^{\otimes 2}\rightarrow \mathbf{A}%
^{\otimes 2}$ such that $\Psi \,\Pi ^{\otimes 2}=\Pi ^{\otimes 2}\,\psi $,
which gives rise to the twisted product $m_{\psi }=m\,\Psi $. And $\theta $
defines, by the equation $\Theta \,\Pi =\Pi \,\theta $, a linear
automorphism $\Theta :\mathbf{A}\rightarrow \mathbf{A}$ satisfying $\Theta
\left( \mathbf{A}_{1}\right) \subset \mathbf{A}_{1}$. Then, 
\begin{equation*}
m_{\psi }\,\Pi ^{\otimes 2}=m\,\Pi ^{\otimes 2}\,\psi =m\,\Pi ^{\otimes
2}\,\delta _{1}\theta \,\left( \theta \otimes \theta \right) ^{-1}=\Pi
\,m_{\otimes }\,\delta _{1}\theta \,\left( \theta \otimes \theta \right)
^{-1}=\Pi \,\theta \,m_{\otimes }\,\left( \theta \otimes \theta \right)
^{-1},
\end{equation*}
where Eq. $\left( \ref{im}\right) $ was used in the last equality.
Multiplying to the right by $\theta \otimes \theta $ and using the relation
between $\theta $ and $\Theta $ we arrive at $m_{\psi }\,\left( \Theta
\otimes \Theta \right) \,\Pi ^{\otimes 2}=\Theta \,m\,\Pi ^{\otimes 2}$, and
since $\Pi $ is right invertible, equality $m_{\psi }\,\left( \Theta \otimes
\Theta \right) =\Theta \,m$ follows. This means the linear automorphism $%
\Theta $ defines an algebra map $\Theta :\mathbf{A}\rightarrow \mathbf{A}%
_{\psi }$ such that $\Theta \left( \mathbf{A}_{1}\right) \subset \mathbf{A}%
_{1}$, therefore it defines a quantum space isomorphism $\mathcal{A}%
\backsimeq \mathcal{A}_{\psi }$.

Now, consider a pair of admissible 2-cochains $\varphi $ and $\psi $ such
that $\varphi \backsim _{\theta }\psi $, being $\theta $ admissible. By
definition of cohomologous $\delta _{1}\theta \,\varphi \,\left( \theta
\otimes \theta \right) ^{-1}=\psi $. Hence, reproducing the calculations
above, $m_{\psi }\,\left( \Theta \otimes \Theta \right) =\Theta \,m_{\varphi
}$, being $\Theta $ the linear automorphism defined by $\theta $ and $\Pi $.
That is to say $\mathcal{A}_{\varphi }\backsimeq \mathcal{A}_{\psi }$.
\end{proof}

Consider on the set of $\mathcal{A}$-admissible counital 2-cocycles, the
cohomology relation through $\mathcal{A}$-admissible 1-cochains; i.e. the $%
\mathcal{A}$-admissible cochains $\varphi ,\psi \in \frak{Z}^{2}$ are
cohomologous if $\varphi \backsim _{\theta }\psi $, being $\theta $
admissible. It is an equivalence relation since proposition of \S \textbf{%
2.2.3 }and the fact that admissible $n$-cochains form a subgroup of $\mathsf{%
C}^{n}\left[ \mathbf{A}_{1}\right] $ for all $n$. Then, by theorem above,
the resulting quotient set, namely $H_{\mathcal{A}}^{2}\left[ \mathbf{A}_{1}%
\right] $, characterize completely the isomorphism classes of twisting of $%
\mathcal{A}$. Since the last claim of theorem above (related to Eq. $\left( 
\ref{quo}\right) $), the set $H_{\mathcal{A}}^{2}\left[ \mathbf{A}_{1}\right]
$ can be regarded as a quotient of $Aut_{\mathrm{Vct}}\left[ \mathbf{A}%
\otimes \mathbf{A}\right] $.

Following an analogous reasoning, one can define a space $H_{\mathcal{A}}^{1}%
\left[ \mathbf{A}_{1}\right] $, given by $\mathcal{A}$-admissible 1-cocycles
($ipso$ $facto$ counital) quotient by the cohomology relation through
admissible cochains (which gives the equality relation, as follows from the
general case analyzed in \S \textbf{2.2.1}). This space is in bijection with
the group of automorphisms of $\mathcal{A}$, i.e. $H_{\mathcal{A}}^{1}\left[ 
\mathbf{A}_{1}\right] \backsimeq Aut_{\mathrm{FGA}}\left[ \mathcal{A}\right] 
$.\footnote{%
To give an automorphism of $\mathcal{A}=\left( \mathbf{A}_{1},\mathbf{A}%
\right) $ is the same as to give an element $\alpha \in Aut\left[ \mathbf{A}%
_{1}\right] $ such that $\alpha ,\alpha ^{-1}$ can be extended to algebra
maps on $\mathbf{A}$. But these are precisely the admissible 1-cocycles in $%
\mathsf{C}^{1}\left[ \mathbf{A}_{1}\right] $. Since $\mathsf{Z}^{1}\left[ 
\mathbf{A}_{1}\right] =\frak{Z}^{1}\left[ \mathbf{A}_{1}\right] \backsimeq
Aut\left[ \mathbf{A}_{1}\right] $, we have the mentioned bijection.}

It is worth mentioning that \textbf{non} every admissible 2-cocycle is a
coboundary of an admissible 1-cochain. If this were the case, the
commutative algebra $\Bbbk \left[ a,b\right] $ of \S \textbf{3.1.1} would be
isomorphic to the quantum plane $\mathbf{A}_{\hslash }^{2|0}=\Bbbk \left[ a,b%
\right] _{\psi _{\hslash }}$. In particular, we can affirm that there do
exist non trivial twisting of quantum spaces, i.e. $H_{\mathcal{A}}^{2}\left[
\mathbf{A}_{1}\right] \neq \left\{ \mathbb{I}^{\otimes 2}\right\} $.
Nevertheless, for quantum spaces formed out by free algebras, every counital
2-cocycle gives rise to an isomorphic twisted quantum space.\footnote{%
If the twisted space is not also free, then it would have a non cero related
ideal, which would imply the underlying vector space is smaller than the
untwisted one.} This is a consequence of the fact that every counital
2-cocycle is a counital 2-coboundary (and the admissibility condition is
immediate), as we show below.

\subsubsection{The equality $\mathsf{Z}^{2}=\mathsf{B}^{2}$}

We have seen that $\mathsf{B}^{2}\subset \mathsf{Z}^{2}$ for every vector
space. In this subsection we show the other inclusion and some related
consequences.

\begin{theorem}
Every $\psi \in \mathsf{Z}^{2}$ is cohomologous to the identity, i.e. $\psi
\in \mathsf{B}^{2}$. Moreover, for every linear automorphism $\varpi :%
\mathbf{A}_{1}\backsimeq \mathbf{A}_{1}$, there exists a unique $\theta \in 
\mathsf{C}^{1}$ such that $\psi =\partial \theta $ and $\theta _{1}=\varpi $.
\end{theorem}

\begin{proof}
Let us come back to the proof of \textbf{Theor. 8}. Given a cochain $\psi
\in \mathsf{C}^{2}$ (non necessarily a counital 2-cocycle) we have defined
inductively, from Eq. $\left( \ref{ind}\right) $, a family of linear
isomorphisms $\theta _{n}:\mathbf{A}_{1}^{\otimes n}\backsimeq \mathbf{A}%
_{1}^{\otimes n}$, $n\geq 2$. In an analogous way, take a map $\varpi :%
\mathbf{A}_{1}\backsimeq \mathbf{A}_{1}$ and consider $\theta \in \mathsf{C}%
^{1}$ such that $\theta _{0}\doteq \psi _{0,0}^{-1}=\left. 1\right/ \psi
_{0,0}$, $\theta _{1}\doteq \varpi $, and 
\begin{equation*}
\theta _{n+1}\doteq \psi _{n,1}\,\left( \theta _{n}\otimes \varpi \right)
,\;\;n\in \mathbb{N}.
\end{equation*}
Let us show that if $\psi $ is a 2-cocycle then, using the identification
given by Eq. $\left( \ref{iden}\right) $, 
\begin{equation}
\theta _{r+s}\thickapprox \psi _{r,s}\,\left( \theta _{r}\otimes \theta
_{s}\right) ,  \label{equi}
\end{equation}
or equivalently $\psi =\delta _{1}\theta \,\left( \theta \otimes \theta
\right) ^{-1}=\partial \theta $. We will make induction on $s$. Recall that,
since Eq. $\left( \ref{2cc2}\right) $, $\psi _{r,0}\thickapprox \psi
_{0,0}\cdot \mathbb{I}_{r}$ for all 2-cocycle. Then for $s=0$, by definition
of $\theta _{0}$ and $\theta _{1}$, Eq. $\left( \ref{equi}\right) $ holds.
For $s+1$, 
\begin{equation*}
\psi _{r,s+1}\,\left( \theta _{r}\otimes \theta _{s+1}\right) \thickapprox
\psi _{r,s+1}\,\left( \theta _{r}\otimes \psi _{s,1}\,\left( \theta
_{s}\otimes \varpi \right) \right) \,=\psi _{r,s+1}\,\left( \mathbb{I}%
_{r}\otimes \psi _{s,1}\right) \,\left( \theta _{r}\otimes \theta
_{s}\otimes \varpi \right) .
\end{equation*}
Using Eq. $\left( \ref{2cc}\right) $ for $t=1$, 
\begin{equation*}
\psi _{r,s+1}\,\left( \theta _{r}\otimes \theta _{s+1}\right) \thickapprox
\psi _{r+s,1}\,\left( \psi _{r,s}\otimes \mathbb{I}\right) \,\left( \theta
_{r}\otimes \theta _{s}\otimes \varpi \right) =\psi _{r+s,1}\,\left( \psi
_{r,s}\,\left( \theta _{r}\otimes \theta _{s}\right) \otimes \varpi \right) ,
\end{equation*}
and from the inductive hypothesis (i.e. validity of Eq. $\left( \ref{equi}%
\right) $) 
\begin{equation*}
\psi _{r,s+1}\,\left( \theta _{r}\otimes \theta _{s+1}\right) \thickapprox
\psi _{r+s,1}\,\left( \theta _{r+s}\otimes \varpi \right) =\theta _{r+s+1},
\end{equation*}
as we wanted to show.

Now, consider another $\chi $ such that $\psi =\partial \chi $ and $\chi
_{1}=\varpi $. From $\left( \ref{equi}\right) $, $\chi _{0}=\left. 1\right/
\psi _{0,0}$ and the rest of $\chi $ is given inductively by the equation $%
\chi _{n+1}=\psi _{n,1}\,\left( \chi _{n}\otimes \varpi \right) $,
concluding in this way the proof of the our claim.
\end{proof}

\begin{corollary}
If $\psi \in \frak{Z}^{2}$, then $\psi \in \frak{B}^{2}$. In addition,
fixing an automorphism $\varpi $ of $\mathbf{A}_{1}$, there exists a unique $%
\theta \in \frak{C}^{1}$ such that $\psi =\partial \theta $ and $\theta
_{1}=\varpi $.
\end{corollary}

\begin{proof}
From previous theorem, if $\psi \in \frak{Z}^{2}=\mathsf{Z}^{2}\cap \frak{C}%
^{2}$, then $\psi \in \mathsf{B}^{2}\cap \frak{C}^{2}$ and there exists a
unique $\theta \in \mathsf{C}^{1}$ such that $\psi =\partial \theta $ and $%
\theta _{1}=\varpi $. But, as we have shown in \S \textbf{2.2.1} (see Eq. $%
\left( \ref{eqy}\right) $), $\mathsf{B}^{2}\cap \frak{C}^{2}=\frak{B}^{2}$.
Thus $\psi \in \frak{B}^{2}$ and accordingly $\theta \in \frak{C}^{1}$.
\end{proof}

This corollary can also be proven by noting that $\psi $ is counital if and
only if $\psi _{0,0}=1$ (see Eq. $\left( \ref{c2c}\right) $). This makes $%
\theta _{0}=\mathbb{I}_{0}=1$, i.e. $\theta $ is counital.

\begin{corollary}
Every twisting of the quantum space $\mathcal{A}=\left( \mathbf{A}_{1},%
\mathbf{A}_{1}^{\otimes }\right) $ is isomorphic to $\mathcal{A}$.
\end{corollary}

\begin{proof}
Every $n$-cochain in $\mathsf{C}^{n}\left[ \mathbf{A}_{1}\right] $ is
obviously admissible for the quantum space $\mathcal{A}=\left( \mathbf{A}%
_{1},\mathbf{A}_{1}^{\otimes }\right) $ (since $\ker \Pi =0$). Then, from
above results, given a counital 2-cocycle $\psi $ there exists a counital
1-cochain $\theta $ (admissible for $\mathcal{A}$) such that $\psi =\partial
\theta $. In this situation, \textbf{Theor. 9} insures $\mathcal{A}%
\backsimeq \mathcal{A}_{\psi }$.
\end{proof}

These results motive the definition of a subgroup $\frak{P}^{1}\left[ 
\mathbf{A}_{1}\right] \subset \frak{C}^{1}\left[ \mathbf{A}_{1}\right] $
constituted by counital 1-cochains $\theta $ such that $\theta _{1}=\mathbb{I%
}$. We shall call \textbf{primitive }the elements of this subgroup. From now
on, when we write $\psi =\partial \theta $, we are supposing $\theta $ is in 
$\frak{P}^{1}$, unless we say the contrary.

As an example, the cochain 
\begin{equation}
\varsigma _{\mathcal{A}}:a_{k_{1}}...a_{k_{r}}\otimes
a_{k_{r+1}}...a_{k_{r+s}}\mapsto a_{k_{1}}...a_{k_{r}}\otimes \left( \sigma
_{\mathcal{A}}^{-r}\right) _{k_{r+1}}^{j_{1}}...\left( \sigma _{\mathcal{A}%
}^{-r}\right) \,_{k_{r+s}}^{j_{s}}\,a_{j_{1}}...a_{j_{s}}  \label{sa}
\end{equation}
related to the twisted internal coHom objects $\mathcal{B}_{\varsigma
}\triangleright \mathcal{A}_{\varsigma }$, has primitive 
\begin{equation*}
\theta _{\mathcal{A}}:a_{k_{1}}...a_{k_{r}}\mapsto \delta
_{k_{1}}^{j_{1}}\left( \sigma _{\mathcal{A}}^{-1}\right)
_{k_{2}}^{j_{2}}\left( \sigma _{\mathcal{A}}^{-2}\right)
_{k_{3}}^{j_{3}}...\left( \sigma {}_{\mathcal{A}}^{-r+1}\right)
\,_{k_{r}}^{j_{r}}\,a_{j_{1}}...a_{j_{r}}.
\end{equation*}

\bigskip

Since we have shown in \S \textbf{2.2.3} that cohomology relation is an
equivalence relation in $\mathsf{C}^{2}\left[ \mathbf{A}_{1}\right] $,
previous results imply, in particular, that every pair of 2-cocycles are
cohomologous. The following proposition tell us which class of 1-cochains
implement such a relation.

\begin{proposition}
Consider a pair of cochains $\psi ,\varphi \in \frak{Z}^{2}$ with primitive $%
\lambda ,\chi \in \frak{P}^{1}$, respectively.

a\emph{) }They are cohomologous through $\theta $ \emph{iff} $\theta
=\lambda \,\omega \,\chi ^{-1}$, with $\omega \in \frak{Z}^{1}$.

b\emph{) }They are cohomologous through an admissible cochain \emph{iff}
there exists $\omega \in \frak{Z}^{1}$ such that $\lambda \,\omega \,\chi
^{-1}$ is admissible.
\end{proposition}

\begin{proof}
a) If $\psi =\partial \lambda $, $\varphi =\partial \chi $, $\theta =\lambda
\,\omega \,\chi ^{-1}$ and $\partial \omega =\mathbb{I}^{\otimes 2}$ (i.e. $%
\omega $ is a 1-cocycle), then 
\begin{equation*}
\begin{array}{l}
\delta _{1}\theta \,\varphi \,\left( \theta \otimes \theta \right)
^{-1}=\delta _{1}\left( \lambda \,\omega \,\chi ^{-1}\right) \,\left( \delta
_{1}\,\chi \,\left( \chi \otimes \chi \right) ^{-1}\right) \,\left( \lambda
\,\omega \,\chi ^{-1}\otimes \lambda \,\omega \,\chi ^{-1}\right) ^{-1} \\ 
\\ 
=\delta _{1}\left( \lambda \,\omega \,\right) \,\left( \lambda \,\omega
\,\otimes \lambda \,\omega \,\right) ^{-1}=\delta _{1}\lambda \,\left(
\delta _{1}\omega \,\left( \,\omega \,\otimes \,\omega \,\right)
^{-1}\right) \,\left( \lambda \,\,\otimes \lambda \,\,\right) ^{-1}=\delta
_{1}\lambda \,\,\left( \lambda \,\,\otimes \lambda \,\,\right) ^{-1}=\psi .
\end{array}
\end{equation*}
That means $\varphi \backsim _{\theta }\psi $, and one implication follows.
Now suppose $\varphi \backsim _{\theta }\psi $. We always can write $\theta
=\lambda \,\omega \,\chi ^{-1}$, for every $\theta \in \mathsf{C}^{1}$. It
is enough to take $\omega =\lambda ^{-1}\,\theta \,\chi $. Let us see that $%
\omega \in \frak{Z}^{1}$. From equation above if $\varphi \backsim _{\theta
}\psi $ we have that 
\begin{equation*}
\delta _{1}\lambda \,\left( \delta _{1}\omega \,\left( \,\omega \,\otimes
\,\omega \,\right) ^{-1}\right) \,\left( \lambda \,\,\otimes \lambda
\,\,\right) ^{-1}=\delta _{1}\lambda \,\,\left( \lambda \,\,\otimes \lambda
\,\,\right) ^{-1},
\end{equation*}
which is fulfilled if and only if $\delta _{1}\omega \,=\omega \,\otimes
\,\omega $, i.e. $\omega \in \frak{Z}^{1}$, as we wanted to show.

b) This part is an immediate consequence of the first one.
\end{proof}

Recall that to give an element $\omega \in \frak{Z}^{1}\left[ \mathbf{A}_{1}%
\right] $ is the same as giving one in $\varpi \in Aut\left[ \mathbf{A}_{1}%
\right] $, since every 1-cocycle is of the form $\omega =\varpi ^{\otimes }$%
. This connects above proposition with \textbf{Theor. 9}.

\bigskip

A brief comment about bialgebras is in order. The quasicomplex $\mathsf{G}%
^{\bullet }\left[ \mathbf{A}_{1}^{\otimes }\right] $ also satisfies that
every 2-cocycle is a 2-coboundary. In fact, given a cocycle $\chi \in 
\mathsf{G}^{2}$, a cochain $\lambda \in \mathsf{G}^{1}$, such that $d\lambda
=\chi $, can be inductively defined by the formula 
\begin{equation*}
\lambda _{0}=\left. 1\right/ \chi _{0,0};\;\;\lambda _{n+1}=\left( \lambda
_{n}\otimes \varepsilon \right) \ast \chi _{n,1}.
\end{equation*}
Moreover, every pair of 2-cocycles in $\mathsf{G}_{\mathcal{Z}}^{2}\left[ 
\mathbf{A}_{1}^{\otimes }\right] $ are cohomologous.\footnote{%
The cochains that implement the cohomology relation have a form analogous to
the ones given in last proposition for the complex $\mathsf{C}^{\bullet }%
\left[ \mathbf{A}_{1}\right] $.} But for $\mathsf{G}_{\mathcal{Z}}^{\bullet }%
\left[ \mathbf{A}\right] $, given a 2-cocycle $\phi $ there, we can just say 
$\phi \,\Pi ^{\otimes 2}=d\vartheta $; and $\phi $ will be a coboundary 
\emph{iff }$\vartheta =\theta \,\Pi $.

\subsubsection{Role of primitive 1-cochains}

The role of the 1-cochains $\theta $ defining a given 2-cocycle $\psi $, for
a generic quantum space $\mathcal{A}$, is expressed by the following theorem.

\begin{theorem}
Let $\mathcal{A}=\left( \mathbf{A}_{1},\mathbf{A}\right) $ be a quantum
space and consider an admissible $\psi \in \frak{Z}^{2}\left[ \mathbf{A}_{1}%
\right] $. The canonical epimorphism $\mathbf{A}_{1}^{\otimes
}\twoheadrightarrow \mathbf{A}_{\psi }$ associated to $\mathcal{A}_{\psi
}=\left( \mathbf{A}_{1},\mathbf{A}_{\psi }\right) $ is given by $\Pi
\,\theta $, being $\theta \in \frak{P}^{1}\left[ \mathbf{A}_{1}\right] $ the
primitive of $\psi $. Furthermore, if $\mathbf{I}$ is the ideal related to $%
\mathcal{A}$, then $\mathbf{I}_{\psi }$, the one related to $\mathcal{A}%
_{\psi }$, is equal to $\theta ^{-1}\left( \mathbf{I}\right) $.
\end{theorem}

\begin{proof}
Since $\theta _{0,1}=\mathbb{I}_{0,1}$ (which means $\theta \left( \mathbf{A}%
_{1}\right) \equiv \mathbf{A}_{1}$ and $\theta \left( \lambda \right) =\eta
\left( \lambda \right) $), $\Pi \,\theta $ defines the inclusion $\mathbf{A}%
_{1}\hookrightarrow \mathbf{A}_{\psi }$ and a unit preserving map. If we
show that 
\begin{equation}
m_{\otimes }\,\psi \,\left( \theta \otimes \theta \right) =\theta
\,m_{\otimes },  \label{ppal}
\end{equation}
what implies (multiplying by $\Pi $ to the left) $m_{\psi }\,\left( \Pi
\,\theta \otimes \Pi \,\theta \right) =\Pi \,\theta \,m_{\otimes }$, we are
proving $\Pi \,\theta $ is an algebra homomorphism $\mathbf{A}_{1}^{\otimes
}\rightarrow \mathbf{A}_{\psi }$.

Eq. $\left( \ref{ppal}\right) $ restricted to a subspace $\mathbf{A}%
_{1}^{\otimes r}\otimes \mathbf{A}_{1}^{\otimes s}$ means $\psi
_{r,s}\,\left( \theta _{r}\otimes \theta _{s}\right) =\left[ \delta
_{1}\theta \right] _{r,s}$, or equivalently, $\psi _{r,s}\,\left( \theta
_{r}\otimes \theta _{s}\right) \thickapprox \theta _{r+s}$. But this is true
since $\partial \theta =\psi $, hence Eq. $\left( \ref{ppal}\right) $ holds.

Finally, we must show $\theta ^{-1}\left( \mathbf{I}\right) $ is equal to $%
\mathbf{I}_{\psi }$, the ideal related to $\mathcal{A}_{\psi }$. By
definition, $\mathbf{I}=\ker \Pi $ and, from the last result, $\mathbf{I}%
_{\psi }=\ker \Pi \theta $. In addition, it is well known that $\ker \Pi
\theta =\theta ^{-1}\left( \ker \Pi \right) $, and consequently our claim
follows.
\end{proof}

Thus we can write (compare to Eq. $\left( \ref{pitita}\right) $) 
\begin{equation*}
a_{j_{1}}\cdot _{\psi }...\cdot _{\psi }a_{j_{n}}=\Pi \theta \left(
a_{i_{1}}\otimes ...\otimes a_{i_{n}}\right) ,
\end{equation*}
being $"\cdot _{\psi }"=m_{\psi }$. On the other hand, note that if $\theta $
is $\mathcal{A}$-admissible, then $\mathbf{I}_{\psi }=\theta ^{-1}\left( 
\mathbf{I}\right) =\mathbf{I}$ and, therefore, $\mathcal{A}\backsimeq 
\mathcal{A}_{\psi }$, as we proved before.

As an immediate corollary of last theorem we have:

\begin{corollary}
Consider a conic quantum space $\mathcal{A}$ with $\mathbf{I}%
=\bigoplus_{n\geq 2}\mathbf{I}_{n}$, and a counital $\mathcal{A}$-admissible
2-cocycle $\psi =\partial \theta \in \frak{Z}^{2}\left[ \mathbf{A}_{1}\right]
$. The ideal related to the twisted quantum space $\mathcal{A}_{\psi }$ is 
\begin{equation*}
\mathbf{I}_{\psi }=\bigoplus\nolimits_{n\geq 2}\mathbf{I}_{\psi
,n}=\bigoplus\nolimits_{n\geq 2}\theta ^{-1}\left( \mathbf{I}_{n}\right)
.\;\;\;\blacksquare
\end{equation*}
\end{corollary}

In coordinates, if $\mathbf{I}$ is linearly generated by 
\begin{equation*}
\left\{ R_{\lambda _{n}}^{k_{1}...k_{n}}\;a_{k_{1}}...a_{k_{n}}\right\}
_{\lambda _{n}\in \Lambda _{n}}\subset \mathbf{I}_{n},
\end{equation*}
then $\mathbf{I}_{\psi }$ is linearly generated by the elements 
\begin{equation*}
R_{\lambda _{n}}^{k_{1}...k_{n}}\;\theta ^{-1}\left(
a_{k_{1}}...a_{k_{n}}\right) =R_{\lambda _{n}}^{k_{1}...k_{n}}\;\left(
\theta _{n}^{-1}\right)
_{k_{1}...k_{n}}^{j_{1}...j_{n}}\,a_{j_{1}}...a_{j_{n}}.
\end{equation*}

Last corollary gives us another way to see that twisting of conic quantum
spaces are also conic. Further work shows the same for $m$-th quantum
spaces. In what follows, and to avoid any confusion, we shall write $\mathbf{%
X}\cdot \mathbf{Y}\subset $ $\mathbf{A}_{1}^{\otimes }$ for the image under $%
m_{\otimes }$ of a couple of subspaces $\mathbf{X},\mathbf{Y}\subset $ $%
\mathbf{A}_{1}^{\otimes }$.

\begin{proposition}
If $\mathcal{A}\in \mathrm{CA}^{m}$ and $\mathbf{I}=I\left[ \mathbf{Y}\right]
$, the ideal generated by $\mathbf{Y}\subset \mathbf{A}_{1}^{\otimes m}$,
then $\mathcal{A}_{\psi }\in \mathrm{CA}^{m}$ and $\mathbf{I}_{\psi }=I\left[
\theta ^{-1}\left( \mathbf{Y}\right) \right] $. That is to say, twist
transformations of $\mathcal{A}\in \mathrm{CA}^{m}$ give objects of $\mathrm{%
CA}^{m}$. Moreover, if $\mathcal{A}\in \mathrm{CA}$ has an ideal $\mathbf{I}$
generated by a graded vector subspace of $\mathbf{A}_{1}^{\otimes }$ 
\begin{equation*}
\mathbf{S}=\bigoplus\nolimits_{n\in \mathbb{N}_{0}}\mathbf{S}_{n};\;\;%
\mathbf{S}_{0,1}\doteq \left\{ 0\right\} ,
\end{equation*}
i.e. $\mathbf{I=}I\left[ \mathbf{S}\right] =\mathbf{A}_{1}^{\otimes }\cdot 
\mathbf{S}\cdot \mathbf{A}_{1}^{\otimes }$,\footnote{%
This is the case of quantum spaces of the form $\mathcal{A}\bigcirc \mathcal{%
B}$ with $\bigcirc =\bullet ,\triangleleft ,\triangleright ,\diamond $ or $%
\odot $, and $\mathcal{A}^{!}$.} then the ideal related to $\mathcal{A}%
_{\psi }$ is 
\begin{equation*}
\mathbf{I}_{\psi }=\theta ^{-1}\left( I\left[ \mathbf{S}\right] \right) =I%
\left[ \theta ^{-1}\left( \mathbf{S}\right) \right] =\mathbf{A}_{1}^{\otimes
}\cdot \theta ^{-1}\left( \mathbf{S}\right) \cdot \mathbf{A}_{1}^{\otimes
}.\;\;\;\blacksquare
\end{equation*}
\end{proposition}

Before going to the proof, let us show the lemma below.

\begin{lemma}
Let $\mathcal{A}$ be a \textbf{conic} quantum space with ideal $\mathbf{I}%
=\bigoplus_{n\geq 2}\mathbf{I}_{n}$. $\psi =\partial \theta $ is admissible 
\emph{iff} 
\begin{equation}
\theta ^{-1}\left( \mathbf{A}_{1}^{\otimes r}\cdot \mathbf{I}_{s}+\mathbf{I}%
_{r}\cdot \mathbf{A}_{1}^{\otimes s}\right) =\mathbf{A}_{1}^{\otimes r}\cdot
\theta ^{-1}\left( \mathbf{I}_{s}\right) +\theta ^{-1}\left( \mathbf{I}%
_{r}\right) \cdot \mathbf{A}_{1}^{\otimes s}.  \label{app}
\end{equation}
\end{lemma}

\begin{proof}
If $\psi $ is admissible, then (see Eq. $\left( \ref{ap}\right) $) 
\begin{equation*}
\psi _{r,s}\left( \mathbf{A}_{1}^{\otimes r}\otimes \mathbf{I}_{s}+\mathbf{I}%
_{r}\otimes \mathbf{A}_{1}^{\otimes s}\right) =\mathbf{A}_{1}^{\otimes
r}\otimes \mathbf{I}_{s}+\mathbf{I}_{r}\otimes \mathbf{A}_{1}^{\otimes s}.
\end{equation*}
Since $\psi =\delta _{1}\theta \,\left( \theta \otimes \theta \right) ^{-1}$%
, equation above says 
\begin{equation}
\left[ \delta _{1}\theta \,^{-1}\right] _{r,s}\left( \mathbf{A}_{1}^{\otimes
r}\otimes \mathbf{I}_{s}+\mathbf{I}_{r}\otimes \mathbf{A}_{1}^{\otimes
s}\right) =\left( \theta _{r}\otimes \theta _{s}\right) ^{-1}\left( \mathbf{A%
}_{1}^{\otimes r}\otimes \mathbf{I}_{s}+\mathbf{I}_{r}\otimes \mathbf{A}%
_{1}^{\otimes s}\right) .  \label{ap1}
\end{equation}
Applying to the left $m_{\otimes }$, using $m_{\otimes }\,\delta _{1}\theta
=\theta \,m_{\otimes }$ and $\theta \left( \mathbf{A}_{1}^{\otimes r}\right)
=\mathbf{A}_{1}^{\otimes r}$, we arrive precisely at Eq. $\left( \ref{app}%
\right) $. Reciprocally, if $\theta $ satisfies $\left( \ref{app}\right) $,
Eq. $\left( \ref{ap1}\right) $ follows immediately and with it the
admissibility of $\partial \theta $.
\end{proof}

\begin{proof}
\textbf{(of proposition)} Let us consider the more general case. Suppose we
have $\mathcal{A}\in \mathrm{CA}$ with 
\begin{equation*}
I\left[ \mathbf{S}\right] =\bigoplus\nolimits_{n\geq 2}\mathbf{I}_{n};\;\;%
\mathbf{I}_{n}=\sum\nolimits_{r=2}^{n}\sum\nolimits_{i=0}^{n-r}\mathbf{A}%
_{1}^{\otimes n-r-i}\cdot \mathbf{S}_{r}\cdot \mathbf{A}_{1}^{\otimes i},
\end{equation*}
and an $\mathcal{A}$-admissible 2-cocycle $\psi =\partial \theta $. From
lemma above, since $\mathbf{I}_{1}=\left\{ 0\right\} $, we have 
\begin{equation}
\theta ^{-1}\left( \mathbf{A}_{1}\cdot \mathbf{I}_{r}\right) =\mathbf{A}%
_{1}\cdot \theta ^{-1}\left( \mathbf{I}_{r}\right) \;\;and\;\;\theta
^{-1}\left( \mathbf{I}_{r}\cdot \mathbf{A}_{1}\right) =\theta ^{-1}\left( 
\mathbf{I}_{r}\right) \cdot \mathbf{A}_{1},  \label{4}
\end{equation}
for all $r$, and in particular 
\begin{equation*}
\theta ^{-1}\left( \mathbf{A}_{1}\cdot \mathbf{S}_{2}\right) =\mathbf{A}%
_{1}\cdot \theta ^{-1}\left( \mathbf{S}_{2}\right) \;\;and\;\;\theta
^{-1}\left( \mathbf{S}_{2}\cdot \mathbf{A}_{1}\right) =\theta ^{-1}\left( 
\mathbf{S}_{2}\right) \cdot \mathbf{A}_{1}.
\end{equation*}
Consequently, $\theta ^{-1}\left( \mathbf{I}_{3}\right) =\mathbf{A}_{1}\cdot
\theta ^{-1}\left( \mathbf{S}_{2}\right) +\theta ^{-1}\left( \mathbf{S}%
_{2}\right) \cdot \mathbf{A}_{1}+\theta ^{-1}\left( \mathbf{S}_{3}\right) $.
We can show by induction on $n\geq 3$ that 
\begin{equation*}
\theta ^{-1}\left( \mathbf{I}_{n}\right)
=\sum\nolimits_{r=2}^{n}\sum\nolimits_{i=0}^{n-r}\mathbf{A}_{1}^{\otimes
n-r-i}\cdot \theta ^{-1}\left( \mathbf{S}_{r}\right) \cdot \mathbf{A}%
_{1}^{\otimes i}.
\end{equation*}
In fact, because $\mathbf{I}_{n+1}=\mathbf{A}_{1}\cdot \mathbf{I}_{n}+%
\mathbf{I}_{n}\cdot \mathbf{A}_{1}+\mathbf{S}_{n+1}$, and using $\left( \ref
{4}\right) $%
\begin{equation*}
\theta ^{-1}\left( \mathbf{I}_{n+1}\right) =\theta ^{-1}\left( \mathbf{A}%
_{1}\cdot \mathbf{I}_{n}+\mathbf{I}_{n}\cdot \mathbf{A}_{1}+\mathbf{S}%
_{n+1}\right) =\mathbf{A}_{1}\cdot \theta ^{-1}\left( \mathbf{I}_{n}\right)
+\theta ^{-1}\left( \mathbf{I}_{n}\right) \cdot \mathbf{A}_{1}+\theta
^{-1}\left( \mathbf{S}_{n+1}\right) ,
\end{equation*}
we have from inductive hypothesis 
\begin{equation*}
\begin{array}{l}
\theta ^{-1}\left( \mathbf{I}_{n+1}\right)
=\sum\nolimits_{r=2}^{n}\sum\nolimits_{i=0}^{n-r}\mathbf{A}_{1}^{\otimes
n+1-r-i}\cdot \theta ^{-1}\left( \mathbf{S}_{r}\right) \cdot \mathbf{A}%
_{1}^{\otimes i} \\ 
\\ 
+\sum\nolimits_{r=2}^{n}\sum\nolimits_{i=0}^{n-r}\mathbf{A}_{1}^{\otimes
n-r-i}\cdot \theta ^{-1}\left( \mathbf{S}_{r}\right) \cdot \mathbf{A}%
_{1}^{\otimes i+1}+\theta ^{-1}\left( \mathbf{S}_{n+1}\right)
=\sum\nolimits_{r=2}^{n+1}\sum\nolimits_{i=0}^{n+1-r}\mathbf{A}_{1}^{\otimes
n+1-r-i}\cdot \theta ^{-1}\left( \mathbf{S}_{r}\right) \cdot \mathbf{A}%
_{1}^{\otimes i}.
\end{array}
\end{equation*}
Thus, $\theta ^{-1}\left( I\left[ \mathbf{S}\right] \right) =I\left[ \theta
^{-1}\left( \mathbf{S}\right) \right] $, as we wanted to see.
\end{proof}

\textbf{Example: }In the quadratic case $\mathbf{I}=I\left[ \mathbf{Y}\right]
$ is an ideal generated by $\mathbf{Y}\subset \mathbf{A}_{1}^{\otimes 2}$.
Then, since $\theta _{2}\doteq \psi _{1,1}$, we have that $\mathbf{I}_{\psi
}=I\left[ \psi _{1,1}^{-1}\left( Y\right) \right] $. Let us suppose $\mathbf{%
Y}$ is given by a set $Y_{\lambda }^{kl}\,a_{k}\otimes a_{l}$, $\lambda \in
\Lambda $. If $\psi _{1,1}\doteq M\in GL\left( n^{2}\right) $, $\mathbf{I}%
_{\psi }$ is algebraically generated by the elements $Y_{\lambda
}^{kl}\,N_{kl}^{ij}\,a_{i}\otimes a_{j}$, with $N=M^{-1}$. For a freely
commutative algebra, we have $\Lambda =\left[ n\right] \times \left[ n\right]
$ and $Y_{\lambda }^{kl}=Y_{ij}^{kl}=\delta _{i}^{k}\,\delta
_{j}^{l}-P_{ij}^{kl}$. From all that Equation $\left( \ref{dr}\right) $
follows.\ \ \ $\blacksquare $

\bigskip

To end this subsection, let us prove the other implication of \textbf{Theor.
9}.

\begin{proof}
\textbf{(of Theor. 9, part 2)} Suppose $\mathcal{A}_{\varphi }\backsimeq 
\mathcal{A}_{\psi }$. That is to say, there exists a linear automorphism $%
\Theta :\mathbf{A}\backsimeq \mathbf{A}$ such that $\Theta \left( 1\right)
=1 $ and, in a basis $\left\{ a_{i}\right\} $ of $\mathbf{A}_{1}$, $\Theta
\left( a_{i}\right) =\varpi _{i}^{j}\,a_{j}$ and 
\begin{equation*}
\Theta \left( a_{i_{1}}\cdot _{\varphi }...\cdot _{\varphi }a_{i_{k}}\right)
=\varpi _{i_{1}}^{j_{1}}...\varpi _{i_{k}}^{j_{k}}\,a_{j_{1}}\cdot _{\psi
}...\cdot _{\psi }a_{j_{k}},
\end{equation*}
being $\varpi _{i}^{j}$ the matrix elements of an automorphism $\varpi :%
\mathbf{A}_{1}\backsimeq $ $\mathbf{A}_{1}$. If $\chi $ and $\lambda $ are
the primitive of $\varphi $ and $\psi $, resp., the last equation can be
written $\Theta \,\Pi \,\chi _{k}=\Pi \,\lambda _{k}\,\varpi ^{\otimes k}$, $%
k\in \mathbb{N}_{0}$. Thus, the 1-cochain $\theta $ given by $\theta
_{k}=\lambda _{k}\,\varpi ^{\otimes k}\,\chi _{k}^{-1}$ is $\mathcal{A}$%
-admissible. Calling $\omega $ the 1-cocycle such that $\omega _{r}=\varpi
^{\otimes r}$, then $\theta =\lambda \,\omega \,\chi ^{-1}$. From \textbf{%
Prop.} \textbf{10} it follows that $\varphi \backsim _{\theta }\psi $. In
particular, if $\Pi ^{\otimes 2}\left( \varphi -\psi \right) =0$ (see Eq. $%
\left( \ref{quo}\right) $), then $m_{\varphi }=m_{\psi }$ and, in
consequence, $\mathcal{A}_{\varphi }=\mathcal{A}_{\psi }$. That implies $%
\varphi $ and $\psi $ are cohomologous through the admissible 1-cochain $%
\lambda \,\chi ^{-1}$ (since $\varpi =\mathbb{I}$).
\end{proof}

\subsection{The gauge equivalence}

\subsubsection{Composition and inversion of twist transformations}

In this subsection we show that consecutive applications of twist
transformations is again a twist transformation, and that twist
transformations have inverse. To start with, the following results will be
crucial.

\begin{proposition}
The map $\partial $ defines a surjection $\frak{C}^{1}\twoheadrightarrow 
\frak{Z}^{2}$ which restricted to the subgroup $\frak{P}^{1}\subset \frak{C}%
^{1}$ became a bijection $\frak{P}^{1}\backsimeq \frak{Z}^{2}$.\ \ \ \ $%
\blacksquare $
\end{proposition}

This is an immediate corollary of \textbf{Theor}. \textbf{10}.

\begin{proposition}
The set $\frak{Z}^{2}$ is a group under the product 
\begin{equation*}
\left( \psi ,\varphi \right) \mapsto \psi \star \varphi \doteq \psi \,\left(
\theta \otimes \theta \right) \,\varphi \,\left( \theta \otimes \theta
\right) ^{-1},
\end{equation*}
being $\theta \in \frak{P}^{1}$ the primitive of $\psi $. With respect to
this group structure the bijection $\frak{P}^{1}\backsimeq \frak{Z}^{2}$ is
a group isomorphism.
\end{proposition}

\begin{proof}
From $\partial :\frak{P}^{1}\backsimeq \frak{Z}^{2}$ we have that $\psi
\star \varphi =\partial \theta \star \partial \chi $ is equal to 
\begin{equation*}
\partial \theta \,\left( \theta \otimes \theta \right) \,\partial \chi
\,\left( \theta \otimes \theta \right) ^{-1}=\delta _{1}\theta \,\delta
_{1}\chi \,\left( \theta \otimes \theta \right) ^{-1}\,\left( \chi \otimes
\chi \right) ^{-1}=\delta _{1}\theta \,\left( \theta \,\chi \otimes \theta
\,\chi \right) ^{-1}=\partial \left( \theta \,\chi \right) ,
\end{equation*}
thus $\theta \,\chi \in \frak{P}^{1}$ is sent to $\psi \star \varphi $ \emph{%
via }$\partial $. In other words, $\partial $ restricted to $\frak{P}^{1}$
translates the product on its domain into the map $\star $. Hence, $\star $
is an associative product with unit $\mathbb{I}^{\otimes 2}=\partial \mathbb{%
I}$, and each $\psi \in \frak{Z}^{2}$ has inverse 
\begin{equation*}
\frak{i}\psi =\partial \left( \theta ^{-1}\right) =\left( \delta _{1}\theta
\right) ^{-1}\,\left( \theta \otimes \theta \right) =\left( \theta \otimes
\theta \right) ^{-1}\,\left( \partial \theta \right) ^{-1}\,\left( \theta
\otimes \theta \right) =\left( \theta \otimes \theta \right) ^{-1}\,\psi
^{-1}\,\left( \theta \otimes \theta \right) .
\end{equation*}
This concludes the proof.
\end{proof}

Now we are ready to prove our claim.

\begin{proposition}
Given a quantum space $\mathcal{A}=\left( \mathbf{A}_{1},\mathbf{A}\right) $
and a couple of counital 2-cocycles $\psi $ and $\varphi $, admissible for $%
\mathcal{A}$ and $\mathcal{A}_{\psi }$, respectively, then $\psi \star
\varphi $ is $\mathcal{A}$-admissible and $\left( \mathcal{A}_{\psi }\right)
_{\varphi }=\mathcal{A}_{\psi \star \varphi }$. In addition $\frak{i}\psi $,
the inverse of $\psi $ under the product $\star $, is admissible for $%
\mathcal{A}_{\psi }$, and accordingly $\left( \mathcal{A}_{\psi }\right) _{%
\frak{i}\psi }=\mathcal{A}$.
\end{proposition}

\begin{proof}
Let $\Psi $ and $\Phi $ be the automorphisms of $\mathbf{A}^{\otimes 2}$
related to $\psi $ and $\varphi $, respectively, i.e. 
\begin{equation*}
\Psi \,\Pi ^{\otimes 2}=\Pi ^{\otimes 2}\,\psi \;\;and\;\;\;\Phi \,\left(
\Pi \,\theta \right) ^{\otimes 2}=\left( \Pi \,\theta \right) ^{\otimes
2}\,\varphi ,
\end{equation*}
if $\psi =\partial \theta $. The last equality implies $\Phi \,\Pi ^{\otimes
2}=\Pi ^{\otimes 2}\,\left( \theta \otimes \theta \right) \,\varphi \,\left(
\theta \otimes \theta \right) ^{-1}$ and 
\begin{equation*}
\Psi \,\Phi \,\Pi ^{\otimes 2}=\Pi ^{\otimes 2}\,\psi \,\left( \theta
\otimes \theta \right) \,\varphi \,\left( \theta \otimes \theta \right)
^{-1}=\Pi ^{\otimes 2}\,\psi \star \varphi .
\end{equation*}
Thus $\psi \star \varphi $ is $\mathcal{A}$-admissible, and since the
product of $\left( \mathcal{A}_{\psi }\right) _{\varphi }$ is $\left(
m\,\Psi \right) \,\Phi =m\,\left( \Psi \,\Phi \right) $, the first claim of
the proposition follows. Now, we must show that the 2-cocycle $\frak{i}\psi
=\partial \left( \theta ^{-1}\right) =\left( \theta \otimes \theta \right)
^{-1}\,\psi ^{-1}\,\left( \theta \otimes \theta \right) $ is admissible for $%
\mathcal{A}_{\psi }$. But 
\begin{equation}
\left( \Pi \,\theta \right) ^{\otimes 2}\,\frak{i}\psi =\Pi ^{\otimes
2}\,\psi ^{-1}\,\left( \theta \otimes \theta \right) =\Psi ^{-1}\,\left( \Pi
\,\theta \right) ^{\otimes 2},  \label{dis}
\end{equation}
hence $\frak{i}\psi \left( \ker \left( \Pi \,\theta \right) ^{\otimes
2}\right) \subset \ker \left( \Pi \,\theta \right) ^{\otimes 2}$. Then, the
twisting of $\mathcal{A}_{\psi }$ by $\frak{i}\psi $ can be defined, and
from the above result $\left( \mathcal{A}_{\psi }\right) _{\frak{i}\psi }=%
\mathcal{A}_{\psi \star \frak{i}\psi }=\mathcal{A}$. In particular, the
multiplication of $\left( \mathcal{A}_{\psi }\right) _{\frak{i}\psi }$ will
be, since Eq. $\left( \ref{dis}\right) $, $\left( m\,\Psi \right) \,\Psi
^{-1}=m\,\left( \Psi \,\Psi ^{-1}\right) =m$.
\end{proof}

Let us recall that something similar happens in the bialgebra case. It is
well-known that if $\chi ,\zeta $ are 2-cocycles in $\mathsf{G}^{2}\left[ 
\mathbf{A}\right] $ and $\mathsf{G}^{2}\left[ \mathbf{A}_{\chi }\right] $,
respectively, then $\zeta \ast \chi $ is a 2-cocycle in $\mathsf{G}^{2}\left[
\mathbf{A}\right] $, and $\left( \mathbf{A}_{\chi }\right) _{\zeta }=\mathbf{%
A}_{\zeta \ast \chi }$. In addition, $\chi ^{-1}$ is a 2-cocycle in $\mathsf{%
G}^{2}\left[ \mathbf{A}_{\chi }\right] $, therefore $\left( \mathbf{A}_{\chi
}\right) _{\chi ^{-1}}=\mathbf{A}$. The main difference is that in $\mathsf{G%
}^{2}\left[ \mathbf{A}\right] $ the convolution does not define a subgroup
structure for its 2-cocycles.

\subsubsection{Gauge transformations}

Let us call $\mathbb{T}_{\psi }$ the twist transformation related to $\psi
\in \frak{Z}^{2}$ over a given quantum space. The results of the previous
subsection imply the twist transformations, associated to an inclusion of
vector spaces $\mathbf{A}_{1}\subset \mathbf{A}$, define a groupoid with
(partial) composition $\mathbb{T}_{\psi }\,\mathbb{T}_{\varphi }=\mathbb{T}%
_{\psi \star \varphi }$ and inversion $\mathbb{T}_{\psi }^{-1}=\mathbb{T}_{%
\frak{i}\psi }$. As a generalization of this fact, let us consider the
following definition.

\begin{definition}
We shall say a couple of quantum spaces $\mathcal{A}$ and $\mathcal{B}$ are 
\textbf{twist or gauge related}, namely $\mathcal{A}\backsim \mathcal{B}$,%
\textbf{\ }if there exists an $\mathcal{A}$-admissible $\psi \in \frak{Z}^{2}%
\left[ \mathbf{A}_{1}\right] $ and an isomorphism of quantum spaces $\alpha :%
\mathcal{A}_{\psi }\backsimeq \mathcal{B}$. The pairs $\left( \alpha ,\psi
\right) $ will be called \textbf{gauge transformations} between $\mathcal{A}$
and $\mathcal{B}$.\ \ \ $\blacksquare $
\end{definition}

\begin{proposition}
The twist relation is an equivalence relation.
\end{proposition}

\begin{proof}
\emph{Reflexivity:\ }Since $\mathbb{I}^{\otimes 2}\in \frak{Z}^{2}\left[ 
\mathbf{A}_{1}\right] $ is trivially $\mathcal{A}$-admissible, the identity
morphism insures $\mathcal{A}\backsim \mathcal{A}$.

\emph{Symmetry:} Suppose $\mathcal{A}\backsim \mathcal{B}$ \emph{via} a
cochain $\psi $ and an isomorphism $\alpha $. Let $\theta $ be the primitive
of $\psi $, and consider the 2-cocycle $\frak{i}\psi =\left( \theta \otimes
\theta \right) ^{-1}\,\psi ^{-1}\,\left( \theta \otimes \theta \right) \in 
\frak{Z}^{2}\left[ \mathbf{A}_{1}\right] $. Since the restriction $\alpha
_{1}=\left. \alpha \right| _{\mathbf{A}_{1}}$ defines an isomorphism $%
\mathbf{A}_{1}\backsimeq \mathbf{B}_{1}$, through the functor given in \S 
\textbf{2.3} we have a bijection $\frak{Z}^{2}\left[ \mathbf{A}_{1}\right]
\backsimeq \frak{Z}^{2}\left[ \mathbf{B}_{1}\right] $, such that 
\begin{equation*}
\psi \mapsto \psi ^{\alpha _{1}}=\left( \alpha _{1}^{\otimes }\otimes \alpha
_{1}^{\otimes }\right) \,\psi \,\left( \alpha _{1}^{\otimes }\otimes \alpha
_{1}^{\otimes }\right) ^{-1}.
\end{equation*}
In particular, $\frak{i}\psi ^{\alpha _{1}}=\left( \alpha _{1}^{\otimes
}\otimes \alpha _{1}^{\otimes }\right) \,\frak{i}\psi \,\left( \alpha
_{1}^{\otimes }\otimes \alpha _{1}^{\otimes }\right) ^{-1}\in \frak{Z}^{2}%
\left[ \mathbf{B}_{1}\right] $. Using that $\alpha \,\Pi _{\mathbf{A}%
}\,\theta =\Pi _{\mathbf{B}}\,\alpha _{1}^{\otimes }$ and $\Psi \,\Pi _{%
\mathbf{A}}^{\otimes 2}=\Pi _{\mathbf{A}}^{\otimes 2}\,\psi $, we have 
\begin{equation*}
\Pi _{\mathbf{B}}^{\otimes 2}\,\frak{i}\psi ^{\alpha _{1}}=\left( \alpha
\otimes \alpha \right) \,\Psi ^{-1}\,\left( \alpha \otimes \alpha \right)
^{-1}\,\Pi _{\mathbf{B}}^{\otimes 2}.
\end{equation*}
Thus, $\frak{i}\psi ^{\alpha _{1}}$ is $\mathcal{B}$-admissible and defines
the twisted product 
\begin{equation*}
m_{\mathbf{B}}\,\left( \alpha \otimes \alpha \right) \,\Psi ^{-1}\,\left(
\alpha \otimes \alpha \right) ^{-1}.
\end{equation*}
On the other hand, since $m_{\mathbf{B}}\,\left( \alpha \otimes \alpha
\right) =\alpha \,m_{\mathbf{A}}\,\Psi $, 
\begin{equation*}
\alpha ^{-1}\,\left[ m_{\mathbf{B}}\,\left( \alpha \otimes \alpha \right)
\,\Psi ^{-1}\,\left( \alpha \otimes \alpha \right) ^{-1}\right] =m_{\mathbf{A%
}}\,\left( \alpha \otimes \alpha \right) ^{-1},
\end{equation*}
therefore $\alpha ^{-1}$ defines a quantum space isomorphism $\mathcal{B}_{%
\frak{i}\psi ^{\alpha _{1}}}\backsimeq \mathcal{A}$, and consequently $%
\mathcal{B}\backsim \mathcal{A}$.

\emph{Transitivity: }Consider another quantum space $\mathcal{C}$ such that $%
\mathcal{B}\backsim \mathcal{C}$ \emph{via }a cochain $\varphi $ and an
isomorphism $\beta $. Using again $\alpha \,\Pi _{\mathbf{A}}\,\theta =\Pi _{%
\mathbf{B}}\,\alpha _{1}^{\otimes }$ and $\Psi \,\Pi _{\mathbf{A}}^{\otimes
2}=\Pi _{\mathbf{A}}^{\otimes 2}\,\psi $, it can see that 
\begin{equation*}
\psi \star \varphi ^{\alpha _{1}^{-1}}=\psi \,\left( \theta \otimes \theta
\right) \,\left( \alpha _{1}^{\otimes }\otimes \alpha _{1}^{\otimes }\right)
^{-1}\,\varphi \,\left( \alpha _{1}^{\otimes }\otimes \alpha _{1}^{\otimes
}\right) \,\left( \theta \otimes \theta \right) ^{-1}
\end{equation*}
is $\mathcal{A}$-admissible and has a related automorphism $\Psi \,\left(
\alpha \otimes \alpha \right) ^{-1}\,\Phi \,\left( \alpha \otimes \alpha
\right) $, if $\Phi \,\Pi _{\mathbf{B}}^{\otimes 2}=\Pi _{\mathbf{B}%
}^{\otimes 2}\,\varphi $. It defines on $\mathcal{A}$ the twisted product $%
m_{\mathbf{A}}\,\Psi \,\left( \alpha \otimes \alpha \right) ^{-1}\,\Phi
\,\left( \alpha \otimes \alpha \right) $. Then, 
\begin{eqnarray*}
\beta \,\alpha \,\left[ m_{\mathbf{A}}\,\Psi \,\left( \alpha \otimes \alpha
\right) ^{-1}\,\Phi \,\left( \alpha \otimes \alpha \right) \right] =\beta
\,m_{\mathbf{B}}\,\Phi \,\left( \alpha \otimes \alpha \right) \\
=m_{\mathbf{C}}\,\left( \beta \,\alpha \otimes \beta \,\alpha \right) ,
\end{eqnarray*}
where we have used $\beta \,m_{\mathbf{B}}\,\Phi =m_{\mathbf{C}}\,\left(
\beta \,\otimes \beta \,\right) $. It follows that $\mathcal{A}\backsim 
\mathcal{C}$.
\end{proof}

Denoting by $\mathbb{T}_{\left( \alpha ,\psi \right) }$ the gauge
transformation defined by $\left( \alpha ,\psi \right) $, the theorem above
says those transformations form a groupoid (or a category based on quantum
spaces) with composition and inverse 
\begin{equation*}
\mathbb{T}_{\left( \beta ,\varphi \right) }\,\mathbb{T}_{\left( \alpha ,\psi
\right) }\doteq \mathbb{T}_{\left( \beta ,\varphi \right) \star \left(
\alpha ,\psi \right) },\;\;\;\mathbb{T}_{\left( \alpha ,\psi \right)
}^{-1}\doteq \mathbb{T}_{\frak{i}\left( \alpha ,\psi \right) },
\end{equation*}
being 
\begin{equation*}
\left( \beta ,\varphi \right) \star \left( \alpha ,\psi \right) \doteq
\left( \beta \,\alpha ,\psi \star \varphi ^{\alpha _{1}^{-1}}\right)
\;\;\;and\;\;\;\frak{i}\left( \alpha ,\psi \right) \doteq \left( \alpha
^{-1},\frak{i}\psi ^{\alpha _{1}}\right) .
\end{equation*}
As examples of gauge equivalence we have $\underline{hom}^{\Upsilon }\left[ 
\mathcal{B},\mathcal{A}\right] \backsim \underline{hom}\left[ \mathcal{B},%
\mathcal{A}\right] $ for every pair $\mathcal{B},\mathcal{A}$ in $\mathrm{CA}
$, and $\mathcal{A}\circ _{\tau }\mathcal{B}\backsim \mathcal{A}\circ 
\mathcal{B}$ for every STTP of $\mathcal{A}$ and $\mathcal{B}$.

\bigskip

Geometrically, we are defining an equivalence relation among non commutative
algebraic varieties in $\mathsf{QLS}$.

A characterization of gauge equivalence between conic quantum spaces can be
given in terms of their corresponding ideals. More precisely,

\begin{theorem}
Let $\mathcal{A},\mathcal{B}$ be objects of $\mathrm{CA}$, with related
ideals $\mathbf{I}$ and $\mathbf{J}$, respectively. $\mathcal{A}\backsim 
\mathcal{B}$ \emph{iff }there exists an homogeneous \emph{(}of degree cero%
\emph{)} linear isomorphism $\vartheta :\mathbf{A}_{1}^{\otimes }\backsimeq 
\mathbf{B}_{1}^{\otimes }$ such that 
\begin{equation}
\vartheta \left( \mathbf{A}_{1}^{\otimes r}\cdot \mathbf{I}_{s}+\mathbf{I}%
_{r}\cdot \mathbf{A}_{1}^{\otimes s}\right) =\mathbf{B}_{1}^{\otimes r}\cdot 
\mathbf{J}_{s}+\mathbf{J}_{r}\cdot \mathbf{B}_{1}^{\otimes s};\;\;r,s\in 
\mathbb{N}_{0}.  \label{char}
\end{equation}
\end{theorem}

\begin{proof}
Suppose $\mathcal{A}$ is gauge equivalent to $\mathcal{B}$, i.e. there
exists an $\mathcal{A}$-admissible $\psi \in \frak{Z}^{2}\left[ \mathbf{A}%
_{1}\right] $ and an isomorphism of quantum spaces $\alpha :\mathcal{A}%
_{\psi }\backsimeq \mathcal{B}$. If $\psi =\partial \theta $, from \textbf{%
Lemma 3}, admissibility condition is equivalent to 
\begin{equation}
\theta ^{-1}\left( \mathbf{A}_{1}^{\otimes r}\cdot \mathbf{I}_{s}+\mathbf{I}%
_{r}\cdot \mathbf{A}_{1}^{\otimes s}\right) =\mathbf{A}_{1}^{\otimes r}\cdot
\theta ^{-1}\left( \mathbf{I}_{s}\right) +\theta ^{-1}\left( \mathbf{I}%
_{r}\right) \cdot \mathbf{A}_{1}^{\otimes s}.  \label{char2}
\end{equation}
Because $\alpha $ defines an algebra homomorphism $\mathbf{A}_{\psi
}\backsimeq \mathbf{B}$, that is to say 
\begin{equation*}
\alpha _{1}^{\otimes }\left( \mathbf{I}_{\psi ,r}\right) =\alpha
_{1}^{\otimes }\,\theta ^{-1}\left( \mathbf{I}_{r}\right) =\mathbf{J}_{r},
\end{equation*}
we have 
\begin{equation}
\begin{array}{l}
\alpha _{1}^{\otimes }\,\theta ^{-1}\left( \mathbf{A}_{1}^{\otimes r}\cdot 
\mathbf{I}_{s}+\mathbf{I}_{r}\cdot \mathbf{A}_{1}^{\otimes s}\right) =\alpha
_{1}^{\otimes }\left( \mathbf{A}_{1}^{\otimes r}\cdot \theta ^{-1}\left( 
\mathbf{I}_{s}\right) +\theta ^{-1}\left( \mathbf{I}_{r}\right) \cdot 
\mathbf{A}_{1}^{\otimes s}\right) \\ 
\\ 
=\alpha _{1}^{\otimes }\left( \mathbf{A}_{1}^{\otimes r}\right) \cdot \alpha
_{1}^{\otimes }\,\theta ^{-1}\left( \mathbf{I}_{s}\right) +\alpha
_{1}^{\otimes }\,\theta ^{-1}\left( \mathbf{I}_{r}\right) \cdot \alpha
_{1}^{\otimes }\left( \mathbf{A}_{1}^{\otimes s}\right) \\ 
\\ 
=\mathbf{B}_{1}^{\otimes r}\cdot \mathbf{J}_{s}+\mathbf{J}_{r}\cdot \mathbf{B%
}_{1}^{\otimes s},
\end{array}
\label{char3}
\end{equation}
The linear map $\vartheta =\alpha _{1}^{\otimes }\theta ^{-1}:$ $\mathbf{A}%
_{1}^{\otimes }\rightarrow \mathbf{B}_{1}^{\otimes }$ is obviously the map
we are looking for.

Reciprocally, consider an homogeneous linear bijection $\vartheta :\mathbf{A}%
_{1}^{\otimes }\backsimeq \mathbf{B}_{1}^{\otimes }$ satisfying $\left( \ref
{char}\right) $. Given a basis $\left\{ a_{i}\right\} $ of $\mathbf{A}_{1}$,
if $\vartheta \left( a_{i}\right) =b_{i}$, then $\left\{ b_{i}\right\} $ is
a basis of $\mathbf{B}_{1}$, and we can write 
\begin{equation*}
\vartheta \left( a_{i_{1}}...a_{i_{n}}\right) =\left( \vartheta _{n}\right)
_{i_{1}...i_{n}}^{j_{1}...j_{n}}\,b_{j_{1}}...b_{j_{n}};\;\vartheta \left(
1\right) =1.
\end{equation*}
The assignment $\alpha _{1}=\left. \vartheta \right| _{\mathbf{A}%
_{1}}:a_{i}\mapsto b_{i}$ defines a linear bijection, and $\vartheta $ can
be decomposed as $\vartheta =\alpha _{1}^{\otimes }\theta ^{-1}$ with 
\begin{equation*}
\theta ^{-1}\left( a_{i_{1}}...a_{i_{n}}\right) =\left( \vartheta
_{n}\right) _{i_{1}...i_{n}}^{j_{1}...j_{n}}\,a_{j_{1}}...a_{j_{n}};\;\theta
\left( 1\right) =1.
\end{equation*}
In particular, since $\left( \vartheta _{1}\right) _{i}^{j}=\delta _{i}^{j}$
and $\theta \left( 1\right) =1$, then $\theta _{0,1}=\mathbb{I}_{0,1}$ and,
as a consequence, $\theta \in \frak{P}^{1}\left[ \mathbf{A}_{1}\right] $.
Let us consider the counital 2-cocycle $\psi =\partial \theta $. We shall
see $\psi $ is $\mathcal{A}$-admissible and $\alpha _{1}^{\otimes }$ defines
an isomorphism $\mathcal{A}_{\psi }\backsimeq \mathcal{B}$.

Multiplying Eq. $\left( \ref{char}\right) $ by $\left( \alpha _{1}^{\otimes
}\right) ^{-1}$ to the left, we obtain precisely Eq. $\left( \ref{char2}%
\right) $. But this is the admissibility condition for $\psi $, provided $%
\psi _{r,s}\thickapprox \theta _{r+s}\,\left( \theta _{r}\otimes \theta
_{s}\right) ^{-1}$. Now, using again $\left( \ref{char}\right) $, but in the
form $\left( \ref{char3}\right) $ and fixing $s=0$, equation $\alpha
_{1}^{\otimes }\left( \mathbf{I}_{\psi ,r}\right) =\mathbf{J}_{r}$ follows.
\end{proof}

In the general case we can say,

\begin{proposition}
Given $\mathcal{A},\mathcal{B}\in \mathrm{FGA}$, with related ideals $%
\mathbf{I}$ and $\mathbf{J}$, if $\mathcal{A}\backsim \mathcal{B}$, then
there exists an homogeneous \emph{(}of degree cero\emph{)} linear bijection $%
\vartheta :\mathbf{A}_{1}^{\otimes }\backsimeq \mathbf{B}_{1}^{\otimes }$
such that $\vartheta \left( \mathbf{I}\right) =\mathbf{J}$.
\end{proposition}

\begin{proof}
If $\mathcal{A}\backsim \mathcal{B}$ through $\psi =\partial \theta \in 
\frak{Z}^{2}\left[ \mathbf{A}_{1}\right] $ and an isomorphism $\alpha :%
\mathcal{A}_{\psi }\backsimeq \mathcal{B}$, then $\vartheta =\alpha
_{1}^{\otimes }\,\theta ^{-1}$ is an homogeneous linear isomorphism such
that $\vartheta \left( \mathbf{I}\right) =\mathbf{J}$.
\end{proof}

Thus, gauge equivalent non commutative spaces in $\mathsf{QLS}$ have
isomorphic defining ideals.

\section{Twisting and functors on $\mathrm{CA}$}

In this chapter we mainly deal with conic quantum spaces, analyzing the
relationship between certain functors defined on their category, and
quasicomplexes naturally related to them. For instance, $\frak{C}^{\bullet }%
\left[ \mathbf{A}_{1}\right] $ and $\frak{C}^{\bullet }\left[ \mathbf{A}%
_{1}^{\ast }\right] $ are quasicomplexes naturally related to the functor $!:%
\mathcal{A}\mapsto \mathcal{A}^{!}$, while $\frak{C}^{\bullet }\left[ 
\mathbf{A}_{1}\right] $, $\frak{C}^{\bullet }\left[ \mathbf{B}_{1}\right] $
and $\frak{C}^{\bullet }\left[ \mathbf{A}_{1}\otimes \mathbf{B}_{1}\right] $
are related to bifunctors $\circ $, $\bullet $ and $\odot $. We shall study
the following problems. Suppose a twisting $\mathcal{A}_{\varphi }$ of $%
\mathcal{A}$ is given (of course, $\varphi \in \frak{Z}^{2}\left[ \mathbf{A}%
_{1}\right] $): when can we insure there exists $\chi \in \frak{Z}^{2}\left[ 
\mathbf{A}_{1}^{\ast }\right] $ such that $\left( \mathcal{A}_{\varphi
}\right) ^{!}=\left( \mathcal{A}^{!}\right) _{\chi }$ ? \ On the other hand,
given another twisted quantum space, namely $\mathcal{B}_{\phi }$ (now $\phi
\in \frak{Z}^{2}\left[ \mathbf{B}_{1}\right] $): can we write 
\begin{equation*}
\mathcal{A}_{\varphi }\bigcirc \mathcal{B}_{\phi }=\left( \mathcal{A}%
\bigcirc \mathcal{B}\right) _{\kappa },\;\;\bigcirc =\circ ,\bullet
\;or\;\odot ,
\end{equation*}
for some $\kappa \in \frak{Z}^{2}\left[ \mathbf{A}_{1}\otimes \mathbf{B}_{1}%
\right] $ ?

In the following, we answer these questions and other corresponding to the
functors $\triangleleft $, $\triangleright $ and $\diamond $. In relation to
latter functors, we study the internal coHom objects of twisted quantum
spaces, comparing them with the untwisted ones.

\subsection{2nd admissibility condition}

For conic quantum spaces, a stronger condition than admissibility is
defined. Its consequences are analyzed along all subsequent subsections.

\begin{definition}
Let $\mathcal{A}$ be a conic quantum space with related ideal $\mathbf{I}%
=\bigoplus_{n\geq 2}\mathbf{I}_{n}$. We shall call \textbf{2nd }$\mathcal{A}$%
\textbf{-admissible}, or \textbf{2nd admissible} for $\mathcal{A}$, the $n$%
-cochains satisfying 
\begin{equation}
\begin{array}{l}
\psi \left( \mathbf{A}_{1}^{\otimes p_{1}}\otimes ...\otimes \left( \mathbf{A%
}_{1}^{\otimes p_{i}}\cdot \mathbf{I}_{q}\cdot \mathbf{A}_{1}^{\otimes
r}\right) \otimes ...\otimes \mathbf{A}_{1}^{\otimes p_{n}}\right) \\ 
\\ 
=\mathbf{A}_{1}^{\otimes p_{1}}\otimes ...\otimes \left( \mathbf{A}%
_{1}^{\otimes p_{i}}\cdot \mathbf{I}_{q}\cdot \mathbf{A}_{1}^{\otimes
r}\right) \otimes ...\otimes \mathbf{A}_{1}^{\otimes p_{n}},
\end{array}
\label{2nd}
\end{equation}
\ \ \ for all $p_{k},q,r\in \mathbb{N}_{0}$, $k,i=1...n$.$%
\;\;\;\;\blacksquare $
\end{definition}

For 2-cochains, 2nd admissibility is equivalent to the inclusions 
\begin{equation}
\begin{array}{l}
\psi \left( \mathbf{A}_{1}^{\otimes s}\otimes \left( \mathbf{A}_{1}^{\otimes
p}\cdot \mathbf{I}_{q}\cdot \mathbf{A}_{1}^{\otimes r}\right) \right)
\subset \mathbf{A}_{1}^{\otimes s}\otimes \left( \mathbf{A}_{1}^{\otimes
p}\cdot \mathbf{I}_{q}\cdot \mathbf{A}_{1}^{\otimes r}\right) , \\ 
\\ 
\psi \left( \left( \mathbf{A}_{1}^{\otimes p}\cdot \mathbf{I}_{q}\cdot 
\mathbf{A}_{1}^{\otimes r}\right) \otimes \mathbf{A}_{1}^{\otimes s}\right)
\subset \left( \mathbf{A}_{1}^{\otimes p}\cdot \mathbf{I}_{q}\cdot \mathbf{A}%
_{1}^{\otimes r}\right) \otimes \mathbf{A}_{1}^{\otimes s};
\end{array}
\label{2a}
\end{equation}
and \textbf{Lemma 3} translates into the next result.

\begin{lemma}
Under the conditions of above definition, $\psi =\partial \theta $ is 2nd
admissible \emph{iff} 
\begin{equation}
\theta ^{-1}\left( \mathbf{A}_{1}^{\otimes p}\cdot \mathbf{I}_{q}\cdot 
\mathbf{A}_{1}^{\otimes r}\right) =\mathbf{A}_{1}^{\otimes p}\cdot \theta
^{-1}\left( \mathbf{I}_{q}\right) \cdot \mathbf{A}_{1}^{\otimes r}.
\label{2at}
\end{equation}
\end{lemma}

\begin{proof}
If $\left( \ref{2at}\right) $ holds for the primitive $\theta $ of $\psi $,
we have, putting $p=x+y$, 
\begin{equation*}
\theta ^{-1}\left( \mathbf{A}_{1}^{\otimes x}\cdot \mathbf{A}_{1}^{\otimes
y}\cdot \mathbf{I}_{q}\cdot \mathbf{A}_{1}^{\otimes r}\right) =\mathbf{A}%
_{1}^{\otimes x}\cdot \mathbf{A}_{1}^{\otimes y}\cdot \theta ^{-1}\left( 
\mathbf{I}_{q}\right) \cdot \mathbf{A}_{1}^{\otimes r}.
\end{equation*}
But $\mathbf{A}_{1}^{\otimes y}\cdot \theta ^{-1}\left( \mathbf{I}%
_{q}\right) \cdot \mathbf{A}_{1}^{\otimes r}=\theta ^{-1}\left( \mathbf{A}%
_{1}^{\otimes y}\cdot \mathbf{I}_{q}\cdot \mathbf{A}_{1}^{\otimes r}\right) $
and $\mathbf{A}_{1}^{\otimes x}=\theta ^{-1}\left( \mathbf{A}_{1}^{\otimes
x}\right) $, thus 
\begin{equation*}
\theta ^{-1}\left( \mathbf{A}_{1}^{\otimes x}\cdot \left( \mathbf{A}%
_{1}^{\otimes y}\cdot \mathbf{I}_{q}\cdot \mathbf{A}_{1}^{\otimes r}\right)
\right) =\theta ^{-1}\left( \mathbf{A}_{1}^{\otimes x}\right) \cdot \theta
^{-1}\left( \mathbf{A}_{1}^{\otimes y}\cdot \mathbf{I}_{q}\cdot \mathbf{A}%
_{1}^{\otimes r}\right) ,
\end{equation*}
and using the fact that $\psi =\partial \theta $, the first part of $\left( 
\ref{2a}\right) $ follows. In a similar way the second part can be shown,
and 2nd admissibility is fulfilled. The converse follows in a similar way.
\end{proof}

Clearly, if $\psi $ is 2nd admissible, then is admissible. Furthermore, the
2nd admissible cochains form a subgroup of the admissible ones.

\subsubsection{Restricted gauge equivalence}

For 2-cocycles we have, in relation to the product $\star $ in $\frak{Z}^{2}$%
, that:

\begin{proposition}
Let $\psi $ and $\varphi $ be a couple of counital 2-cocycles, 2nd
admissible for $\mathcal{A}$ and $\mathcal{A}_{\psi }$, respectively. Then $%
\psi \star \varphi $ is 2nd $\mathcal{A}$-admissible. Furthermore $\frak{i}%
\psi $, the inverse of $\psi $ under the product $\star $, is 2nd admissible
for $\mathcal{A}_{\psi }$.
\end{proposition}

\begin{proof}
Suppose $\psi =\partial \lambda $ and $\varphi =\partial \chi $. From lemma
above, $\lambda $ and $\chi $ satisfies 
\begin{equation}
\lambda ^{-1}\left( \mathbf{A}_{1}^{\otimes p}\cdot \mathbf{I}_{q}\cdot 
\mathbf{A}_{1}^{\otimes r}\right) =\mathbf{A}_{1}^{\otimes p}\cdot \lambda
^{-1}\left( \mathbf{I}_{q}\right) \cdot \mathbf{A}_{1}^{\otimes r}  \label{l}
\end{equation}
and 
\begin{equation*}
\chi ^{-1}\left( \mathbf{A}_{1}^{\otimes p}\cdot \lambda ^{-1}\left( \mathbf{%
I}_{q}\right) \cdot \mathbf{A}_{1}^{\otimes r}\right) =\mathbf{A}%
_{1}^{\otimes p}\cdot \left( \lambda \,\chi \right) ^{-1}\left( \mathbf{I}%
_{q}\right) \cdot \mathbf{A}_{1}^{\otimes r},
\end{equation*}
since $\lambda ^{-1}\left( \mathbf{I}\right) $ is the ideal related to $%
\mathcal{A}_{\psi }$. In consequence, 
\begin{equation*}
\left( \lambda \,\chi \right) ^{-1}\left( \mathbf{A}_{1}^{\otimes p}\cdot 
\mathbf{I}_{q}\cdot \mathbf{A}_{1}^{\otimes r}\right) =\mathbf{A}%
_{1}^{\otimes p}\cdot \left( \lambda \,\chi \right) ^{-1}\left( \mathbf{I}%
_{q}\right) \cdot \mathbf{A}_{1}^{\otimes r}.
\end{equation*}
On the other hand, $\psi \star \varphi =\partial \lambda \star \partial \chi
=\partial \left( \lambda \,\chi \right) $, hence $\psi \star \varphi $ has a
primitive satisfying $\left( \ref{2at}\right) $. Therefore, $\psi \star
\varphi $ is 2nd admissible for $\mathcal{A}$.

For the second claim of our proposition, using that $\frak{i}\psi =\partial
\left( \lambda ^{-1}\right) $ and (from Eq. $\left( \ref{l}\right) $) 
\begin{eqnarray*}
\lambda \left( \mathbf{A}_{1}^{\otimes p}\cdot \lambda ^{-1}\left( \mathbf{I}%
_{q}\right) \cdot \mathbf{A}_{1}^{\otimes r}\right) &=&\lambda \left(
\lambda ^{-1}\left( \mathbf{A}_{1}^{\otimes p}\cdot \mathbf{I}_{q}\cdot 
\mathbf{A}_{1}^{\otimes r}\right) \right) \\
&& \\
&=&\mathbf{A}_{1}^{\otimes p}\cdot \mathbf{I}_{q}\cdot \mathbf{A}%
_{1}^{\otimes r}=\mathbf{A}_{1}^{\otimes p}\cdot \lambda \left( \lambda
^{-1}\left( \mathbf{I}_{q}\right) \right) \cdot \mathbf{A}_{1}^{\otimes r},
\end{eqnarray*}
the 2nd $\mathcal{A}_{\psi }$-admissibility of $\frak{i}\psi $ is immediate.
\end{proof}

Therefore, in terms of 2nd admissible cocycles a restricted twist
equivalence can be defined. A characterization of equivalent quantum spaces
through this gauge transformations follow from the proof of \textbf{Theor. 12%
}, and says:

\begin{theorem}
Let $\mathcal{A},\mathcal{B}$ be objects of $\mathrm{CA}$, with related
ideals $\mathbf{I}$ and $\mathbf{J}$, respectively. $\mathcal{A}\backsim 
\mathcal{B}$ \emph{iff }there exists an homogeneous \emph{(}of degree cero%
\emph{)} linear isomorphism $\vartheta :\mathbf{A}_{1}^{\otimes }\backsimeq 
\mathbf{B}_{1}^{\otimes }$ such that 
\begin{equation*}
\vartheta \left( \mathbf{A}_{1}^{\otimes r}\cdot \mathbf{I}_{s}\cdot \mathbf{%
A}_{1}^{\otimes t}\right) =\mathbf{B}_{1}^{\otimes r}\cdot \mathbf{J}%
_{s}\cdot \mathbf{B}_{1}^{\otimes t};\;\;r,s,t\in \mathbb{N}%
_{0}.\;\;\;\blacksquare
\end{equation*}
\end{theorem}

By last, consider the following property of 2nd admissible 2-cocycles.

\begin{proposition}
If a pair of 2nd admissible 2-cocycles are cohomologous through an
admissible 1-cochain $\theta $, then $\theta $ is 2nd admissible.
\end{proposition}

\begin{proof}
Let $\psi =\partial \lambda $ and $\varphi =\partial \chi $ be cohomologous
2-cocycles through an admissible cochain $\theta $. We know from \textbf{%
Prop.} \textbf{11 }that $\theta =\lambda \,\omega \,\chi ^{-1}$, with $%
\omega $ a 1-cocycle. Using \textbf{Lemma 4} for $\chi $ and then the
admissibility of $\theta $ (and the fact that $\omega $ is a 1-cocycle), 
\begin{equation*}
\begin{array}{l}
\lambda \,\omega \,\chi ^{-1}\left( \mathbf{A}_{1}^{\otimes p}\cdot \mathbf{I%
}_{q}\cdot \mathbf{A}_{1}^{\otimes r}\right) =\lambda \left( \mathbf{A}%
_{1}^{\otimes p}\cdot \omega \,\chi ^{-1}\left( \mathbf{I}_{q}\right) \cdot 
\mathbf{A}_{1}^{\otimes r}\right) \\ 
\\ 
=\lambda \left( \mathbf{A}_{1}^{\otimes p}\cdot \lambda ^{-1}\,\left(
\lambda \,\omega \,\chi ^{-1}\right) \left( \mathbf{I}_{q}\right) \cdot 
\mathbf{A}_{1}^{\otimes r}\right) =\lambda \left( \mathbf{A}_{1}^{\otimes
p}\cdot \lambda ^{-1}\left( \mathbf{I}_{q}\right) \cdot \mathbf{A}%
_{1}^{\otimes r}\right) .
\end{array}
\end{equation*}
Using same lemma for $\lambda $ we have $\theta \left( \mathbf{A}%
_{1}^{\otimes p}\cdot \mathbf{I}_{q}\cdot \mathbf{A}_{1}^{\otimes r}\right) =%
\mathbf{A}_{1}^{\otimes p}\cdot \mathbf{I}_{q}\cdot \mathbf{A}_{1}^{\otimes
r}$, as we have claimed.
\end{proof}

From that, \textbf{Theor. 9} can be rephrased in terms of 2nd admissibility
too.

\subsubsection{The \emph{(}anti\emph{)}bicharacter case}

Examples of 2nd admissible cochains are the admissible (anti)bicharacters.

\begin{proposition}
Let $\psi \in \frak{C}^{2}\left[ \mathbf{A}_{1}\right] $ be an \emph{(}%
\textbf{anti}\emph{)}\textbf{bicharacter}. $\psi $ is admissible \emph{iff }%
is 2nd admissible.\ \ \ $\blacksquare $
\end{proposition}

That follows from proposition below.

\begin{proposition}
Let $\psi \in \frak{C}^{2}\left[ \mathbf{A}_{1}\right] $ be an \emph{(}anti%
\emph{)}bicharacter\emph{\ }for which there exists a graded vector space 
\begin{equation*}
\mathbf{S}=\bigoplus\nolimits_{n\in \mathbb{N}_{0}}\mathbf{S}_{n}\subset 
\mathbf{A}_{1}^{\otimes }=\bigoplus\nolimits_{n\in \mathbb{N}_{0}}\mathbf{A}%
_{1}^{\otimes n};\;\;\mathbf{S}_{0,1}\doteq \left\{ 0\right\} ,
\end{equation*}
such that 
\begin{equation}
\psi \left( \mathbf{A}_{1}^{\otimes m}\otimes \mathbf{S}_{n}+\mathbf{S}%
_{m}\otimes \mathbf{A}_{1}^{\otimes n}\right) =\mathbf{A}_{1}^{\otimes
m}\otimes \mathbf{S}_{n}+\mathbf{S}_{m}\otimes \mathbf{A}_{1}^{\otimes n}.
\label{uu}
\end{equation}
Then, for all $p,q,r,u\in \mathbb{N}_{0}$, we have 
\begin{equation}
\begin{array}{l}
\psi \left( \left( \mathbf{A}_{1}^{\otimes p}\cdot \mathbf{S}_{q}\cdot 
\mathbf{A}_{1}^{\otimes r}\right) \otimes \mathbf{A}_{1}^{\otimes u}\right)
=\left( \mathbf{A}_{1}^{\otimes p}\cdot \mathbf{S}_{q}\cdot \mathbf{A}%
_{1}^{\otimes r}\right) \otimes \mathbf{A}_{1}^{\otimes u}, \\ 
\\ 
\psi \left( \mathbf{A}_{1}^{\otimes u}\otimes \left( \mathbf{A}_{1}^{\otimes
p}\cdot \mathbf{S}_{q}\cdot \mathbf{A}_{1}^{\otimes r}\right) \right) =%
\mathbf{A}_{1}^{\otimes u}\otimes \left( \mathbf{A}_{1}^{\otimes p}\cdot 
\mathbf{S}_{q}\cdot \mathbf{A}_{1}^{\otimes r}\right) ;
\end{array}
\label{iq}
\end{equation}
and if $\psi $ is a 2-cocycle with primitive $\theta $, 
\begin{equation}
\theta ^{-1}\left( \mathbf{A}_{1}^{\otimes p}\cdot \mathbf{S}_{q}\cdot 
\mathbf{A}_{1}^{\otimes r}\right) =\mathbf{A}_{1}^{\otimes p}\cdot \theta
^{-1}\left( \mathbf{S}_{q}\right) \cdot \mathbf{A}_{1}^{\otimes r}.
\label{iq2}
\end{equation}
\end{proposition}

\begin{proof}
Consider an element $a\in \mathbf{A}_{1}$ and $s\in \mathbf{S}_{n}$. Since
Eq. $\left( \ref{uu}\right) $, 
\begin{equation*}
\psi \left( a\otimes s\right) =\psi _{1,n}\left( a\otimes s\right) \in 
\mathbf{A}_{1}\otimes \mathbf{S}_{n}+\mathbf{S}_{1}\otimes \mathbf{A}%
_{1}^{\otimes n}.
\end{equation*}
But $\mathbf{S}_{1}=\left\{ 0\right\} $, hence 
\begin{equation}
\psi _{1,n}\left( \mathbf{A}_{1}\otimes \mathbf{S}_{n}\right) =\mathbf{A}%
_{1}\otimes \mathbf{S}_{n}\;\;\;or\;\;\psi _{1,n}\left( \mathbf{A}%
_{1}\otimes \mathbf{S}_{n}\right) =\mathbf{A}_{1}\otimes \mathbf{S}_{n}.
\label{n1}
\end{equation}
Now, consider an element $a\otimes b=a\cdot b\in \mathbf{A}_{1}^{\otimes 2}$%
. Suppose $\psi $ is a bicharacter. From the first part of Eq. $\left( \ref
{rst}\right) $ we have 
\begin{eqnarray*}
\psi _{2,n}\left( \left( a\cdot b\right) \otimes s\right) =\psi
_{1+1,n}\left( \left( a\cdot b\right) \otimes s\right) \\
=\left( \mathbb{I}\otimes \psi _{1,n}\right) \,\left( \mathbb{I}\otimes
f_{1,n}^{-1}\right) \,\left( \psi _{1,n}\left( a\otimes s\right) \otimes
b\right) ,
\end{eqnarray*}
and from Eq. $\left( \ref{n1}\right) $ it follows that $\psi \left( \mathbf{A%
}_{1}^{\otimes 2}\otimes \mathbf{S}_{n}\right) =\mathbf{A}_{1}^{\otimes
2}\otimes \mathbf{S}_{n}$. An inductive reasoning shows that 
\begin{equation*}
\psi \left( \mathbf{A}_{1}^{\otimes m}\otimes \mathbf{S}_{n}\right) =\psi
_{m,n}\left( \mathbf{A}_{1}^{\otimes m}\otimes \mathbf{S}_{n}\right) =%
\mathbf{A}_{1}^{\otimes m}\otimes \mathbf{S}_{n}
\end{equation*}
for all $n,m$. Now, consider an element $c\otimes \left( a\cdot s\right) \in 
\mathbf{A}_{1}^{\otimes m}\otimes \mathbf{A}_{1}\cdot \mathbf{S}_{n}$. Using
the second part of Eq. $\left( \ref{rst}\right) $, 
\begin{equation*}
\psi _{m,1+n}\left( c\otimes \left( a\cdot s\right) \right) =\left( \psi
_{m,1}\otimes \mathbb{I}_{n}\right) \,\left( \mathbb{I}_{m}\otimes
f_{1,n}^{-1}\right) \,\left( \psi _{m,n}\left( c\otimes s\right) \otimes
a\right) ,
\end{equation*}
hence $\psi \left( \mathbf{A}_{1}^{\otimes m}\otimes \mathbf{A}_{1}\cdot 
\mathbf{S}_{n}\right) =\mathbf{A}_{1}^{\otimes m}\otimes \mathbf{A}_{1}\cdot 
\mathbf{S}_{n}$. Analogous arguments, followed by inductive reasoning, lead
us finally to the equalities $\left( \ref{iq}\right) $.

To prove the second part of the lemma, i.e. Eq. $\left( \ref{iq2}\right) $,
recall that a 2-cocycle $\psi $ and its primitive $\theta $ are related by
the equation 
\begin{equation*}
\theta _{m+n}^{-1}\thickapprox \left( \theta _{m}^{-1}\otimes \theta
_{n}^{-1}\right) \,\psi _{m,n}^{-1}.
\end{equation*}
Thus, the previous result implies $\theta ^{-1}\left( \mathbf{S}_{m}\cdot 
\mathbf{A}_{1}^{\otimes n}\right) =\theta ^{-1}\left( \mathbf{S}_{m}\right)
\cdot \theta ^{-1}\left( \mathbf{A}_{1}^{\otimes n}\right) $, and the other
equalities follow in a similar way. This ends our proof.
\end{proof}

Each bicharacter $\varsigma _{\mathcal{A}}$ appearing in twisted coHom
objects (see $\left( \ref{sa}\right) $) is admissible, $ipso$ $facto$ 2nd
admissible, thanks to $\sigma _{\mathcal{A}}$ defines an automorphism. In
fact, 
\begin{equation*}
\begin{array}{l}
\varsigma _{\mathcal{A}}\left( \mathbf{A}_{1}^{\otimes s}\otimes \left( 
\mathbf{A}_{1}^{\otimes p}\cdot \mathbf{I}_{q}\cdot \mathbf{A}_{1}^{\otimes
r}\right) \right) =\mathbf{A}_{1}^{\otimes s}\otimes \left[ \sigma _{%
\mathcal{A}}^{-s}\right] ^{\otimes }\left( \mathbf{A}_{1}^{\otimes p}\cdot 
\mathbf{I}_{q}\cdot \mathbf{A}_{1}^{\otimes r}\right) \\ 
\\ 
=\mathbf{A}_{1}^{\otimes s}\otimes \left[ \sigma _{\mathcal{A}}^{-s}\right]
^{\otimes }\left( \mathbf{A}_{1}^{\otimes p}\right) \cdot \left[ \sigma _{%
\mathcal{A}}^{-s}\right] ^{\otimes }\left( \mathbf{I}_{q}\right) \cdot \left[
\sigma _{\mathcal{A}}^{-s}\right] ^{\otimes }\left( \mathbf{A}_{1}^{\otimes
r}\right) =\mathbf{A}_{1}^{\otimes s}\otimes \mathbf{A}_{1}^{\otimes p}\cdot 
\mathbf{I}_{q}\cdot \mathbf{A}_{1}^{\otimes r}.
\end{array}
\end{equation*}

\subsubsection{Quasicomplexes of 2nd admissible cochains}

From Eq. $\left( \ref{2nd}\right) $, easy calculations show that the coface
operators and codegeneracies on $\mathsf{C}^{\bullet }$ send the subgroup of
2nd admissible cochains to itself. Thus, we can define for each conic
quantum space $\mathcal{A}=\left( \mathbf{A}_{1},\mathbf{A}\right) $ a
subquasicomplex $\mathsf{C}^{\bullet }\left[ \mathcal{A}\right] $ of $%
\mathsf{C}^{\bullet }\left[ \mathbf{A}_{1}\right] $ formed out by 2nd $%
\mathcal{A}$-admissible cochains. We are not going to employ explicitly this
quasicomplex in any calculation. However, for the sake of completeness, we
enumerate some of its properties bellow (without proof).

$\bullet $ As in \S \textbf{2.3}, the cosimplicial objects $\left[ n+1\right]
\mapsto \mathsf{C}^{n}\left[ \mathcal{A}\right] $ and $\left[ n+1\right]
\mapsto \mathsf{C}^{n}\left[ \mathcal{B}\right] $ are naturally equivalent
when $\mathcal{A}\backsimeq \mathcal{B}$.

$\bullet $ The map $\mathcal{A}\mapsto \mathsf{C}^{\bullet }\left[ \mathcal{A%
}\right] $ defines a functor $\mathsf{C}^{\bullet }:\mathcal{G}\left[ 
\mathrm{CA}\right] \rightarrow \mathrm{Grp}_{\ast q}$ and a morphism between
the groupoids $\mathcal{G}\left[ \mathrm{CA}\right] $ and $\mathcal{G}\left[ 
\mathrm{Grp}_{\ast q}\right] $ such that to every isomorphism $\alpha :%
\mathcal{A}\backsimeq \mathcal{B}$ it assigns the isomorphism of
quasicomplexes $\alpha ^{\bullet }$ defined as in Eq. $\left( \ref{fiefe}%
\right) $, replacing $f$ by $\alpha _{1}$.

$\bullet $ The subgroups of counital cochains give rise to a subquasicomplex 
$\frak{C}^{\bullet }\left[ \mathcal{A}\right] $ satisfying Eq. $\left( \ref
{eqy}\right) $.

$\bullet $ The first and second cohomology spaces, namely $H^{1,2}\left[ 
\mathcal{A}\right] $, can be defined as in \S \textbf{3.2}. Now, they are
usual cohomological spaces in the sense that are given by quotients $\left. 
\mathsf{Z}^{1,2}\left[ \mathcal{A}\right] \right/ \backsim _{\mathsf{Coh}}$.
It follows that 
\begin{equation*}
H^{1}\left[ \mathcal{A}\right] =\mathsf{Z}^{1}\left[ \mathcal{A}\right]
\backsimeq Aut_{\mathrm{CA}}\left[ \mathcal{A}\right] \backsimeq H_{\mathcal{%
A}}^{1}\left[ \mathbf{A}_{1}\right] ,
\end{equation*}
and from \textbf{Prop. 18} the inclusion $H^{2}\left[ \mathcal{A}\right]
\subset H_{\mathcal{A}}^{2}\left[ \mathbf{A}_{1}\right] $ is immediate.

$\bullet $ From Eq. $\left( \ref{cg}\right) $ (putting $\mathbf{V}=\mathbf{A}%
_{1}$), it follows by restriction the isomorphism 
\begin{equation}
\mathsf{C}^{\bullet }\left[ \mathcal{A}\right] \backsimeq \mathsf{G}%
^{\bullet }\left[ \underline{end}\left[ \mathcal{A}\right] \right] =\mathsf{G%
}^{\bullet }\left[ \mathcal{A}\triangleright \mathcal{A}\right] .
\label{cgad}
\end{equation}
In particular, suppose $\mathcal{A}$ is quadratic and consider the algebras $%
\underline{e}\left[ \mathcal{A}\right] $ and $A\left( R\right) $ defined in 
\cite{man0} and \cite{frt},\footnote{%
If the ideal related to $\mathcal{A}$ is generated by elements $%
R_{ij}^{kl}\,a_{k}a_{l}$, being $R_{ij}^{kl}$ the coefficients of a map $R:%
\mathbf{A}_{1}^{\otimes 2}\rightarrow \mathbf{A}_{1}^{\otimes 2}$, then the
ideal related to $frt\left[ \mathcal{A},R\right] $ is generated by $%
R_{ij}^{kl}\,z_{k}^{n}z_{l}^{m}-z_{i}^{k}z_{j}^{l}\,R_{kl}^{nm}$.}
respectively. The well known associated epimorphisms $A\left( R\right)
\twoheadleftarrow \underline{end}\left[ \mathcal{A}\right]
\twoheadrightarrow \underline{e}\left[ \mathcal{A}\right] $, which follows
from initiality of $\underline{end}\left[ \mathcal{A}\right] $, gives rise
to monics 
\begin{equation*}
\mathsf{G}^{\bullet }\left[ A\left( R\right) \right] \hookrightarrow \mathsf{%
G}^{\bullet }\left[ \underline{end}\left[ \mathcal{A}\right] \right]
\hookleftarrow \mathsf{G}^{\bullet }\left[ \underline{e}\left[ \mathcal{A}%
\right] \right] .
\end{equation*}
Thus, the bialgebra twisting of $A\left( R\right) $ and $\underline{e}\left[ 
\mathcal{A}\right] $ are given by a particular subgroup of 2nd admissible
twisting of the quantum space $\mathcal{A}$.

\smallskip

Anther properties of $\mathsf{C}^{\bullet }\left[ \mathcal{A}\right] $, in
relation with the involution and the product between quasicomplexes in $%
\mathsf{C}^{\bullet }\left[ \mathrm{Vct}\right] $, will be briefly commented
in the following section.

\subsection{Coadjoints, products and internal coHom objects}

\subsubsection{Dual quantum spaces and coadjoint cochains}

Let us first enunciate the following well-known result (without proof).

\begin{lemma}
Consider a pair of finite dimensional $\Bbbk $-vector spaces $\mathbf{U}$
and $\mathbf{V}$, supplied with a non degenerated pairing $\mathbf{U}\times 
\mathbf{V}\rightarrow \Bbbk $. Let $\alpha $ be an automorphism in $\mathbf{V%
}$. Then, for every subset $\mathbf{S}\subset \mathbf{V}$ we have the
equality of vector spaces\footnote{%
As usual, given a subset $\mathbf{S}\subset \mathbf{V}$, $\mathbf{S}^{\perp
} $ is the vector space $\mathbf{S}^{\perp }=\left\{ x\in \mathbf{U}%
:\left\langle x,y\right\rangle =0,\;\forall y\in \mathbf{S}\right\} .$} 
\begin{equation*}
\alpha \left( \mathbf{S}\right) ^{\perp }=\alpha ^{\ast -1}\left( \mathbf{S}%
^{\perp }\right) ,
\end{equation*}
where $\alpha ^{\ast }$ denotes the \emph{(}transpose\emph{)} automorphism
induced in $\mathbf{U}$ by the mentioned pairing.\ \ \ $\blacksquare $
\end{lemma}

Consider a conic quantum space $\mathcal{A}$ and a counital $\mathcal{A}$%
-admissible 2-cocycle $\psi =\partial \theta \in \frak{Z}^{2}\left[ \mathbf{A%
}_{1}\right] $. The ideal $\mathbf{I}_{\psi }$ related to $\mathcal{A}_{\psi
}$ is determined by $\mathbf{I}_{\psi ,n}=\theta ^{-1}\left( \mathbf{I}%
_{n}\right) $, hence (from lemma above) 
\begin{equation*}
\left( \mathbf{I}_{\psi ,n}\right) ^{\perp }=\theta ^{-1}\left( \mathbf{I}%
_{n}\right) ^{\perp }=\theta ^{\ast }\left( \mathbf{I}_{n}^{\perp }\right) ,
\end{equation*}
and recalling Eq. $\left( \ref{idag}\right) $, $\left( \mathbf{I}_{\psi
}\right) ^{\dagger }=I\left[ \bigoplus\nolimits_{n\geq 2}\left( \mathbf{I}%
_{\psi ,n}\right) ^{\perp }\right] =I\left[ \bigoplus\nolimits_{n\geq
2}\theta ^{\ast }\left( \mathbf{I}_{n}^{\perp }\right) \right] $. On the
other hand, consider the 2-cocycle\footnote{%
Recall the isomorphism of quasicomplexes $\mathsf{C}^{\bullet }\left[ 
\mathbf{V}\right] \backsimeq \mathsf{C}^{\bullet }\left[ \mathbf{V}\right]
^{!}=\mathsf{C}^{\bullet }\left[ \mathbf{V}^{\ast }\right] $.} 
\begin{equation*}
\psi ^{\ast -1}=\psi ^{!}=\left( \partial \theta \right) ^{!}=\partial
^{!}\left( \theta ^{!}\right) \in \frak{Z}^{2}\left[ \mathbf{A}_{1}^{\ast }%
\right] =\frak{Z}^{2}\left[ \mathbf{A}_{1}\right] ^{!}.
\end{equation*}
If $\psi ^{!}$ were $\mathcal{A}^{!}$-admissible, it follows from \textbf{%
Prop. 12} (for $\mathbf{S}=\bigoplus\nolimits_{n\geq 2}\mathbf{I}_{n}^{\perp
}$) that the ideal related to $\left( \mathcal{A}^{!}\right) _{\psi ^{!}}$
would be 
\begin{equation*}
\left( \mathbf{I}^{\dagger }\right) _{\psi ^{!}}=\theta ^{\ast }\left( 
\mathbf{I}^{\dagger }\right) =I\left[ \bigoplus\nolimits_{n\geq 2}\theta
^{\ast }\left( \mathbf{I}_{n}^{\perp }\right) \right] .
\end{equation*}
That is to say, $\left( \mathcal{A}_{\psi }\right) ^{!}$ and $\left( 
\mathcal{A}^{!}\right) _{\psi ^{!}}$ coincide. The problem is to insure the
admissibility of $\psi ^{!}$. The things change when we suppose $\psi $ is
2nd admissible or, in particular, if $\psi $ is an admissible
(anti)bicharacter.

\begin{proposition}
If $\psi \in \frak{C}^{n}\left[ \mathbf{A}_{1}\right] $ is 2nd $\mathcal{A}$%
-admissible, then $\psi ^{!}\in \frak{C}^{n}\left[ \mathbf{A}_{1}\right]
^{!} $ is 2nd $\mathcal{A}^{!}$-admissible. For $\mathcal{A}\in \mathrm{CA}%
^{m}$, the converse is also valid.
\end{proposition}

\begin{proof}
We just consider $n=2$, since the other cases can be shown in a similar way.
From \textbf{Lemma} \textbf{5}, 
\begin{eqnarray*}
\psi ^{!}\left( \mathbf{A}_{1}^{\ast \otimes u}\otimes \left( \mathbf{A}%
_{1}^{\ast \otimes p}\cdot \mathbf{I}_{q}^{\perp }\cdot \mathbf{A}_{1}^{\ast
\otimes r}\right) \right) =\psi ^{!}\left( \left( \mathbf{A}_{1}^{\otimes
u}\otimes \left( \mathbf{A}_{1}^{\otimes p}\cdot \mathbf{I}_{q}\cdot \mathbf{%
A}_{1}^{\otimes r}\right) \right) ^{\perp }\right) \\
&& \\
=\psi \left( \mathbf{A}_{1}^{\otimes u}\otimes \left( \mathbf{A}%
_{1}^{\otimes p}\cdot \mathbf{I}_{q}\cdot \mathbf{A}_{1}^{\otimes r}\right)
\right) ^{\perp },
\end{eqnarray*}
and from 2nd $\mathcal{A}$-admissibility of $\psi $, i.e. 
\begin{equation}
\psi \left( \mathbf{A}_{1}^{\otimes u}\otimes \left( \mathbf{A}_{1}^{\otimes
p}\cdot \mathbf{I}_{q}\cdot \mathbf{A}_{1}^{\otimes r}\right) \right) =%
\mathbf{A}_{1}^{\otimes u}\otimes \left( \mathbf{A}_{1}^{\otimes p}\cdot 
\mathbf{I}_{q}\cdot \mathbf{A}_{1}^{\otimes r}\right) ,  \label{ke1}
\end{equation}
we have 
\begin{equation*}
\psi ^{!}\left( \mathbf{A}_{1}^{\ast \otimes u}\otimes \left( \mathbf{A}%
_{1}^{\ast \otimes p}\cdot \mathbf{I}_{q}^{\perp }\cdot \mathbf{A}_{1}^{\ast
\otimes r}\right) \right) =\mathbf{A}_{1}^{\ast \otimes u}\otimes \left( 
\mathbf{A}_{1}^{\ast \otimes p}\cdot \mathbf{I}_{q}^{\perp }\cdot \mathbf{A}%
_{1}^{\ast \otimes r}\right) .
\end{equation*}
Straightforwardly, because $\mathbf{I}_{n}^{\dagger
}=\sum\nolimits_{r=2}^{n}\sum\nolimits_{i=0}^{n-r}\mathbf{A}_{1}^{\ast
\otimes n-r-i}\cdot \mathbf{I}_{r}^{\perp }\cdot \mathbf{A}_{1}^{\ast
\otimes i}$, $\psi ^{!}$ is 2nd $\mathcal{A}^{!}$-admissible. Now, let us
shown the other claim. Suppose $\mathcal{A}\in \mathrm{CA}^{m}$. From the
last result, if $\psi ^{!}\in \frak{C}^{n}\left[ \mathbf{A}_{1}\right] ^{!}$
is 2nd $\mathcal{A}^{!}$-admissible, then $\psi ^{!!}$ is 2nd $\mathcal{A}%
^{!!}$-admissible. But $\psi ^{!!}=\psi $ and the ideal related to $\mathcal{%
A}^{!!}$ is the same as the one related to $\mathcal{A}$. So, the
proposition have been proven.
\end{proof}

In terms of complexes $\mathsf{C}^{\bullet }\left[ \mathcal{A}\right] $,
last proposition says isomorphism $\mathsf{C}^{\bullet }\left[ \mathbf{A}_{1}%
\right] \backsimeq \mathsf{C}^{\bullet }\left[ \mathbf{A}_{1}\right] ^{!}$
gives rise to a monic $\mathsf{C}^{\bullet }\left[ \mathcal{A}\right]
\hookrightarrow \mathsf{C}^{\bullet }\left[ \mathcal{A}^{!}\right] $ which
define a natural transformation $\mathsf{C}^{\bullet }\rightarrow \mathsf{C}%
^{\bullet }\,!$. In addition, when $\mathsf{C}^{\bullet }$ is restricted to $%
\mathrm{CA}^{m}$, the isomorphisms $\mathsf{C}^{\bullet }\left[ \mathcal{A}%
\right] \backsimeq \mathsf{C}^{\bullet }\left[ \mathcal{A}^{!}\right] $ and
the related natural equivalence $\mathsf{C}^{\bullet }\backsimeq \mathsf{C}%
^{\bullet }\,!$ follow.

In resume, we have shown:

\begin{theorem}
Let $\mathcal{A}$ be an object in $\mathrm{CA}$. For every counital 2nd $%
\mathcal{A}$-admissible $\psi \in \frak{Z}^{2}\left[ \mathbf{A}_{1}\right] $%
, $\psi ^{!}$ is a 2nd $\mathcal{A}^{!}$-admissible element of $\frak{Z}^{2}%
\left[ \mathbf{A}_{1}\right] ^{!}$ and $\left( \mathcal{A}_{\psi }\right)
^{!}=\left( \mathcal{A}^{!}\right) _{\psi ^{!}}.\;\;\;\blacksquare $
\end{theorem}

\subsubsection{Compositions of quantum spaces and quasicomplexes}

Given $\mathcal{A},\mathcal{B}\in \mathrm{FGA}$, consider a couple of
related admissible cochains $\psi \in \frak{Z}^{2}\left[ \mathbf{A}_{1}%
\right] $ and $\varphi \in \frak{Z}^{2}\left[ \mathbf{B}_{1}\right] $. From
the discussion we made in \S \textbf{2.3.2}, they give rise to a 2-cocycle $%
\left( \psi ,\varphi \right) $ in the quasicomplex $\frak{C}^{\bullet }\left[
\mathbf{A}_{1}\right] \times \frak{C}^{\bullet }\left[ \mathbf{B}_{1}\right] 
$ which, through the monic $\frak{j}:\frak{C}^{\bullet }\left[ \mathbf{A}_{1}%
\right] \times \frak{C}^{\bullet }\left[ \mathbf{B}_{1}\right]
\hookrightarrow \frak{C}^{\bullet }\left[ \mathbf{A}_{1}\otimes \mathbf{B}%
_{1}\right] $, can be identified with a map (recall Eq. $\left( \ref{moni}%
\right) $) 
\begin{equation*}
\frak{j}\left( \psi ,\varphi \right) =\bigoplus\nolimits_{r,s\in \mathbb{N}%
_{0}}\psi _{r,s}\otimes \varphi _{r,s}\in \frak{Z}^{2}\left[ \mathbf{A}%
_{1}\otimes \mathbf{B}_{1}\right] .
\end{equation*}
Also, they can be seen as a restriction of $\psi \otimes \varphi $. In these
terms, admissibility condition for $\psi $ and $\varphi $ lead us
immediately to $\mathcal{A}\circ \mathcal{B}$-admissibility of $\frak{j}%
\left( \psi ,\varphi \right) \subset \psi \otimes \varphi $. If $\psi
=\partial \theta $ and $\varphi =\partial \chi $, then the primitive of $%
\frak{j}\left( \psi ,\varphi \right) $ is $\frak{j}\left( \theta ,\chi
\right) =\bigoplus_{r\in \mathbb{N}_{0}}\theta _{r}\otimes \chi _{r}$, and
consequently $\frak{j}\left( \theta ,\chi \right) \subset \theta \otimes
\chi $. From that, it follows the ideal related to $\left( \mathcal{A}\circ 
\mathcal{B}\right) _{\frak{j}\left( \psi ,\varphi \right) }$ is precisely
the one related to $\mathcal{A}_{\psi }\circ \mathcal{B}_{\varphi }$. So, we
have

\begin{theorem}
Given $\mathcal{A},\mathcal{B}\in \mathrm{FGA}$ and a couple of counital
admissible cochains $\psi \in \frak{Z}^{2}\left[ \mathbf{A}_{1}\right] $ and 
$\varphi \in \frak{Z}^{2}\left[ \mathbf{B}_{1}\right] $, then $\frak{j}%
\left( \psi ,\varphi \right) \in \frak{Z}^{2}\left[ \mathbf{A}_{1}\otimes 
\mathbf{B}_{1}\right] $ is $\mathcal{A}\circ \mathcal{B}$-admissible and $%
\mathcal{A}_{\psi }\circ \mathcal{B}_{\varphi }\backsimeq \left( \mathcal{A}%
\circ \mathcal{B}\right) _{\frak{j}\left( \psi ,\varphi \right) }$.\ \ \ \ $%
\blacksquare $
\end{theorem}

For other composition of quantum spaces, it is not enough for $\psi $ and $%
\varphi $ to be admissible. For instance, the ideal related to $\mathcal{A}%
\odot \mathcal{B}$, $\bigoplus_{n\in \mathbb{N}_{0}}\mathbf{I}_{n}\otimes 
\mathbf{J}_{n}$, is not preserved by $\frak{j}\left( \psi ,\varphi \right) $%
. Nevertheless, if $\psi $ and $\varphi $ are 2nd admissible we can show,
using similar technics to the ones developed in the previous subsection,
that:

$\bullet $ $\frak{j}\left( \psi ,\varphi \right) $ is 2nd $\mathcal{A}\odot 
\mathcal{B}$-admissible and the related ideal of $\left( \mathcal{A}\odot 
\mathcal{B}\right) _{\frak{j}\left( \psi ,\varphi \right) }$ is 
\begin{equation}
\frak{j}\left( \theta ,\chi \right) ^{-1}\left( \bigoplus\nolimits_{n\in 
\mathbb{N}_{0}}\mathbf{I}_{n}\otimes \mathbf{J}_{n}\right)
=\bigoplus\nolimits_{n\in \mathbb{N}_{0}}\theta ^{-1}\left( \mathbf{I}%
_{n}\right) \otimes \chi ^{-1}\left( \mathbf{J}_{n}\right) ;  \label{X}
\end{equation}

$\bullet $ $\frak{j}\left( \psi ^{!},\varphi \right) $ is 2nd $\mathcal{A}%
\triangleright \mathcal{B}$-admissible and the related ideal of $\left( 
\mathcal{A}\triangleright \mathcal{B}\right) _{\frak{j}\left( \psi
^{!},\varphi \right) }$ is 
\begin{equation}
\begin{array}{l}
\frak{j}\left( \theta ^{!},\chi \right) ^{-1}\left( I\left[
\bigoplus\nolimits_{n\in \mathbb{N}_{0}}\mathbf{I}_{n}^{\perp }\otimes 
\mathbf{J}_{n}\right] \right) =I\left[ \bigoplus\nolimits_{n\in \mathbb{N}%
_{0}}\theta ^{\ast }\left( \mathbf{I}_{n}^{\perp }\right) \otimes \chi
^{-1}\left( \mathbf{J}_{n}\right) \right] \\ 
\\ 
=I\left[ \bigoplus\nolimits_{n\in \mathbb{N}_{0}}\theta ^{-1}\left( \mathbf{I%
}_{n}\right) ^{\perp }\otimes \chi ^{-1}\left( \mathbf{J}_{n}\right) \right]
;
\end{array}
\label{Y}
\end{equation}
and analogous results for the functors $\triangleleft $ and $\diamond $. The
same can be done if $\psi $ and $\varphi $ are admissible bicharacters,
since, for instance, cochains like $\psi ^{!}$ and $\frak{j}\left( \psi
,\varphi \right) $ will also have that property. Note that $\left( \ref{X}%
\right) $ and $\left( \ref{Y}\right) $ are precisely the ideals related to $%
\mathcal{A}_{\psi }\odot \mathcal{B}_{\varphi }$ and $\mathcal{A}_{\psi
}\triangleright \mathcal{B}_{\varphi }$, respectively. Summing up,

\begin{theorem}
Let $\psi \in \frak{Z}^{2}\left[ \mathbf{A}_{1}\right] $ and $\varphi \in 
\frak{Z}^{2}\left[ \mathbf{B}_{1}\right] $ be \textbf{2nd} $\mathcal{A}$ and 
$\mathcal{B}$-\textbf{admissible}, respectively. Then, 
\begin{equation*}
\begin{array}{ll}
\mathcal{A}_{\psi }\odot \mathcal{B}_{\varphi }=\left( \mathcal{A}\odot 
\mathcal{B}\right) _{\frak{j}\left( \psi ,\varphi \right) };\;\; & \mathcal{A%
}_{\psi }\triangleright \mathcal{B}_{\varphi }=\left( \mathcal{A}%
\triangleright \mathcal{B}\right) _{\frak{j}\left( \psi ^{!},\varphi \right)
}; \\ 
&  \\ 
\mathcal{A}_{\psi }\triangleleft \mathcal{B}_{\varphi }=\left( \mathcal{A}%
\triangleleft \mathcal{B}\right) _{\frak{j}\left( \psi ,\varphi ^{!}\right)
};\;\; & \mathcal{A}_{\psi }\diamond \mathcal{B}_{\varphi }=\left( \mathcal{A%
}\diamond \mathcal{B}\right) _{\frak{j}\left( \psi ,\varphi \right) ^{!}}.
\end{array}
\end{equation*}
And for quadratic and $m$-th quantum spaces, $\mathcal{A}_{\psi }\bullet 
\mathcal{B}_{\varphi }=\left( \mathcal{A}\bullet \mathcal{B}\right) _{\frak{j%
}\left( \psi ,\varphi \right) }$.\ \ \ \ $\blacksquare $
\end{theorem}

Among other things, results above imply the existence of functorial
injections 
\begin{equation*}
\mathsf{C}^{\bullet }\left[ \mathcal{A}\right] \times \mathsf{C}^{\bullet }%
\left[ \mathcal{B}\right] \hookrightarrow \mathsf{C}^{\bullet }\left[ 
\mathcal{A}\bigcirc \mathcal{B}\right] ,\;\;\bigcirc =\circ ,\bullet ,\odot
,\triangleleft ,\triangleright \;and\;\diamond .
\end{equation*}

\subsubsection{Twist transformations and the coHom objects}

As an application of formulae above, given $\mathcal{A}$ and $\mathcal{B}$
in $\mathrm{CA}$ and a pair of 2nd admissible cocycles $\psi \in \frak{Z}^{2}%
\left[ \mathbf{A}_{1}\right] $ and $\varphi \in \frak{Z}^{2}\left[ \mathbf{B}%
_{1}\right] $, $\underline{hom}\left[ \mathcal{B}_{\varphi },\mathcal{A}%
_{\psi }\right] =\underline{hom}\left[ \mathcal{B},\mathcal{A}\right] _{%
\frak{j}\left( \psi ^{!},\varphi \right) }$. Then, for twisted coHom objects
we have 
\begin{equation*}
\underline{hom}^{\Upsilon }\left[ \mathcal{B}_{\varphi },\mathcal{A}_{\psi }%
\right] =\underline{hom}\left[ \mathcal{B}_{\varphi },\mathcal{A}_{\psi }%
\right] _{\varsigma }=\underline{hom}\left[ \mathcal{B},\mathcal{A}\right] _{%
\frak{j}\left( \psi ^{!},\varphi \right) \star \varsigma }.
\end{equation*}
On the other hand, from $\underline{end}\left[ \mathcal{A}^{n|m}\right] $
and $\psi _{\hslash }$ (see \S \textbf{3.1.1}) we can compute the bialgebra $%
\underline{end}\left[ \mathcal{A}_{\hslash }^{n|m}\right] $ by making a
twist transformation $\frak{j}\left( \psi _{\hslash }^{!},\psi _{\hslash
}\right) $ (recall $\psi _{\hslash }$ is a bicharacter). In general, as a
direct consequence of \textbf{Theor. 16} for $\triangleright $:

\begin{proposition}
If $\mathcal{A}\backsim \mathcal{C}$ and $\mathcal{B}\backsim \mathcal{D}$
through 2nd admissible cocycles, then $\underline{hom}\left[ \mathcal{B},%
\mathcal{A}\right] \backsim \underline{hom}\left[ \mathcal{D},\mathcal{C}%
\right] $.\ \ \ $\blacksquare $
\end{proposition}

Now, we shall study twisted coHom objects under the perspective of twist
transformation. From \S \textbf{3.1.4}, we know that there exist admissible
bicharacters (see Eqs. $\left( \ref{sed}\right) $ and $\left( \ref{sedd}%
\right) $), \emph{ipso facto} 2nd admissible cochains, $\varsigma _{\mathcal{%
A}}\in \frak{Z}^{2}\left[ \mathbf{A}_{1}\right] $, $\varsigma _{\mathcal{B}%
}\in \frak{Z}^{2}\left[ \mathbf{B}_{1}\right] $ and $\varsigma _{\underline{%
hom}\left[ \mathcal{B},\mathcal{A}\right] }\in \frak{Z}^{2}\left[ \mathbf{B}%
_{1}^{\ast }\otimes \mathbf{A}_{1}\right] $ such that 
\begin{equation*}
\underline{hom}^{\Upsilon }\left[ \mathcal{B},\mathcal{A}\right] \backsimeq 
\mathcal{B}_{\varsigma }\triangleright \mathcal{A}_{\varsigma }=\underline{%
hom}\left[ \mathcal{B}_{\varsigma },\mathcal{A}_{\varsigma }\right]
\backsimeq \underline{hom}\left[ \mathcal{B},\mathcal{A}\right] _{\varsigma
}.
\end{equation*}
By direct calculations, it can be seen that $\varsigma _{\underline{hom}%
\left[ \mathcal{B},\mathcal{A}\right] }=\frak{j}\left( \varsigma _{\mathcal{B%
}}^{!},\varsigma _{\mathcal{A}}\right) $. On the other hand, the product $%
\underline{hom}^{\Upsilon }\left[ \mathcal{B},\mathcal{A}\right] \circ
_{\tau }\mathcal{B}$ appearing in twisted coevaluation map $\mathcal{A}%
\rightarrow \underline{hom}^{\Upsilon }\left[ \mathcal{B},\mathcal{A}\right]
\circ _{\tau }\mathcal{B}$, can be regarded as a twisting of the quantum
space $\underline{hom}^{\Upsilon }\left[ \mathcal{B},\mathcal{A}\right]
\circ \mathcal{B}$ by a cochain $\omega \in \frak{Z}^{2}\left[ \mathbf{B}%
_{1}^{\ast }\otimes \mathbf{A}_{1}\otimes \mathbf{B}_{1}\right] $. The
latter is an admissible anti-bicharacter defined by $\omega _{1,1}$ as
(recall Eq. $\left( \ref{tau}\right) $) 
\begin{equation*}
\omega \left( z_{i}^{j}\otimes z_{r}^{s}\otimes b_{k}\otimes b_{l}\right)
=\phi _{i}^{a}\phi _{r}^{b}\;\left( z_{a}^{c}\otimes z_{b}^{d}\right) \;\rho
_{c}^{j}\rho _{d}^{s}\otimes b_{k}\otimes b_{l},
\end{equation*}
writing $z_{i}^{j}=b^{j}\otimes a_{i}$ and making usual identifications.
Using $\left( \ref{abica}\right) $ we find that $\omega =\frak{j}\left( 
\frak{i}\varsigma _{\mathcal{B}}^{!},\frak{i}\varsigma _{\mathcal{A}},%
\mathbb{I}^{\otimes 2}\right) $. Since for all 2-cocycles $\psi $ and $%
\varphi $, 
\begin{equation*}
\frak{j}\left( \frak{i}\psi ,\frak{i}\varphi \right) =\frak{j}\partial
\left( \theta ^{-1},\chi ^{-1}\right) =\partial \left( \frak{j}\left( \theta
,\chi \right) ^{-1}\right) =\frak{ij}\left( \psi ,\varphi \right) ,
\end{equation*}
then $\omega =\frak{ij}\left( \varsigma _{\mathcal{B}}^{!},\varsigma _{%
\mathcal{A}},\mathbb{I}^{\otimes 2}\right) $, and consequently 
\begin{eqnarray*}
\underline{hom}^{\Upsilon }\left[ \mathcal{B},\mathcal{A}\right] \circ
_{\tau }\mathcal{B} &=&\left( \underline{hom}\left[ \mathcal{B},\mathcal{A}%
\right] _{\frak{j}\left( \varsigma _{\mathcal{B}}^{!},\varsigma _{\mathcal{A}%
}\right) }\circ \mathcal{B}\right) _{\omega }=\left( \left( \underline{hom}%
\left[ \mathcal{B},\mathcal{A}\right] \circ \mathcal{B}\right) _{\frak{i}%
\omega }\right) _{\omega } \\
&& \\
&=&\left( \underline{hom}\left[ \mathcal{B},\mathcal{A}\right] \circ 
\mathcal{B}\right) _{\frak{i}\omega \star \omega }=\underline{hom}\left[ 
\mathcal{B},\mathcal{A}\right] \circ \mathcal{B}.
\end{eqnarray*}
That is to say, the initiality of $\underline{hom}^{\Upsilon }\left[ 
\mathcal{B},\mathcal{A}\right] $ follows from that of $\underline{hom}\left[ 
\mathcal{B},\mathcal{A}\right] $, because the related coevaluation maps are
given by a same arrow. This way to see the construction of objects $%
\underline{hom}^{\Upsilon }\left[ \mathcal{B},\mathcal{A}\right] $ enable us
to generalize it to more general twist transformations. In fact, instead of
maps $\mathcal{A}\rightarrow \mathcal{H}\circ _{\tau }\mathcal{B}$ and
twisting $\tau _{\mathcal{A},\mathcal{B}}$, we can study a class of arrows $%
\mathcal{A}\rightarrow \left( \mathcal{H}\circ \mathcal{B}\right) _{\omega }$%
, where $\omega $ defines an element of $\frak{Z}^{2}\left[ \mathbf{H}%
_{1}\otimes \mathbf{B}_{1}\right] $. In terms of them, for each pair $%
\mathcal{A},\mathcal{B}\in \mathrm{CA}$, a full subcategory of $\left( \frak{%
H}\left( \mathcal{A}\right) \downarrow \frak{H}\left( \mathrm{CA}\circ 
\mathcal{B}\right) \right) $, being $\frak{H}:\mathrm{CA}\hookrightarrow 
\mathrm{GrVct}$ the forgetful functor to graded vector spaces, can be
defined as in \cite{gm}. If such a category, namely $\Omega ^{\mathcal{A},%
\mathcal{B}}$, is related to a cochain of the form $\omega _{\mathcal{A},%
\mathcal{B}}=\frak{ij}\left( \varsigma _{\mathcal{B}}^{!},\varsigma _{%
\mathcal{A}},\mathbb{I}^{\otimes 2}\right) \in \frak{Z}^{2}\left[ \mathbf{B}%
_{1}^{\ast }\otimes \mathbf{A}_{1}\otimes \mathbf{B}_{1}\right] $ (in a
similar way that $\Upsilon ^{\mathcal{A},\mathcal{B}}$ is related to $\tau _{%
\mathcal{A},\mathcal{B}}$),\footnote{%
To define $\omega _{\mathcal{A},\mathcal{B}}$ in a diagram $\left\langle
\varphi ,\mathcal{H}\right\rangle $, \textbf{Prop. 4} must be taken into
account for the inclusions $\mathbf{H}_{1}\otimes \mathbf{B}_{1}\subset 
\mathbf{B}_{1}^{\ast }\otimes \mathbf{A}_{1}\otimes \mathbf{B}_{1}$.} where $%
\varsigma _{\mathcal{X}}$ denotes a 2nd $\mathcal{X}$-admissible 2-cocycle,
then it will have initial object 
\begin{equation*}
\underline{hom}^{\Omega }\left[ \mathcal{B},\mathcal{A}\right] =\underline{%
hom}\left[ \mathcal{B},\mathcal{A}\right] _{\frak{j}\left( \varsigma _{%
\mathcal{B}}^{!},\varsigma _{\mathcal{A}}\right) }=\mathcal{B}_{\varsigma
}\triangleright \mathcal{A}_{\varsigma }.
\end{equation*}
Also, the disjoint union $\Omega ^{\cdot }=\bigvee_{\mathcal{A},\mathcal{B}%
}\Omega ^{\mathcal{A},\mathcal{B}}$ has a semigroupoid structure, and there
exists a related embedding $\Omega ^{\cdot }\hookrightarrow \mathrm{CA}$
that preserves the corresponding (partial) products. Accordingly, the
assignment $\left( \mathcal{B},\mathcal{A}\right) \mapsto \underline{hom}%
^{\Omega }\left[ \mathcal{B},\mathcal{A}\right] $ defines an $\mathrm{CA}$%
-cobased and (by duality) a $\mathsf{QLS}$-based category for each
collection of cocycles $\left\{ \omega _{\mathcal{A},\mathcal{B}}\right\} $.
These results are developed in \cite{gg}.

Finally, let us say the bialgebras $\underline{end}^{\Omega }\left[ \mathcal{%
A}\right] $ are twisting of $\underline{end}\left[ \mathcal{A}\right] $ also
as bialgebras (i.e. in the Drinfeld's sense). In fact, the isomorphism given
by Eq. $\left( \ref{cgad}\right) $ ensure the twisting $\frak{j}\left(
\varsigma ^{!},\varsigma \right) $ on each bialgebra $\underline{end}\left[ 
\mathcal{A}\right] $ is in the image of the map $\digamma $, i.e. there
exists $\chi \in \mathsf{G}^{\bullet }\left[ \mathbf{A}\triangleright 
\mathbf{A}\right] $ such that $\digamma _{\chi }=\frak{j}\left( \varsigma
^{!},\varsigma \right) $. Explicitly, 
\begin{equation*}
\chi _{\mathcal{A}}\left( z_{n_{1}}^{i_{1}}...z_{n_{r}}^{i_{r}}\otimes
z_{m_{1}}^{k_{1}}...z_{m_{s}}^{k_{s}}\right) =\left( \varsigma _{\mathcal{A}%
}\right) _{n_{1}...n_{r},m_{1}...m_{s}}^{i_{1}...i_{r},k_{1}...k_{s}}.
\end{equation*}
The fact that $\varsigma _{\mathcal{A}}$ is 2nd admissible is crucial for $%
\chi _{\mathcal{A}}$ to be well defined on $\left( \mathbf{A}\triangleright 
\mathbf{A}\right) ^{\otimes 2}$.

\section*{Conclusions}

Motivated by the idea of twisted internal coHom objects, we have defined a
non formal algebra deformation process over the category of quantum spaces.
Such deformations, or twisting, are controlled by a cosimplicial
multiplicative quasicomplex structure $\mathsf{C}^{\bullet }$ in the
category $\mathrm{Grp}_{\ast }$ of groups and unit preserving maps. Their
counital 2-cocycles are the elements implementing the mentioned
deformations. Also, given a quantum space $\mathcal{A}$, the cohomology
classes of 2-cocycles ,through the so called $\mathcal{A}$-admissible
cochains, defines the isomorphisms classes of twisting of $\mathcal{A}$.

We have shown that for each bialgebra in the category of quantum spaces, the
corresponding $\mathsf{C}^{\bullet }$ has the quasicomplex of bialgebra
twist transformations as a subobject. More precisely, if $\mathbf{V}$ is a
coalgebra, then there exists a monic $\mathsf{G}_{\mathcal{Z}}^{\bullet }%
\left[ \mathcal{V}\right] \hookrightarrow \mathsf{C}^{\bullet }\left[ 
\mathbf{V}\right] $ in the category of group quasicomplexes, where $\mathcal{%
V}$ is a bialgebra generated by $\mathbf{V}$.

The twist transformations define a gauge equivalence between quantum spaces.
For instance, every $n$-dimensional quantum plane $\mathcal{A}_{\hslash
}^{n|0}$ is gauge equivalent to the affine space $\mathcal{A}^{n|0}$. Also,
the twisted coHom objects are equivalent, in this sense, to the untwisted
ones.

The quasicomplexes $\mathsf{C}^{\bullet }\left[ \mathbf{V}\right] $, $%
\mathbf{V}\in \mathrm{Vct}_{f}$, generate a monoidal category with
involution. We have shown that, for a particular class of twist
transformations, the monoidal structure is compatible with products $\circ $%
, $\bullet $, $\odot $, $\triangleleft $, $\triangleright $ and $\diamond $
between quantum spaces, and the involution with the functor $!$, in the
sense that these functors (acting on objects) `commute' with the twisting
process.

\bigskip

Now, it is natural to ask what happen if, instead of groups $\mathsf{C}^{n}%
\left[ \mathbf{V}\right] $, we consider the algebras of degree cero $n$%
-homogeneous linear endomorphisms $\mathsf{E}^{n}\left[ \mathbf{V}\right] $.
Straightforwardly, the function $\left[ n+1\right] \mapsto \mathsf{E}^{n}%
\left[ \mathbf{V}\right] $ together with Eqs. $\left( \ref{cl}\right) $ and $%
\left( \ref{dege}\right) $ define a cosimplicial object in $\mathrm{Alg}$.
We just have to repeat the proof of \textbf{Theor. 2}. Moreover, maps $%
\partial _{\pm }$ given in Eq. $\left( \ref{paq}\right) $, supply $\mathsf{E}%
^{\bullet }$ with a structure of multiplicative cosimplicial parity
quasicomplex in $\mathrm{Alg}_{\ast }$ (the category of algebras and unit
preserving functions). From the algebra deformation process given at the
beginning of \S \textbf{2}, new twisting on quantum spaces can be
constructed within the cohomological framework defined by $\left( \mathsf{E}%
^{\bullet },\partial _{+},\partial _{-}\right) $. These new transformations
will be investigated elsewhere.

\section*{Acknowledgments}

The author thanks to CNEA, Argentina, for financial support, and thanks to
H. Montani for useful discussions.

\end{document}